\documentclass[12pt]{article}
\usepackage{amsfonts}
\usepackage{amsthm}
\input{epsf.sty}
\usepackage{latexsym,amsmath,amsfonts,amscd}
\usepackage{epsfig}
\usepackage{changebar}
\usepackage{pstricks}
\usepackage{pst-plot}
\usepackage{multirow}
\usepackage{subfigure}
\usepackage{placeins}
\usepackage{amssymb}
\usepackage{mathrsfs}
\usepackage[title]{appendix}
\usepackage{bm}
\usepackage{hyperref}

\usepackage{multirow,bigstrut}
\usepackage{enumitem}
\usepackage{booktabs}

\numberwithin{equation}{section}

\theoremstyle{plain}
\newtheorem{thm}{Theorem}[section]

\theoremstyle{definition}
\newtheorem{exam}[thm]{Example}

\newcommand{\brac}[1]{\left(#1\right)}
\newcommand{\abs}[1]{\left\vert#1\right\vert}

\def\half{\frac 1 2}

\newcommand{\bk}{\mathbf{k}}
\newcommand{\bx}{\mathbf{x}}

\newfont{\iams}{msbm9}

\newcommand{\commentbis}[1]{}
\newcommand{\be}{\begin{eqnarray}}
\newcommand{\ee}{\end{eqnarray}}
\newcommand{\beno}{\begin{eqnarray*}}
\newcommand{\eeno}{\end{eqnarray*}}

\newcommand{\barr}[1]{\begin{array}{#1}}
\newcommand{\earr}{\end{array}}

\newcommand{\beq}{\begin{equation}}
\newcommand{\eeq}{\end{equation}}
\newcommand{\beqa}{\begin{eqnarray}}
\newcommand{\eeqa}{\end{eqnarray}}

\newcommand{\bv}{{\bf v}}

\newcommand{\bV}{{\bf V}}

\newcommand{\bzero}{\mathbf{0}}
\newcommand{\bone}{\mathbf{1}}
\newcommand{\bl}{\mathbf{l}}
\newcommand{\bi}{\mathbf{i}}

\newcommand{\bj}{\mathbf{j}}
\newcommand{\bW}{\mathbf{W}}

\DeclareMathOperator{\sech}{sech}

\begin{document}
\baselineskip=1.5pc

\vspace{.5in}

\begin{center}

{\large\bf An adaptive multiresolution ultra-weak discontinuous Galerkin method for nonlinear Schr\"{o}dinger equations}

\end{center}

\vspace{.1in}

\centerline{
Zhanjing Tao  \footnote{School of Mathematics, Jilin University, Changchun, Jilin 130012, China. E-mail: zjtao@jlu.edu.cn}  \quad
Juntao Huang \footnote{Department of Mathematics,
Michigan State University, East Lansing, MI 48824, USA.
E-mail: huangj75@msu.edu. Corresponding author} \quad
Yuan Liu\footnote{Department of Mathematics, Statistics and Physics,
Wichita State University, Wichita, KS 67260, USA.
E-mail: liu@math.wichita.edu.
Research supported in part by a grant from the Simons Foundation (426993, Yuan Liu). }\quad
Wei Guo  \footnote{Department of Mathematics and Statistics, Texas Tech University, Lubbock, TX, 70409. E-mail:
weimath.guo@ttu.edu. Research is supported by NSF grant DMS-1830838} \quad
Yingda Cheng  \footnote{Department of Mathematics, Department of Computational Mathematics, Science and Engineering, Michigan State University, East Lansing, MI 48824, USA. E-mail: ycheng@msu.edu. Research is supported by NSF grants DMS-1453661 and DMS-1720023}
}

\vspace{.4in}

\centerline{\bf Abstract}

%\vspace{.1in}
This paper develops a high order adaptive scheme for solving nonlinear Schr\"{o}dinger equations. The solutions to such equations often exhibit solitary wave and  local structures, which makes adaptivity essential in improving the simulation efficiency. Our scheme uses the   ultra-weak discontinuous Galerkin (DG)  formulation and belongs to the framework of
adaptive multiresolution schemes. 
 Various numerical experiments are presented to demonstrate the excellent capability of capturing the soliton waves and the blow-up phenomenon.
%  ultra-weak discontinuous Galerkin(DG) method to solve nonlinear Schr\"{o}dinger equations. The two key ingredients of the scheme include an ultra-weak DG formulation in an adaptive function space and two classes of multiwavelets for achieving multiresolution. %In particular, the orthonormal Alpert’s multiwavelets are used to express the DG solution in terms of a hierarchical structure, and the interpolatory multiwavelets are further incorporated to enhance computational efficiency in the presence of nonlinear source. 

\bigskip

\bigskip
%\vfill


%{\bf AMS(MOS) subject classification:} 65M99

\pagenumbering{arabic}

%section 1
\section{Introduction}\label{intro}

\setcounter{equation}{0}
\setcounter{figure}{0}
\setcounter{table}{0}

In this paper, we develop a class of adaptive multiresolution ultra-weak discontinuous Galerkin (DG) method to solve the nonlinear Schr\"{o}dinger (NLS) equations in $d$-dimensional space
\begin{equation}\label{eq:NLS}
    i u_t + \Delta u + f(\abs{u}^2)u = 0,
\end{equation}
where $u$ is a complex function, $f$ is a smooth nonlinear real function.
%and the coupled nonlinear Schr\"{o}dinger equations in one-dimensional space
%\begin{subequations}\label{eq:couple-NLS}
%\begin{align}
%    & i u_t + i \alpha u_x + u_{xx} + \beta u + \kappa v + f(\abs{u}^2, \abs{v}^2)u = 0, \\
%    & i v_t - i \alpha v_x + v_{xx} - \beta u + \kappa v + g(\abs{u}^2, \abs{v}^2)v = 0.
%\end{align}
%\end{subequations}
%where $u$ and $v$ are complex functions, $f$ and $g$ are smooth nonlinear real functions, and $\alpha$, $\beta$, $\kappa$ are real constants. 
The Schr\"{o}dinger equation is of fundamental importance in quantum mechanics, reaching out to many important applications describing the physical phenomena including nonlinear optics, semiconductor electronics, quantum fluids and plasma physics \cite{newell1985solitons, whitham2011linear, chiao1964self}. 
Numerical methods for solving the NLS equations have been investigated extensively, including finite difference \cite{sanz1986conerservative,chang1999difference,ismail2001numerical,sheng2001solving,taha1984analytical}, finite element \cite{griffiths1984numerical,karakashian1998space,xu2005schrodinger,chen2019ultra}, and spectral methods \cite{de2008exponential,pathria1990pseudo,sulem1984numerical}, to name a few. In this paper, we consider the DG method \cite{reed1973triangular,cockburn2000development,cockburn2001runge}, which is a class of finite element methods using piecewise polynomial spaces for the numerical solutions and the test functions. The last several decades have seen tremendous developments of DG methods in approximating partial differential equations (PDEs) in large part due to their distinguished advantages in handling geometry, boundary conditions and accommodating adaptivity. %For a detailed description of the DG methods as well as the implementation and applications, we refer readers to the review papers .
Various types of DG methods have been proposed to compute the NLS equations. In \cite{xu2005schrodinger}, an LDG method using alternating fluxes was developed with $L^2$ stability and the sub-optimal error estimates. An LDG method with various numerical fluxes was analyzed in \cite{liang2015fourth}. An analysis of the LDG method for the NLS equation with wave operator was carried out in \cite{guo2015energy}. The direct DG (DDG) method was applied to Schr\"{o}dinger equation in \cite{lu2015mass}, and the optimal accuracy was further established in \cite{liu2019accuracy}. In \cite{xiong2018uniform}, an hybridized DG (HDG) method was applied to a linear Schrödinger equation.  In this paper, we use the ultra-weak DG method \cite{cheng2008discontinuous}, which is a class of DG methods use repeated integration by parts for calculating higher order derivatives. The ultra-weak DG schemes include the DDG and interior penalty DG methods, and have been investigated in \cite{chen2019ultra,chen2020superconvergence} for convergence and superconvergence.

%, though our framework can be extended to incorporate other types of DG schemes.

The solutions to NLS equations present   solitary waves, blow-up and other localized structures. Therefore, benefits of adaptivity in simulations are self-evident \cite{sanz1986simple,chang1990multigrid,kormann2016time}. 
In this paper, we  consider the adaptive multiresolution approach \cite{gerhard2016adaptive,guo2017adaptive,huang2019adaptive}. By exploring the inherent mesh hierarchy and the associated nestedness of the polynomial approximation spaces, multiresolution analysis (MRA) \cite{mallat1999wavelet} is able to accelerate the computation and avoid the need for a posteriori error indicators. MRA is closely related to popular sparse grid methods \cite{bungartz2004sparse} for solving high dimensional problems.
It is also related to the adaptive mesh refinement (AMR) technique \cite{berger1989local,burstedde2011p4est}, which adjusts the computational grid adaptively to track small scale features of the underlying problems and improves computational efficiency.   
% 
% 
% Multiresolution analysis (MRA) \cite{mallat1999wavelet} explores mesh hierarchy, which induces nested polynomial approximation spaces to accelerate the computation and in the mean time avoid the need for a posteriori error indicators.  
% 
 As a continuation of our previous research for adaptive multiresolution (also called adaptive sparse grid) DG methods \cite{guo2017adaptive,huang2019adaptive,huang2020adaptive,guo2020adaptive}, this paper develops an adaptive multiresolution ultra-weak DG solver for NLS equations \eqref{eq:NLS} and the coupled NLS equations. First, the Alpert’s multiwavelets are employed as the DG bases in the weak formulation, and then the interpolatory multiwavelets are introduced for efficiently computing nonlinear source which has been successfully applied to nonlinear hyperbolic conservation laws \cite{huang2019adaptive} and Hamilton-Jacobi equations \cite{guo2020adaptive}. We refer the readers to \cite{guo2017adaptive,huang2019adaptive} for more details on the background of adaptive multiresolution DG methods. Numerical experiments verify the accuracy of the methods. In particular, the adaptive scheme is demonstrated to capture the moving solitons and also the blow-up phenomenon very well.

% The major numerical challenge in designing numerical methods for nonlinear Schr\"{o}dinger equation is to capture the solution structures at different scales with affordable computational cost. In real applications, the soliton waves and the singular solutions of the NLS equation usually attract a lot of interests. Here, we focus on the adaptive multiresolution DG method for solving this equation. 

The rest of the paper is organized as follows. In Section 2, we review Alpert’s multiwavelets. Section 3 describes the numerical schemes. Section 4 contains numerical examples. We make conclusions in Section 5.

% section 2
\section{Multiresolution analysis and multiwavelets}
\label{sec:mra}

In this section, we briefly review the fundamentals of MRA of DG approximation spaces and the associated multiwavelets. Two classes of multiwavelets, namely  the $L^2$ orthonormal Alpert's multiwavelets \cite{alpert1993class} and the interpolatory multiwavelets \cite{tao2019collocation}, are used to construct our ultra-weak DG scheme. We also introduce a set of key notations used throughout the paper by following \cite{wang2016elliptic}. 

% \subsection{Alpert's multiwavelets}\label{subsec:alpert-basis}

 Alpert's multiwavelets \cite{alpert1993class} have been employed to develop a class of sparse grid DG methods for solving high dimensional PDEs \cite{wang2016elliptic,guo2016sparse}.
Considering a unit sized interval $\Omega=[0,1]$ for simplicity,  we define a set of nested grids $\Omega_0,\,\Omega_1,\ldots$, for which the $n$-th level grid $\Omega_n$ consists of $2^n$ uniform cells
\begin{equation*}
I_{n}^j=(2^{-n}j, 2^{-n}(j+1)], \quad j=0, \ldots, 2^n-1,
\end{equation*}
 Denote $I_{-1}=[0,1].$
The piecewise polynomial space of degree at most $k\ge1$ on grid $\Omega_n$ for $n\ge 0$ is denoted by
\begin{equation}\label{eq:DG-space-Vn}
V_n^k:=\{v: v \in P^k(I_{n}^j),\, \forall \,j=0, \ldots, 2^n-1\}.
\end{equation}
Observing the nested structure
$$V_0^k \subset V_1^k \subset V_2^k \subset V_3^k \subset  \cdots,$$
we can define the multiwavelet subspace $W_n^k$, $n=1, 2, \ldots $ as the orthogonal complement of $V_{n-1}^k$ in $V_{n}^k$ with respect to the $L^2$ inner product on $[0,1]$, i.e.,
\begin{equation*}
V_{n-1}^k \oplus W_n^k=V_{n}^k, \quad W_n^k \perp V_{n-1}^k.
\end{equation*}
By letting $W_0^k:=V_0^k$, we obtain a hierarchical decomposition $V_n^k=\bigoplus_{0 \leq l \leq n} W_l^k$, i.e., MRA of space $V_n^k$.
A set of orthonormal basis can be defined on $W_l^k$ as follows. When $l=0$, the basis $v^0_{i,0}(x)$, $ i=0,\ldots,k$ are the normalized shifted Legendre polynomials in $[0,1]$. When $l>0$, the Alpert's orthonormal multiwavelets \cite{alpert1993class} are employed as the bases and denoted by 
$$v^j_{i,l}(x),\quad i=0,\ldots,k,\quad j=0,\ldots,2^{l-1}-1.$$

We then follow a tensor-product approach to construct the hierarchical finite element space in multi-dimensional space.  Denote $\bl=(l_1,\cdots,l_d)\in\mathbb{N}_0^d$ as the mesh level in a multivariate sense, where $\mathbb{N}_0$  denotes the set of nonnegative integers, we can define the tensor-product mesh grid $\Omega_\bl=\Omega_{l_1}\otimes\cdots\otimes\Omega_{l_d}$ and the corresponding mesh size $h_\bl=(h_{l_1},\cdots,h_{l_d}).$ Based on the grid $\Omega_\bl$, we denote  $I_\bl^\bj=\{\bx:x_m\in(h_mj_m,h_m(j_{m}+1)),m=1,\cdots,d\}$ as an elementary cell, and 
$$\bV_\bl^k:=\{\bv:  \bv \in Q^k(I^{\bj}_{\bl}), \,\,  \bzero \leq \bj  \leq 2^{\bl}-\bone \}= V_{l_1,x_1}^k\times\cdots\times  V_{l_d,x_d}^k$$
as the tensor-product piecewise polynomial space, where $Q^k(I^{\bj}_{\bl})$ represents the collection of polynomials of degree up to $k$ in each dimension on cell $I^{\bj}_{\bl}$. 
If we use equal mesh refinement of size $h_N=2^{-N}$ in each coordinate direction, the  grid and space will be denoted by $\Omega_N$ and $\bV_N^k$, respectively.  
Based on a tensor-product construction, the multi-dimensional increment space can be  defined as
$$\bW_\bl^k=W_{l_1,x_1}^k\times\cdots\times  W_{l_d,x_d}^k.$$
The basis functions in multi-dimensions are defined as
\begin{equation}\label{eq:multidim-basis}
v^\bj_{\bi,\bl}(\bx) := \prod_{m=1}^d v^{j_m}_{i_m,l_m}(x_m),
\end{equation}
for $\bl \in \mathbb{N}_0^d$, $\bj \in B_\bl := \{\bj\in\mathbb{N}_0^d: \,\mathbf{0}\leq\bj\leq\max(2^{\bl-\mathbf{1}}-\mathbf{1},\mathbf{0}) \}$ and $\mathbf{1}\leq\bi\leq \bk+\mathbf{1}$. 
% {\color{red} Note that the basis construction is different from \cite{gerhard2016adaptive}: $v^\bj_{\bi,\bl}(\bx)$ may have an anisotropic support, while the basis introduced in \cite{gerhard2016adaptive} does not.}

 %The orthonormality of the bases can be easily verified.

Introducing the standard norms for the multi-index $$
|\bl|_1:=\sum_{m=1}^d l_m, \qquad   |\bl|_\infty:=\max_{1\leq m \leq d} l_m,
$$
together with  the same component-wise arithmetic operations and relations   as defined in \cite{wang2016elliptic},  we reach the decomposition
\begin{equation}\label{eq:hiere_tp}
\bV_N^k=\bigoplus_{\substack{ |\bl|_\infty \leq N\\\bl \in \mathbb{N}_0^d}} \bW_\bl^k.
\end{equation}
%On the other hand, a standard choice of sparse grid   space 
Further, by a standard truncation of $\bV_N^k$ \cite{wang2016elliptic, guo2016sparse}, we obtain the sparse grid   space 
\begin{equation}
\label{eq:hiere_sg}
\hat{\bV}_N^k=\bigoplus_{\substack{ |\bl|_1 \leq N\\\bl \in \mathbb{N}_0^d}}\bW_\bl^k \subset \bV_N^k.
\end{equation}
We skip the details about the property of the space, but refer the readers to \cite{wang2016elliptic, guo2016sparse}. In Section \ref{sec:DG-adaptive}, we will describe   the adaptive scheme which adapts a subspace of $\bV_N^k$ according to the numerical solution, hence offering more flexibility and efficiency.

Alpert's multiwavelets described above are associated with the $L^2$ projection operator. For nonlinear source terms, we use the interpolatory multiwavelets based on Lagrange interpolations introduced in \cite{tao2019collocation}. For the details, we refer readers to \cite{tao2019collocation,huang2019adaptive}.

\section{Adaptive multiresolution DG scheme}
\label{sec:DG-adaptive}

In this section, we present the adaptive multiresolution ultra-weak DG scheme for solving the NLS equation \eqref{eq:NLS}. We consider periodic boundary conditions for simplicity, while the method can be adapted to other non-periodic boundary conditions. 

For illustrative purposes, we first introduce some basis notation about jumps and averages for piecewise functions defined on a grid $\Omega_N$. Denote by $\Gamma$ the union of the boundaries for all the elements in the partition $\Omega_N$. The jump  and average of $q\in L^2(\Gamma)$ and $\textbf{q}\in [L^2(\Gamma)]^d$ are defined as follows. Suppose $e$ is an edge shared by elements $T^+$ and $T^-$, we define the unit normal vectors $\textbf{n}^+$ and $\textbf{n}^-$ on $e$ pointing exterior to $T^+$ and $T^-$, and then 
\begin{align}
[q]=q^- \textbf{n}^- + q^+ \textbf{n}^+, \qquad & \{ q\} = \frac{1}{2} (q^-+q^+),  \nonumber  \\
[\textbf{q}] = \textbf{q}^- \cdot \textbf{n}^-  + \textbf{q}^+ \cdot \textbf{n}^+ ,  \qquad & \{ \textbf{q}\} = \frac{1}{2}(\textbf{q}^- + \textbf{q}^+). \nonumber
\end{align}
%where $\textbf{n}$ is the unit normal.

For any subspace $\bV$ of $\bV_N^k$, define the corresponding complex-valued finite element space
\begin{equation}
	\mathbb{V}:=\{ v=v_1 + i v_2: v_1, v_2 \in \bV \}
\end{equation}

The semi-discrete ultra-weak DG scheme \cite{chen2019ultra} for \eqref{eq:NLS} is defined as follows: we are looking for $u_h \in \mathbb{V}$ such that for any test function $\phi_h \in \mathbb{V}$,
\begin{align}\label{eq:semi-uwdg}
i \int_{\Omega} (u_h)_t \phi_h d\bx  + \int_{\Omega} u_h \nabla^2 \phi_h d \bx & - \sum_{e\in \Gamma} \int_e \hat{u}_h [\nabla \phi_h] ds  +\sum_{e \in \Gamma} \int_e \widetilde{\nabla u_h} \cdot [\phi_h] ds \\  \nonumber & +\int_{\Omega} f(|u_h|^2)u_h \phi_h d\bx = 0.
\end{align}
We take the following numerical fluxes 
\begin{align}\label{eq:num-flux}
 \widetilde{\nabla u_h} = \{ \nabla u_h \} + \alpha_1[\nabla u_h]{\bf e} + \beta_1[u_h], \quad \hat{u}_h = \{u_h\}  + \alpha_2[u_h]\cdot{\bf e} + \beta_2[\nabla u_h].
\end{align}
Here $\alpha_1$, $\alpha_2$, $\beta_1$ and $\beta_2$ are prescribed complex numbers which may depend on the mesh size $h$ and ${\bf e}=(1,\dots,1)\in\mathbb{R}^d$. In this work, we numerically test two types of numerical fluxes. The first one is the alternating flux corresponding to $\alpha_1=\half$, $\alpha_2=-\half$ and $\beta_1=\beta_2=0$. The second one is a dissipative numerical flux \cite{chen2019ultra} which has the parameters $\alpha_1 = \half$, $\alpha_2 = -\half$, $\beta_1 = 1-i$, $\beta_2 = 1+i$.  

In order to efficiently calculate the nonlinear term $\int_{\Omega} f(|u_h|^2)u_h \phi_h d\bx$ in \eqref{eq:semi-uwdg}, the multiresolution Lagrange interpolation is applied \cite{tao2019collocation,huang2019adaptive}, i.e., we modified the weak formulation of ultra-weak DG as follows. We are looking for $u_h \in \mathbb{V}$ such that for any test function $\phi_h \in \mathbb{V}$,
\begin{align}\label{eq:semi-uwdg-interp}
i \int_{\Omega} (u_h)_t \phi_h d\bx  + \int_{\Omega} u_h \nabla^2 \phi_h d \bx & - \sum_{e\in \Gamma} \int_e \hat{u}_h [\nabla \phi_h] ds  +\sum_{e \in \Gamma} \int_e \widetilde{\nabla u_h} \cdot [\phi_h] ds \\  \nonumber & +\int_{\Omega} \mathcal{I}_h\brac{f(|u_h|^2)u_h} \phi_h d\bx = 0. 
\end{align}
To preserve the accuracy of the original DG scheme \eqref{eq:semi-uwdg}, it is required that Lagrange interpolation of the same order is applied in \eqref{eq:semi-uwdg-interp}. For the details, see the argument in \cite{DG4,huang2017quadrature,huang2019adaptive}. By applying the interpolation, the unidirectional principle and fast algorithm described in \cite{huang2019adaptive} can be employed to further improve efficiency. In numerical experiments, we also consider 
the coupled nonlinear Schr\"{o}dinger equations in one-dimensional space
\begin{subequations}\label{eq:couple-NLS}
\begin{align}
    & i u_t + i \alpha u_x + u_{xx} + \beta u + \kappa v + f(\abs{u}^2, \abs{v}^2)u = 0, \\
    & i v_t - i \alpha v_x + v_{xx} - \beta u + \kappa v + g(\abs{u}^2, \abs{v}^2)v = 0,
\end{align}
\end{subequations}
where $u$ and $v$ are complex functions, $f$ and $g$ are smooth nonlinear real functions, and $\alpha$, $\beta$, $\kappa$ are real constants. 
We use the same DG scheme for solving the coupled NLS equation \eqref{eq:couple-NLS} except that the first order derivatives $u_x$ and $v_x$ are treated by the standard DG scheme with upwind numerical fluxes. The details are omitted here for brevity.

% Alpert's multiwavelets described in Section \ref{subsec:alpert-basis} are associated with the $L^2$ projection operator. The interpolatory multiwavelets introduced in \cite{tao2019collocation} are constructed based on interpolation operators and also essential for efficient computation of integrals in the DG formulation, especially in high dimensions. In this work, only Lagrange interpolation is considered, while we note that Hermite interpolation can also be used but its implementation is more involved. The details are provided below.

%\begin{align}\label{udg4}
%i\int_{\Omega} u_t \phi dx & - i \alpha \int_{\Omega} u \phi_x dx +\int_{\Omega} u \phi_{xx} dx + \int_{\Omega} (\beta u  + \kappa v + f u) \phi dx \\ \nonumber
%& - i \alpha \sum^N_{j=1} \widetilde{u}_{j+\frac{1}{2}} [\phi]_{j+\frac{1}{2}} - \sum^N_{j=1} (\widetilde{u}_x)_{j+\frac{1}{2}} [\phi]_{j+\frac{1}{2}} + \sum^N_{j=1} \hat{u}_{j+\frac{1}{2}}[\phi_x]_{j+\frac{1}{2}} = 0 \\ \nonumber
%i\int_{\Omega} v_t \psi dx & + i \alpha \int_{\Omega} v \psi_x dx + \int_{\Omega} v \phi_{xx} dx + \int_{\Omega} (-\beta u + \kappa v + gv)\psi dx \\
%& + i \alpha \sum^N_{j=1} \widetilde{v}_{j+\frac{1}{2}} [\psi]_{j+\frac{1}{2}} - \sum^N_{j=1} (\widetilde{v}_x)_{j+\frac{1}{2}} [\psi]_{j+\frac{1}{2}} + \sum^N_{j=1} \hat{v}_{j+\frac{1}{2}}[\psi_x]_{j+\frac{1}{2}} = 0
%\end{align}

For time discretization, we employ the third order implicit-explicit (IMEX) Runge-Kutta (RK) scheme \cite{pareschi2005implicit} to advance the semi-discrete scheme \eqref{eq:semi-uwdg-interp}. Specifically, the second derivative term $u_{xx}$ is treated implicitly to avoid the severe CFL time constraint, while the nonlinear source $f(\abs{u}^2)u$ is treated explicitly for efficiency. The adaptive procedure follows the technique developed in \cite{guo2020adaptive,huang2019adaptive} to determine the space $\bV$ that evolves dynamically over time. The only difference is that the first-order Euler forward and Euler backward scheme is applied for the prediction procedure. The main idea is that in light of the distinguished property of multiwavelets, we keep track of multiwavelet coefficients, i.e. $L^2$ norms of $u_h$, as an error indicator for refining and coarsening, aiming to efficiently capture the solitons or singular solution of \eqref{eq:NLS}. We also remark that other types of time discretizations, e.g., exponential time differencing (ETD) or Krylov implicit integration factor (IIF) methods, could also be applied here. The efficiency of different types of time stepping remains to be investigated.

\section{Numerical examples}\label{sec:numer}

In this section, we perform numerical experiments to validate the performance of our scheme. We consider the NLS equation \eqref{eq:NLS} in 1D and 2D, and the coupled NLS equation \eqref{eq:couple-NLS} in 1D, with computational domain being $[0, 1]^d$ with $d = 1, 2$. We employ the third order IMEX RK  scheme in \cite{pareschi2005implicit}. The CFL number is taken to be 0.1, i.e., $\Delta t=0.1 \Delta x$, unless otherwise stated. All adaptive calculations are obtained by $k = 3$. $\textrm{DoF}=\textrm{dim}(\bV)$ refers to the number of Alperts' multiwavelets basis functions in the adaptive grids.

\subsection{Accuracy test for NLS equation}

\begin{exam}\label{ex:accuracy}
We start with the accuracy test for the NLS equation on the domain $[0,1]^d$:
\begin{equation}
  i u_t + \Delta u + \abs{u}^2 u + \abs{u}^4 u = 0,
\end{equation}
with periodic boundary conditions. The exact solution is taken to be
\begin{equation}
    u(x,t) = \exp(i(2\pi\sum_{i=1}^dx_i - \omega t)),
\end{equation}
with $\omega=4d \pi^2-2$.
\end{exam}

We first test the accuracy of sparse grid in 2D. The results with $k=1,2,3$ are presented in Table \ref{tb:accuracy-NLS-2D-sparse}. To match the accuracy of time, we take $\Delta t=0.1\Delta x^{4/3}$ for $k=3$. As expected, the convergence order in average is between $k$ and $k+1$.

\begin{table}[htbp]
  \centering
\caption{Example \ref{ex:accuracy}, accuracy test for NLS equation, $d=2$, sparse grid, $t=0.1$.}
\label{tb:accuracy-NLS-2D-sparse}
    \begin{tabular}{c|c|c|c|c|c}
    \hline
    \multirow{2}[4]{*}{} & \multirow{2}[4]{*}{} & \multicolumn{2}{c|}{Real part} & \multicolumn{2}{c}{Imaginary part} \bigstrut\\
\cline{3-6}          &   $N$    &   $L^2$-error    &    order   &   $L^2$-error   &  order \bigstrut\\
    \hline
    \multirow{5}{3em}{$k=1$} 
		&  5    &  2.82e-01     &   -    &  2.90e-01     &  - \bigstrut[t]\\
		&  6    &  1.28e-01     &   1.15    &  1.35e-01     & 1.11 \\
		&  7    &  1.90e-02     &   2.74    &  1.90e-02     & 2.83 \\
		&  8    &  5.37e-03     &   1.82    &  5.27e-03     & 1.85 \\
		&  9    &  1.13e-03     &   2.25    &  1.13e-03     & 2.22 \bigstrut[b]\\
    \hline
    \multirow{5}{3em}{$k=2$} 
		&  3    &  3.20e-02     &   -    &  4.33e-02     &  - \bigstrut[t]\\
		&  4    &  7.91e-03     &   2.02    &  1.43e-02     & 1.60 \\
		&  5    &  7.74e-04     &   3.35    &  7.77e-04     & 4.20 \\
		&  6    &  1.88e-04     &   2.04    &  2.66e-04     & 1.55 \\
		&  7    &  1.46e-05     &   3.68    &  1.47e-05     & 4.18 \bigstrut[b]\\
    \hline
    \multirow{5}{3em}{$k=3$}
		&  3    &  9.82e-03     &   -    &  2.67e-02     &  - \bigstrut[t]\\
		&  4    &  1.96e-04     &   5.64    &  2.29e-04     & 6.87 \\
		&  5    &  2.05e-05     &   3.26    &  1.46e-05     & 3.97 \\
		&  6    &  2.60e-06     &   2.98    &  9.40e-07     & 3.95 \\
		&  7    &  5.99e-08     &   5.44    &  5.89e-08     & 4.00 \bigstrut[b]\\
    \hline    
    \end{tabular}%
  \label{tab:addlabel}%
\end{table}%

We then test the accuracy of adaptive method in 2D in Table \ref{tb:accuracy-NLS-2D-adaptive}. We observe that it takes much less DoF with higher order polynomial degrees than lower order ones.

\begin{table}[!hbp]
	\centering
	\caption{Example \ref{ex:accuracy}, accuracy test for NLS equation, $d=2$. Adaptive. $t=0.1$.}
	\label{tb:accuracy-NLS-2D-adaptive}
	\begin{tabular}{c|c|c|c|c|c|c|c|c}
		\hline
		 \multirow{2}[4]{*}{} & \multirow{2}[4]{*}{$\epsilon$} & \multirow{2}[4]{*}{DoF} & \multicolumn{3}{c|}{Real part of $u$} & \multicolumn{3}{c}{Imaginary part of $u$} \\
		\cline{4-9} 
		&  &  & L$^2$-error & $R_{\textrm{DoF}}$ & $R_{\epsilon}$ & L$^2$-error & $R_{\textrm{DoF}}$ & $R_{\epsilon}$\\
		\hline
		\multirow{4}{3em}{$k=1$}
	   & 1e-01 & 192 & 9.33e-01 & - & - & 9.59e-01 & - & - \\
	   & 1e-02 & 960 & 9.39e-02 & 1.43 & 1.00 & 9.25e-02 & 1.45 & 1.02 \\
	   & 1e-03 & 1792 & 2.33e-02 & 2.23 & 0.61 & 2.29e-02 & 2.23 & 0.61 \\
	   & 1e-04 & 7168 & 4.80e-03 & 1.14 & 0.69 & 4.78e-03 & 1.13 & 0.68 \\	
		\hline
		\multirow{4}{3em}{$k=2$}
		& 1e-01 & 108 & 7.93e-02 & - & - & 8.68e-02 & - & - \\
		& 1e-02 & 432 & 1.01e-02 & 1.49 & 0.89 & 9.60e-03 & 1.59 & 0.96 \\
		& 1e-03 & 720 & 1.10e-03 & 4.35 & 0.96 & 1.11e-03 & 4.23 & 0.94 \\
		& 1e-04 & 1800 & 2.06e-04 & 1.82 & 0.73 & 2.99e-04 & 1.43 & 0.57 \\
		\hline
		\multirow{4}{3em}{$k=3$}
		& 1e-02 & 320 & 1.43e-02 & - & - & 2.86e-02 & - & - \\
		& 1e-03 & 512 & 2.81e-03 & 3.46 & 0.71 & 2.81e-03 & 4.93 & 1.01 \\
		& 1e-04 & 896 & 3.08e-04 & 3.95 & 0.96 & 3.07e-04 & 3.96 & 0.96 \\
		& 1e-05 & 1984 & 3.63e-05 & 2.69 & 0.93 & 3.63e-05 & 2.69 & 0.93 \\
		\hline
	\end{tabular}
\end{table}

\begin{figure}
	\centering
	\subfigure[$k=1$]{
		\begin{minipage}[b]{0.46\textwidth}
			\includegraphics[width=1\textwidth]{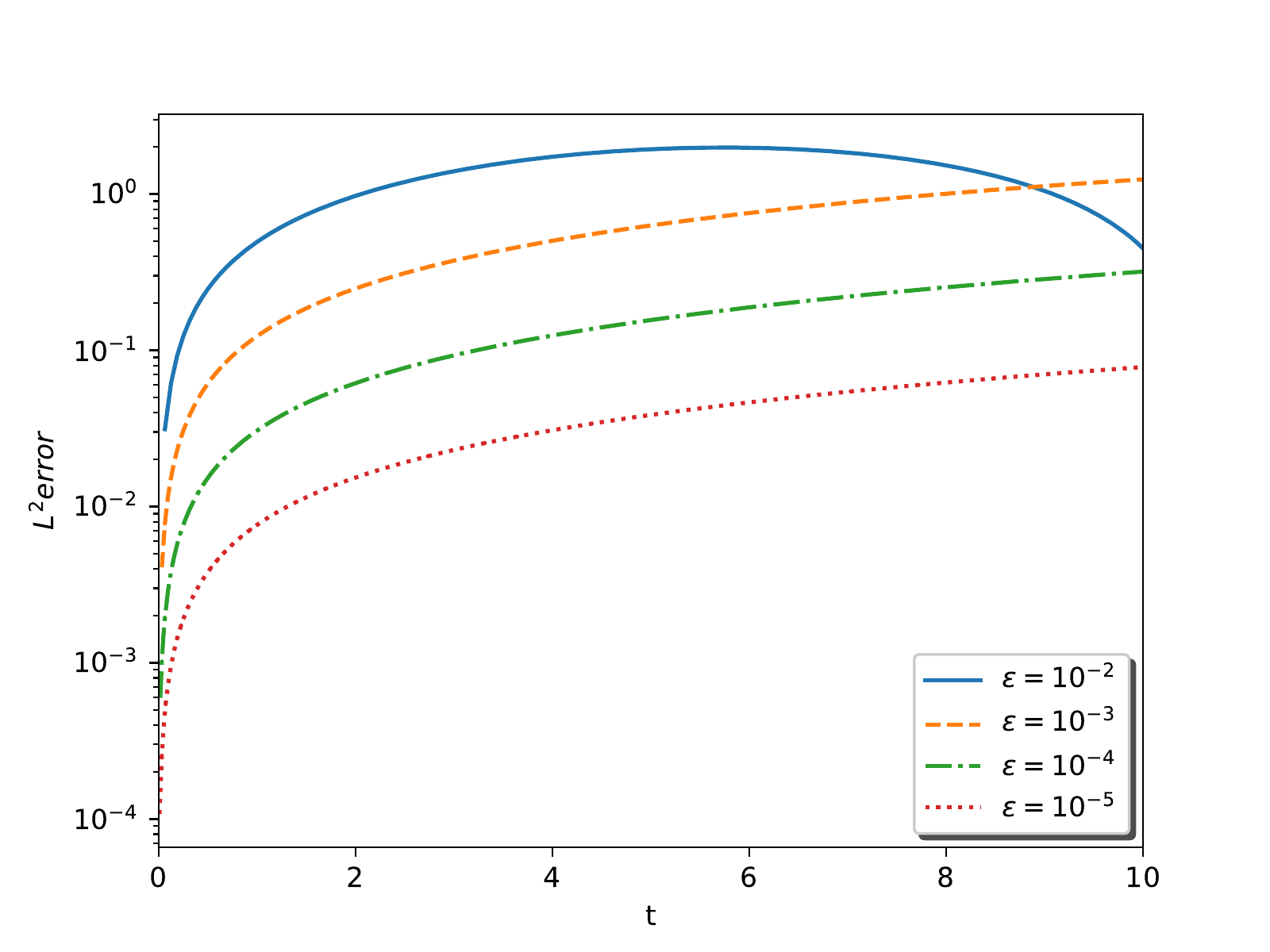}
		\end{minipage}
	}
	\subfigure[$k=1$]{
		\begin{minipage}[b]{0.46\textwidth}    
			\includegraphics[width=1\textwidth]{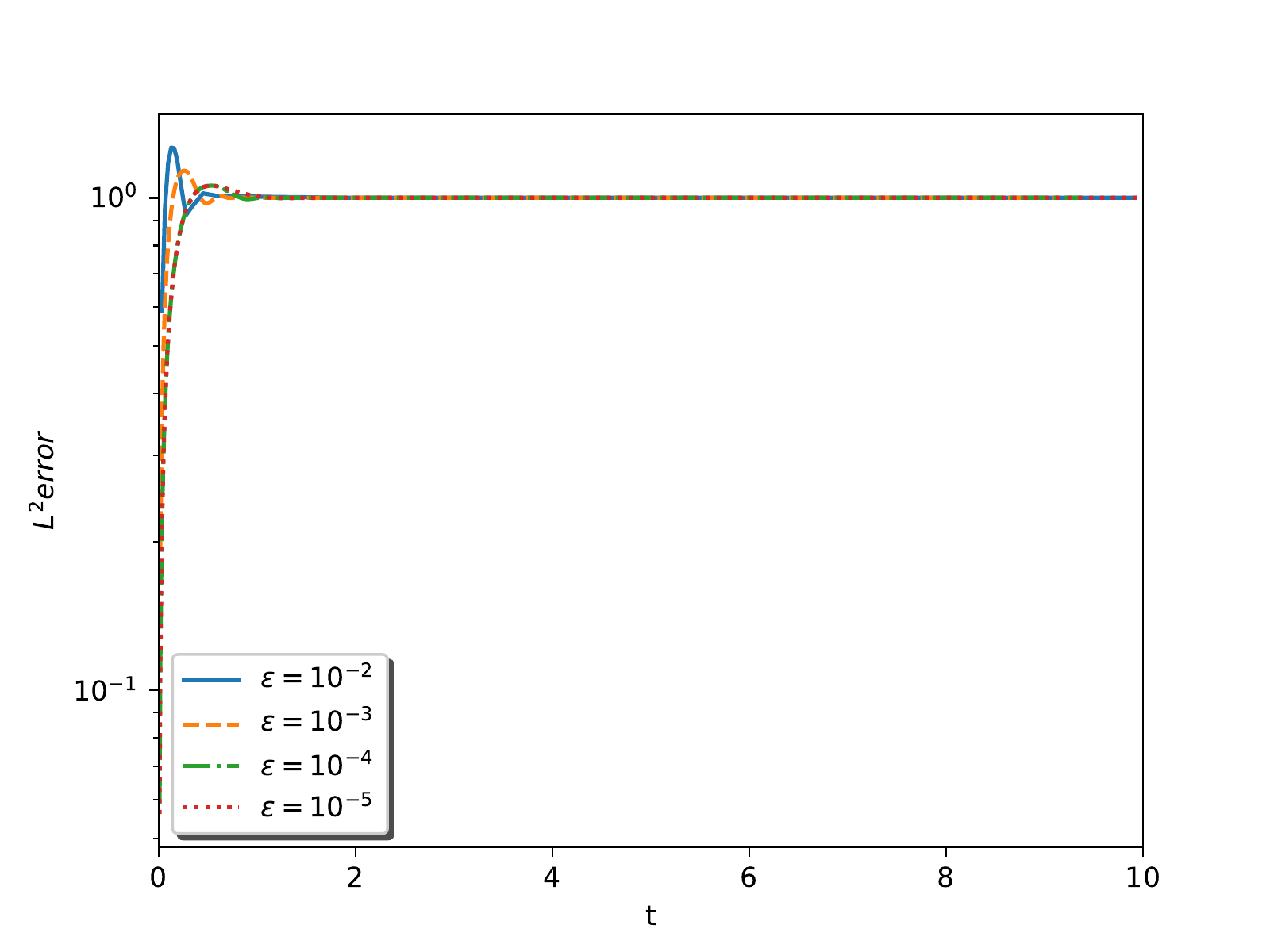}
		\end{minipage}
	}
	\bigskip
	\subfigure[$k=2$]{
		\begin{minipage}[b]{0.46\textwidth}
			\includegraphics[width=1\textwidth]{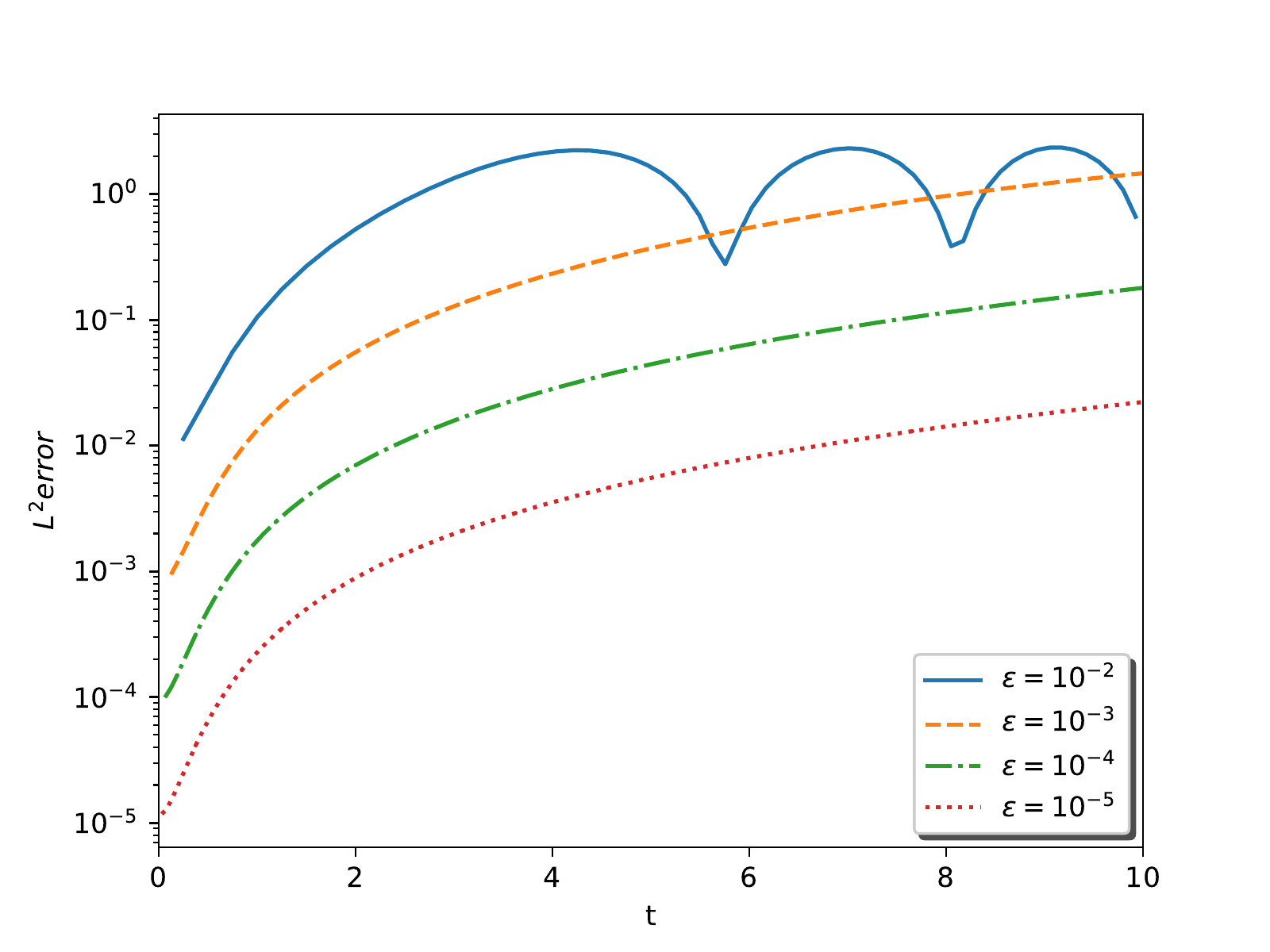}
		\end{minipage}
	}
	\subfigure[$k=2$]{
		\begin{minipage}[b]{0.46\textwidth}    
			\includegraphics[width=1\textwidth]{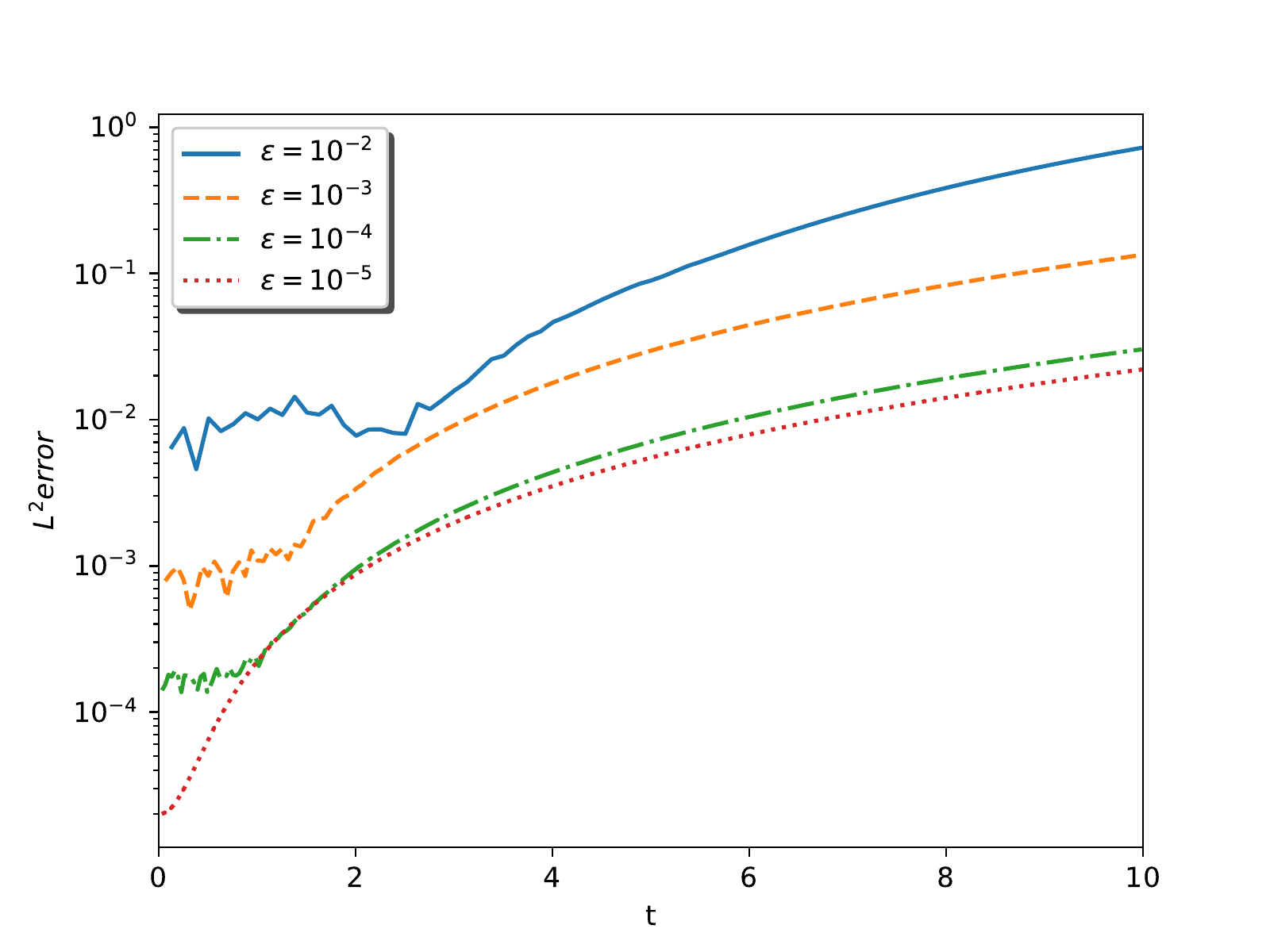}
		\end{minipage}
	}
	\bigskip
	\subfigure[$k=3$]{
		\begin{minipage}[b]{0.46\textwidth}
			\includegraphics[width=1\textwidth]{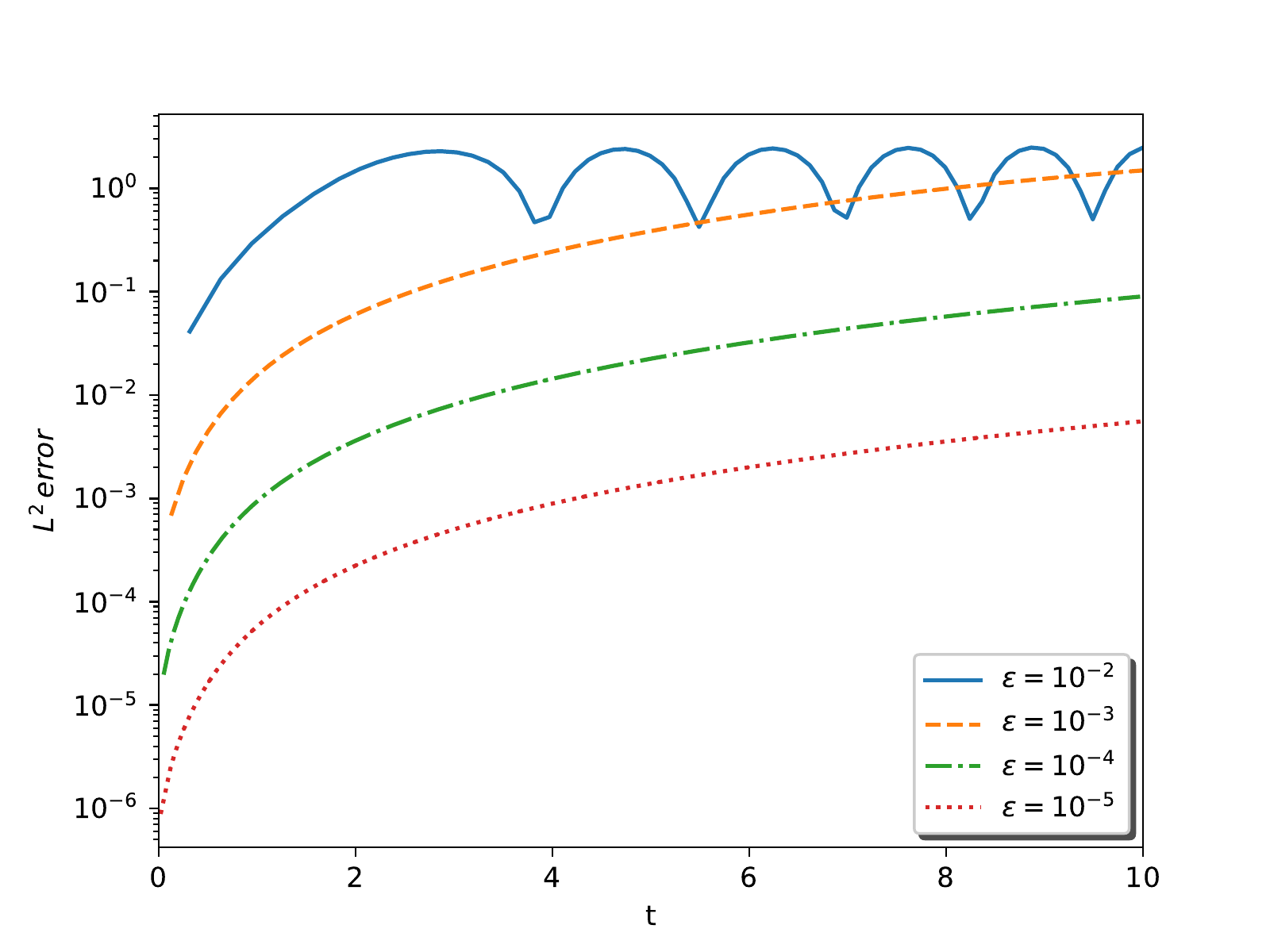}
		\end{minipage}
	}
	\subfigure[$k=3$]{
		\begin{minipage}[b]{0.46\textwidth}    
			\includegraphics[width=1\textwidth]{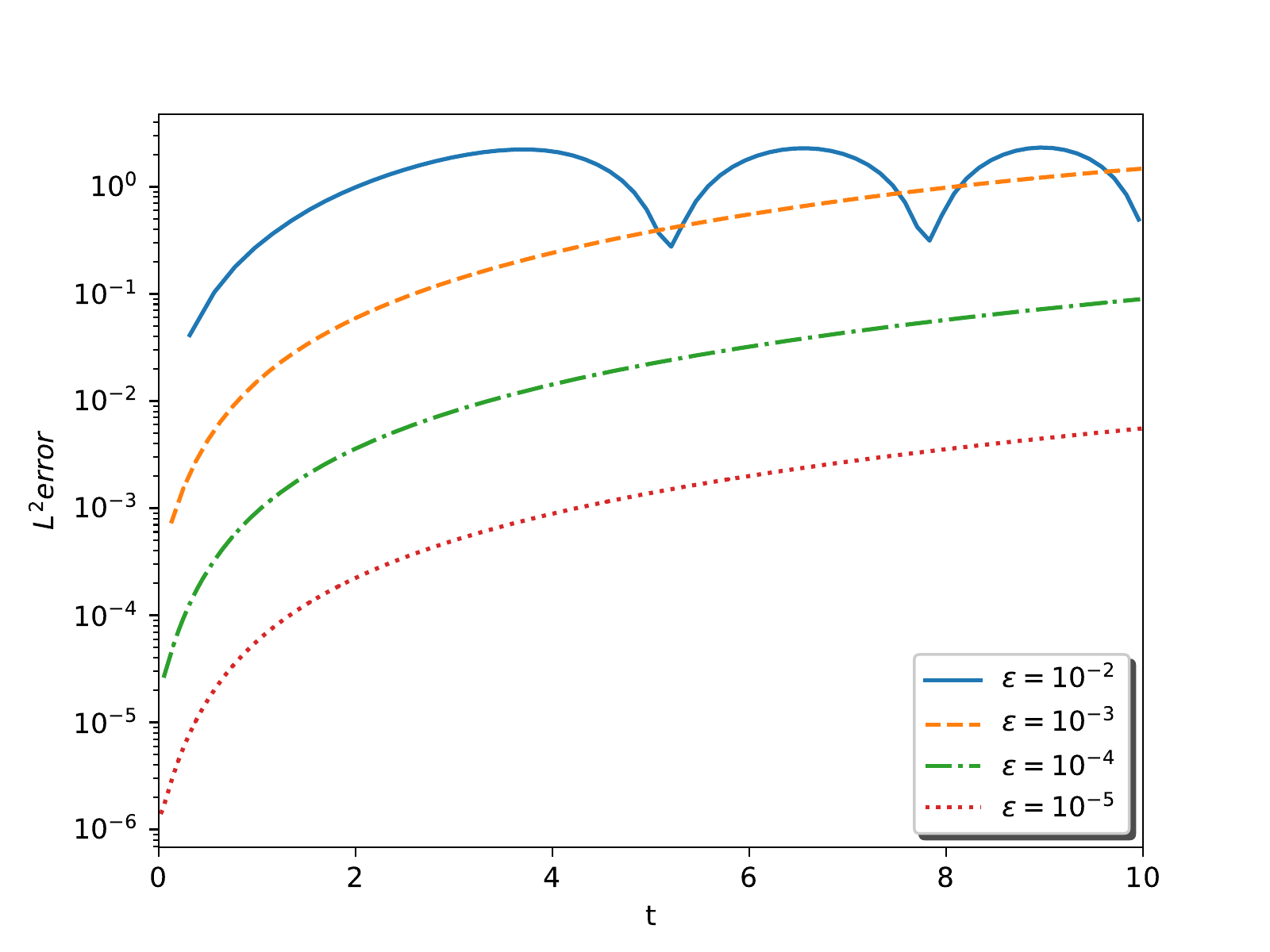}
		\end{minipage}
	}
	\caption{Example \ref{ex:accuracy}: $L^2$-error vs time. $d=1$, $t=10$. $N=8$ and $\eta=\epsilon/10$. Left: conservative numerical flux; right: dissipative numerical flux.}
	\label{fig:error_time_1D}
\end{figure}

Next, we compare the performance of our numerical scheme with the alternating (conservative) numerical flux and the dissipative numerical flux. The time history of the $L^2$-error with different values of polynomial degrees $k$ and error tolerance $\epsilon$ is shown in Figure \ref{fig:error_time_1D}. Note that the adaptive scheme with dissipative flux and $k=1$ does not converge since the corresponding full grid DG is not consistent. In general, the two kind of numerical fluxes has the similar magnitude of errors. Since the conservative numerical flux performs better in regular DG \cite{chen2019ultra}, we will take the conservative numerical flux in the following examples.

\subsection{NLS equation in 1D}

% \begin{exam}[Soliton propogation for NLS equation in 1D]\label{exam:soliton-1D}
\begin{exam}\label{exam:soliton-1D}
In this example, we show the soliton propagation of the NLS equation \eqref{eq:NLS} in the domain $[0, 1]$:
\begin{equation}
    i u_t + \frac{1}{M^2}u_{xx} + 2\abs{u}^2 u = 0,
\end{equation}
with the initial conditions corresponding to the single soliton \cite{xu2005schrodinger}
\begin{equation}\label{eq:1D-single-init}
    u(x,0) = \sech(X-x_0)\exp(2i(X-x_0)),
\end{equation}
and the double soliton \cite{xu2005schrodinger}
\begin{equation}\label{eq:1D-double-init}
    u(x,0) = \sum_{j=1}^2 \sech(X-x_j)\exp(\frac{1}{2}ic_j(X-x_j)),
\end{equation}
with $X=M(x-\frac12)$. Here the parameters are taken as $M=50$, $x_0=25$, $x_1=-10$, $x_2=10$, $c_1=4$ and $c_2=-4$.

\end{exam}

\begin{figure}
	\centering
	\subfigure[numerical solution at $t=0$]{
		\begin{minipage}[b]{0.46\textwidth}
			\includegraphics[width=1\textwidth]{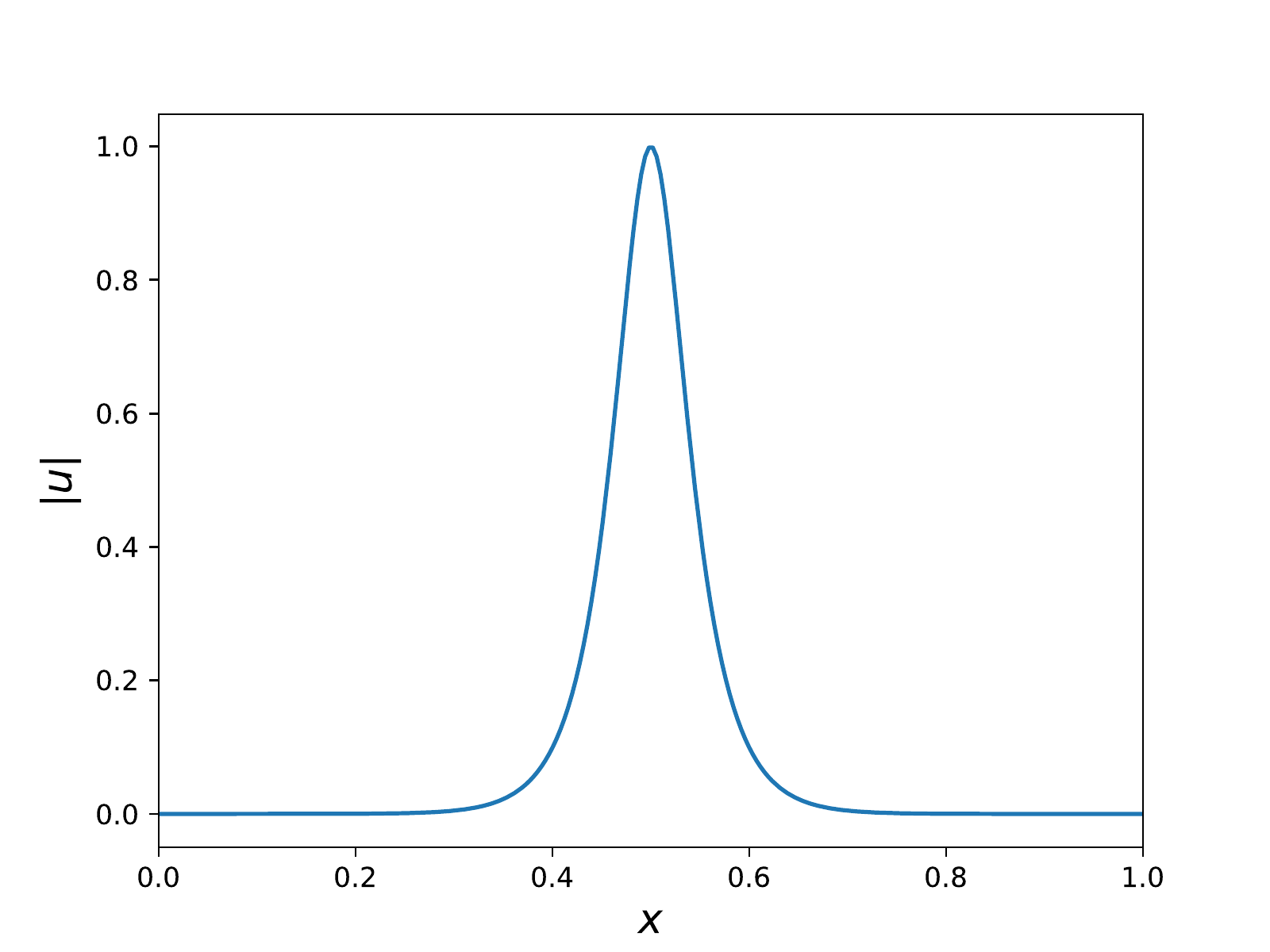}
		\end{minipage}
	}
	\subfigure[active elements at $t=0$]{
		\begin{minipage}[b]{0.46\textwidth}    
			\includegraphics[width=1\textwidth]{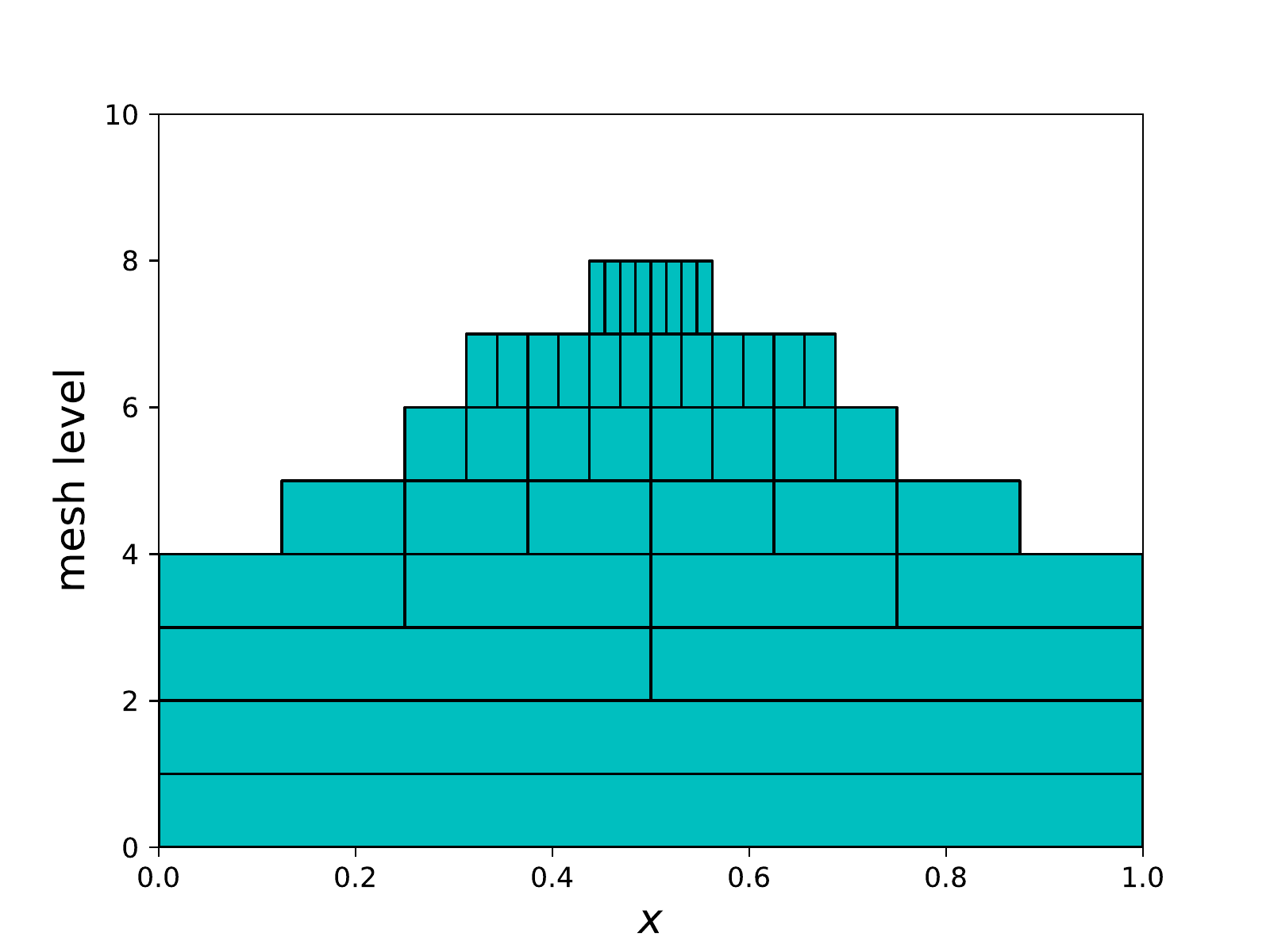}
		\end{minipage}
	}
	\bigskip
	\subfigure[numerical solution at $t=2$]{
		\begin{minipage}[b]{0.46\textwidth}
			\includegraphics[width=1\textwidth]{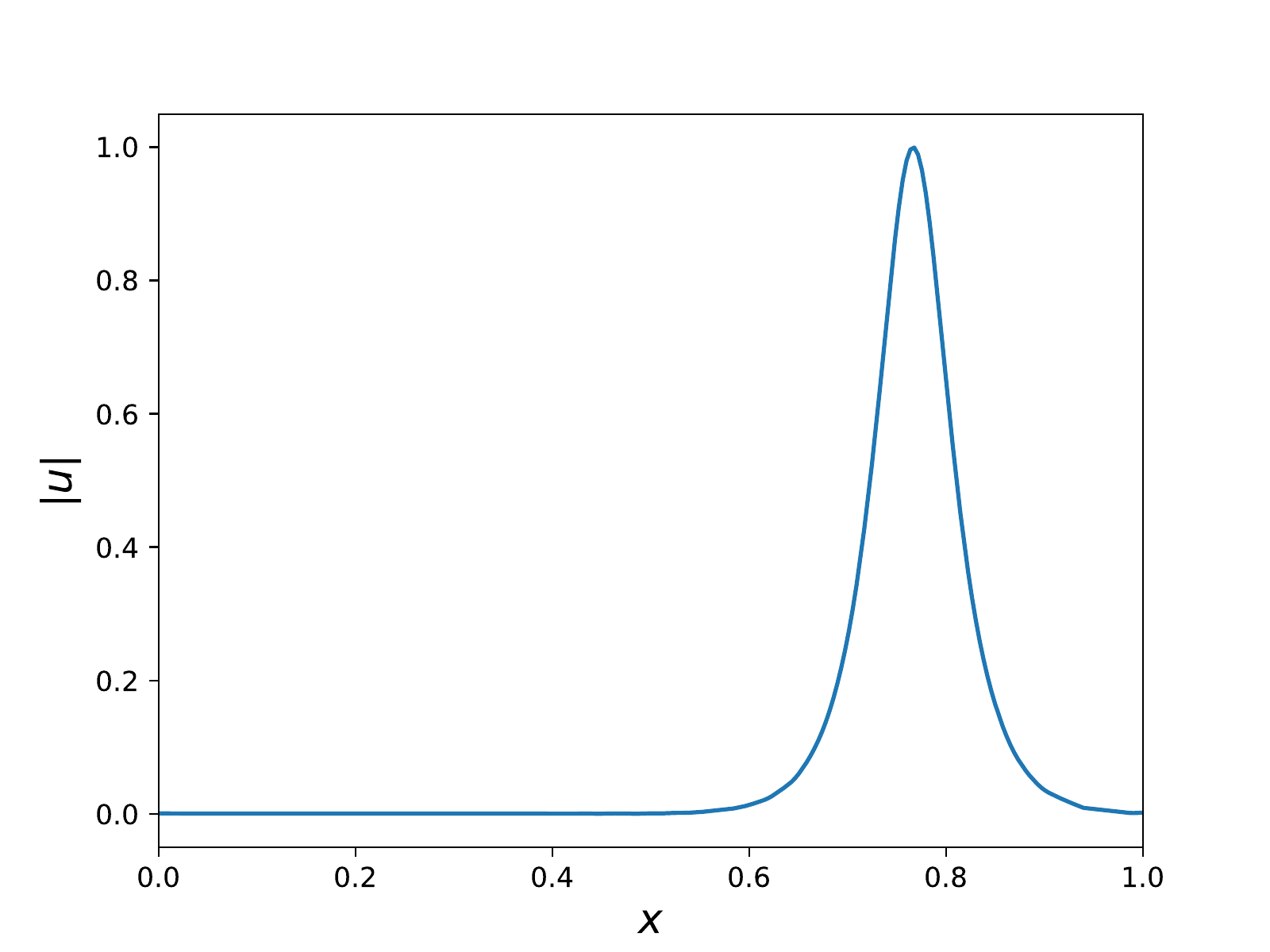}
		\end{minipage}
	}
	\subfigure[active elements at $t=2$]{
		\begin{minipage}[b]{0.46\textwidth}    
			\includegraphics[width=1\textwidth]{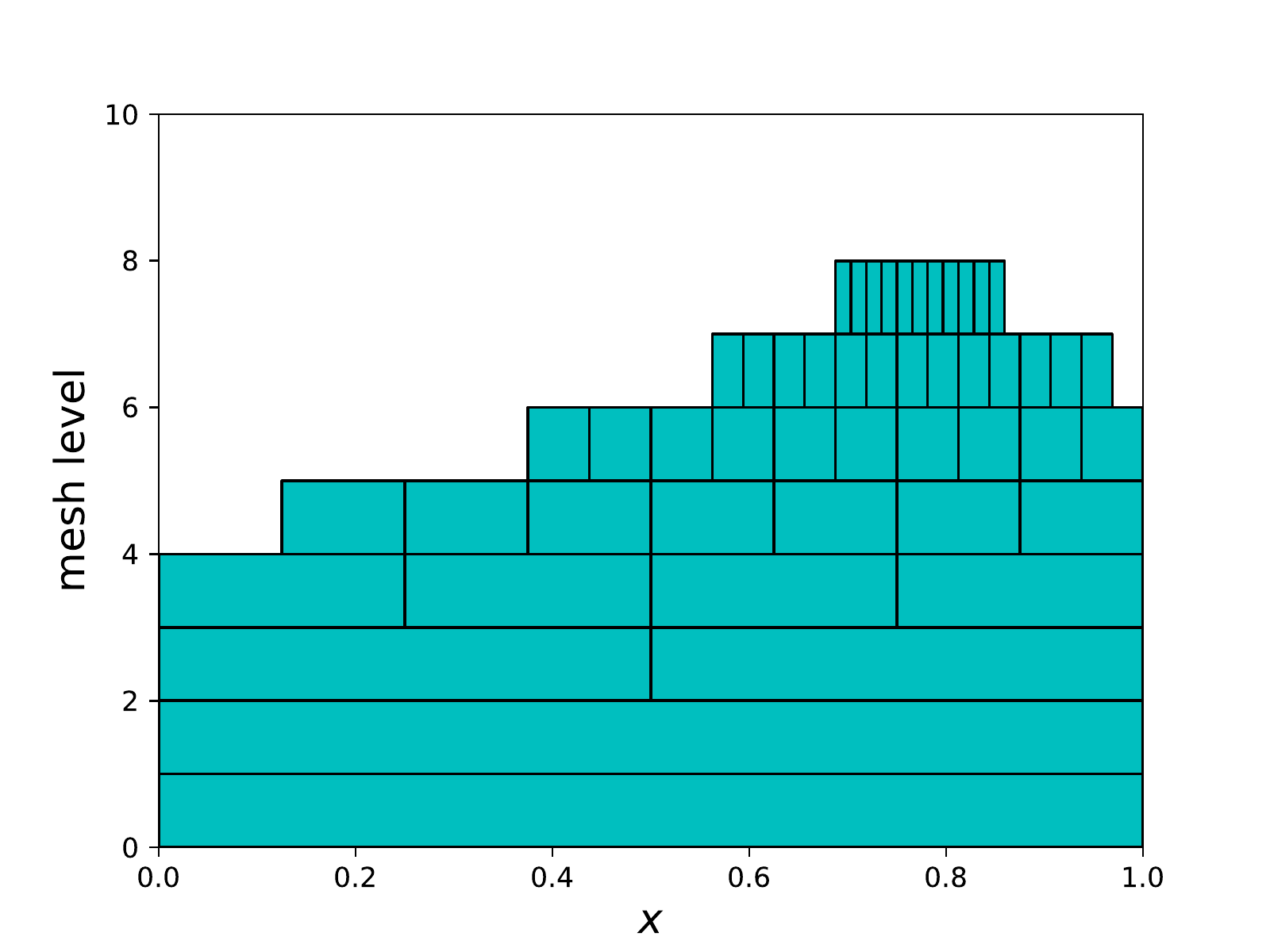}
		\end{minipage}
	}
	\caption{Example \ref{exam:soliton-1D}: 1D NLS equation, single soliton.  Left: numerical solutions; right: active elements. $t=0$ and 2. $N=8$, $\epsilon=10^{-4}$, $\eta=10^{-5}$.}
	\label{fig:single-soliton-1D}
\end{figure}

The numerical solutions and the active elements for the single soliton \eqref{eq:1D-single-init} are shown in Figure \ref{fig:single-soliton-1D}. We observe that the envelope or the modulus $\abs{u}$ are captured by our adaptive scheme quite well. The active elements are also moving with the wave peak.

\begin{figure}
	\centering
	\subfigure[numerical solution at $t=0$]{
		\begin{minipage}[b]{0.46\textwidth}
			\includegraphics[width=1\textwidth]{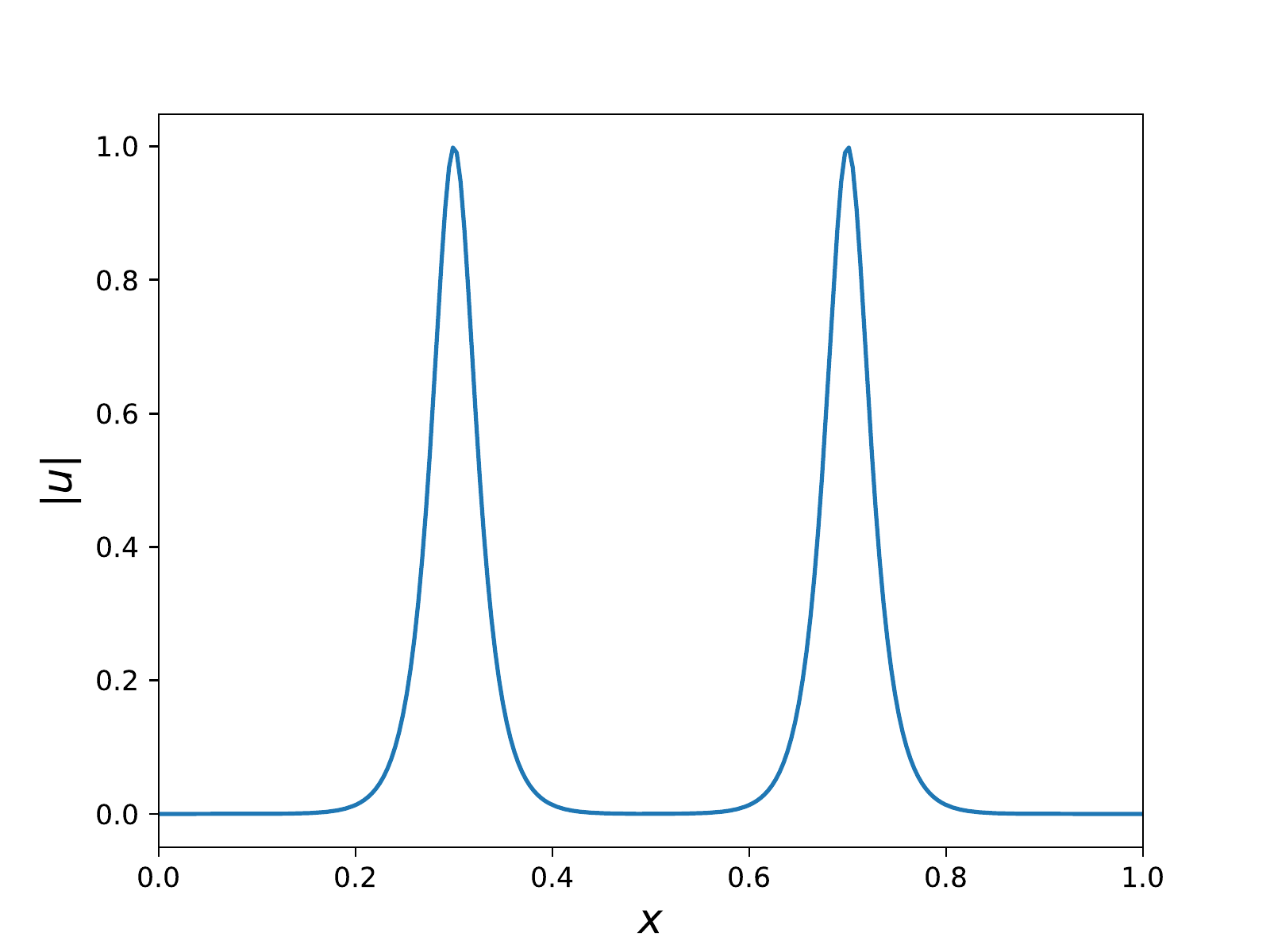}
		\end{minipage}
	}
	\subfigure[active elements at $t=0$]{
		\begin{minipage}[b]{0.46\textwidth}    
			\includegraphics[width=1\textwidth]{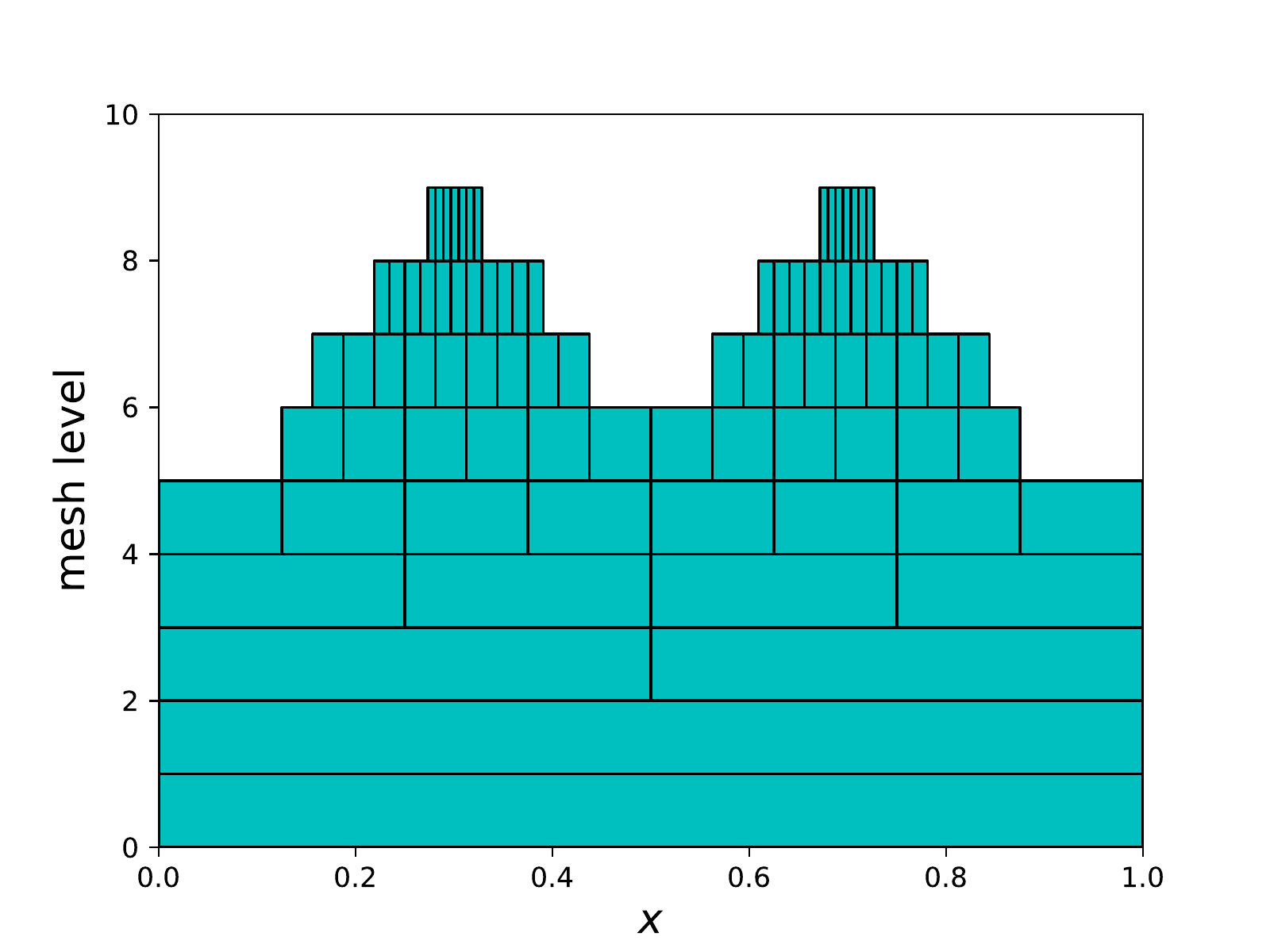}
		\end{minipage}
	}
	\bigskip
	\subfigure[numerical solution at $t=2.5$]{
		\begin{minipage}[b]{0.46\textwidth}
			\includegraphics[width=1\textwidth]{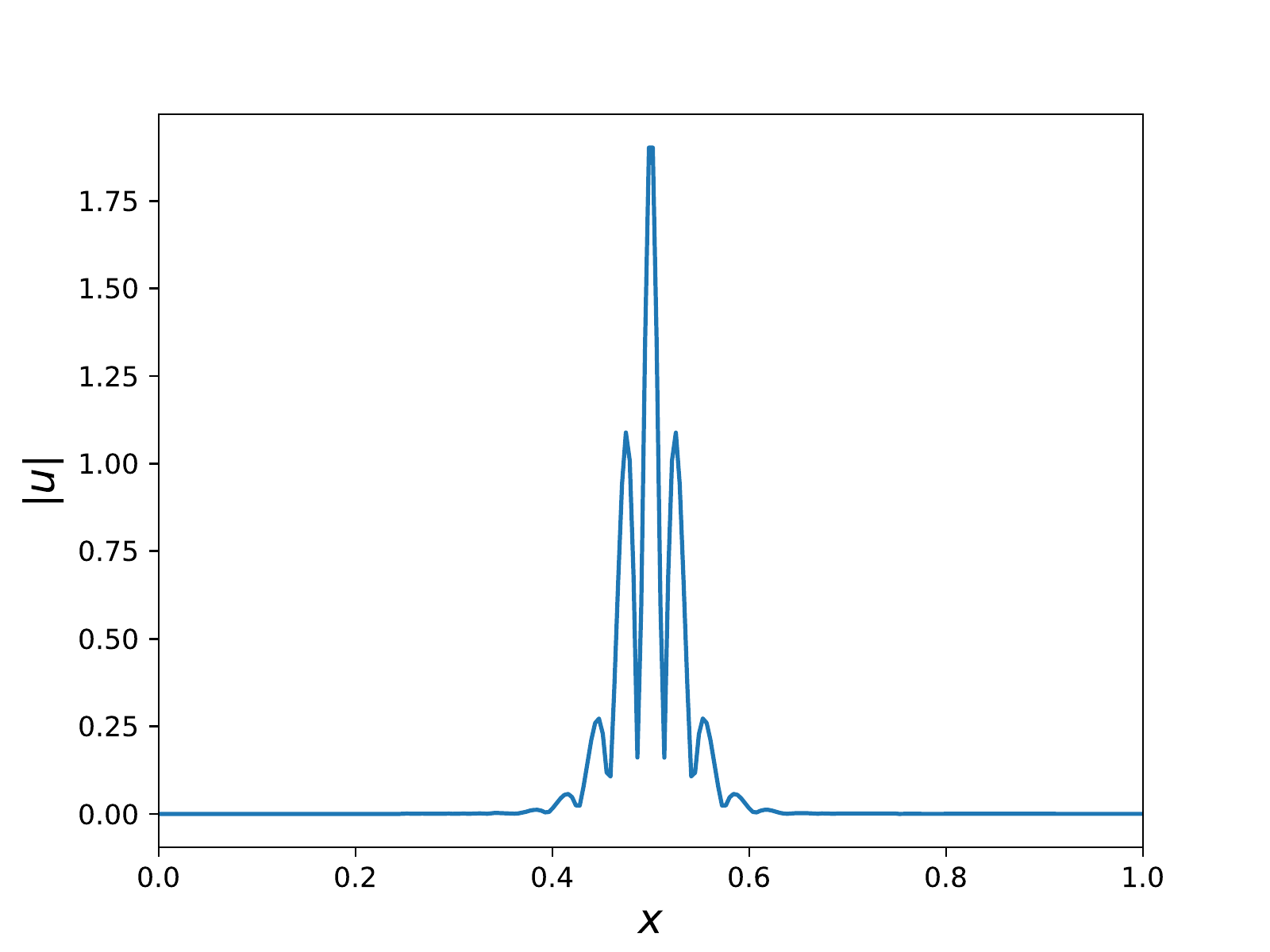}
		\end{minipage}
	}
	\subfigure[active elements at $t=2.5$]{
		\begin{minipage}[b]{0.46\textwidth}    
			\includegraphics[width=1\textwidth]{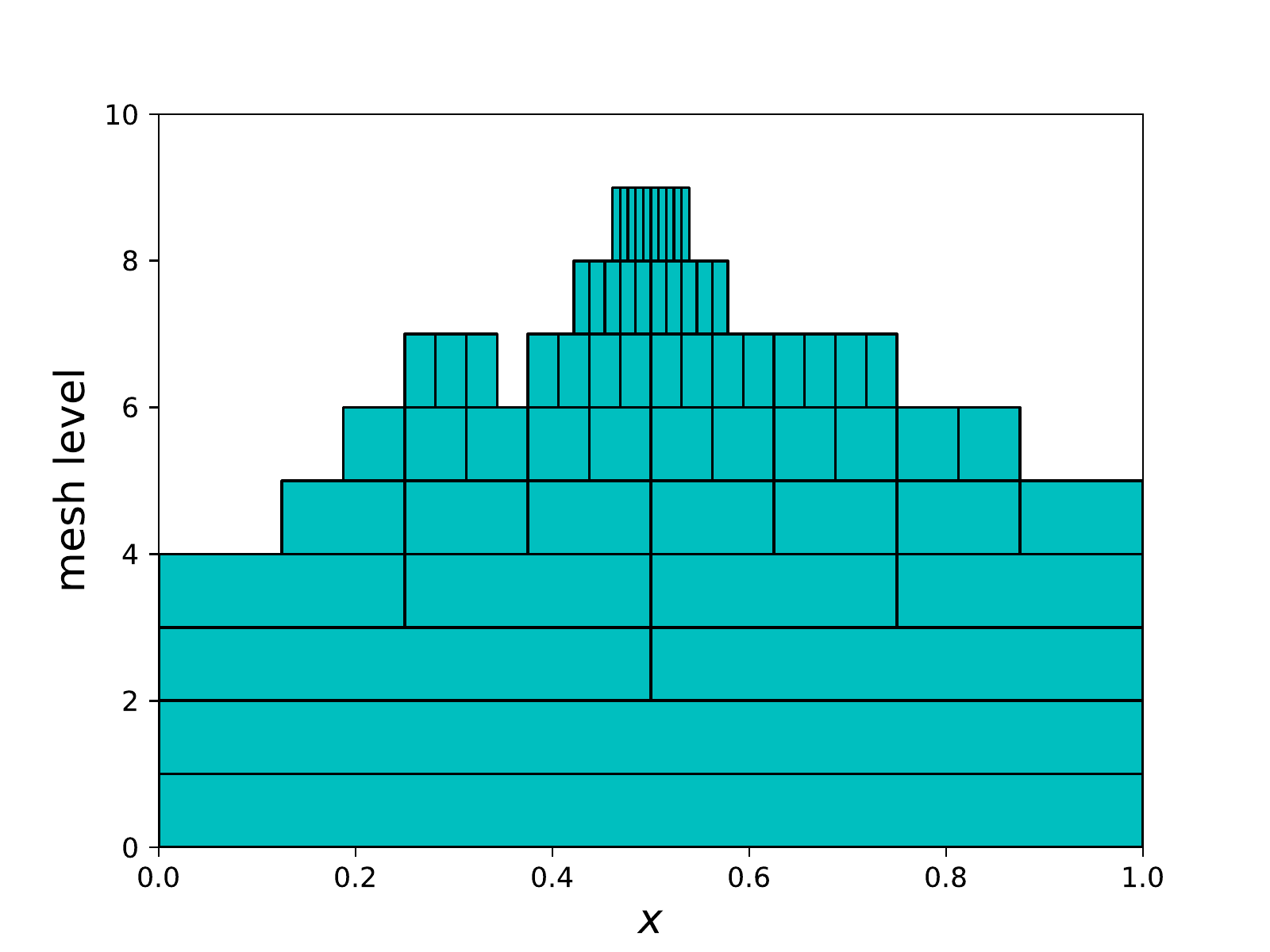}
		\end{minipage}
	}
	\bigskip
	\subfigure[numerical solution at $t=5$]{
		\begin{minipage}[b]{0.46\textwidth}
			\includegraphics[width=1\textwidth]{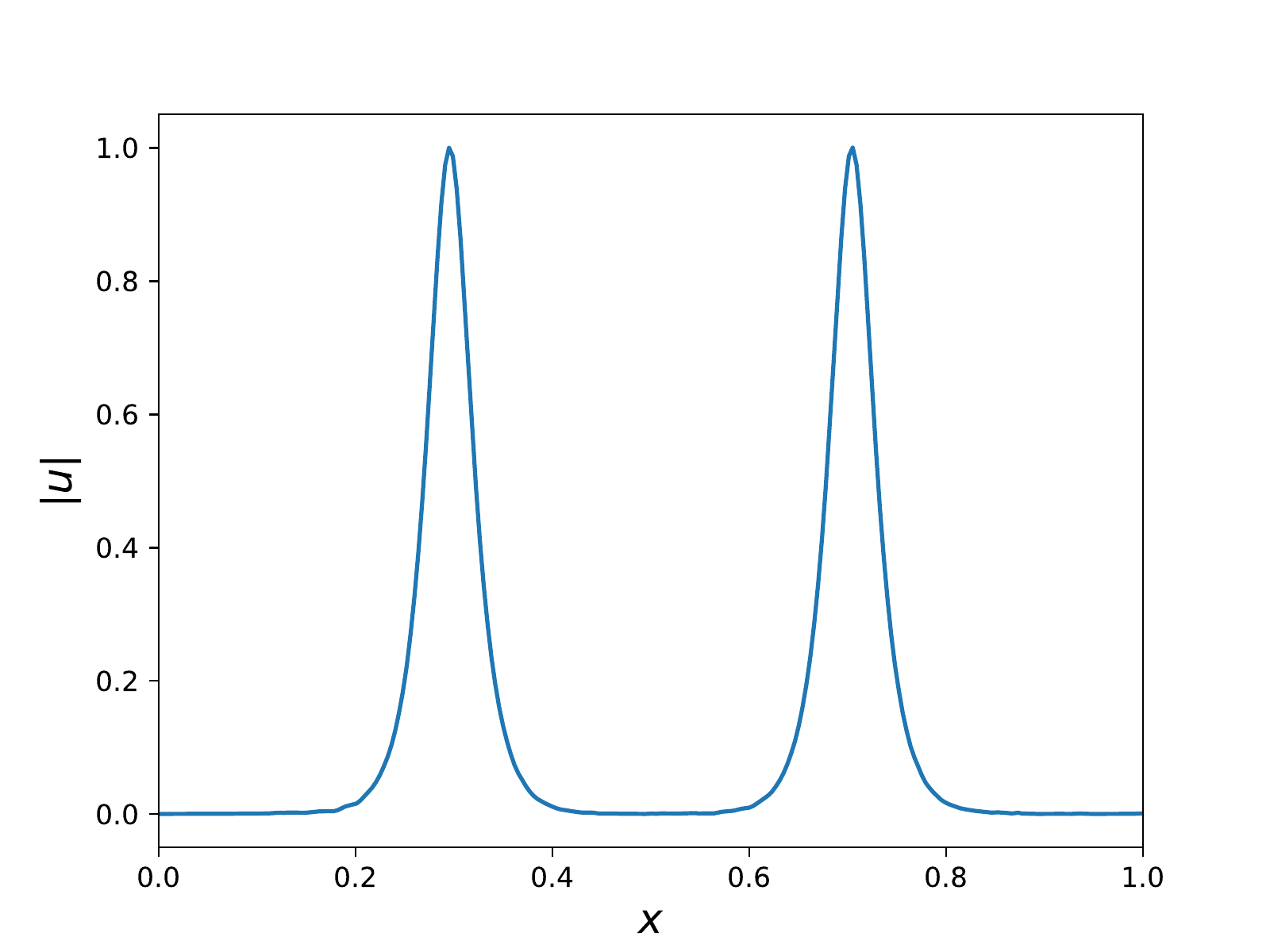}
		\end{minipage}
	}
	\subfigure[active elements at $t=5$]{
		\begin{minipage}[b]{0.46\textwidth}    
			\includegraphics[width=1\textwidth]{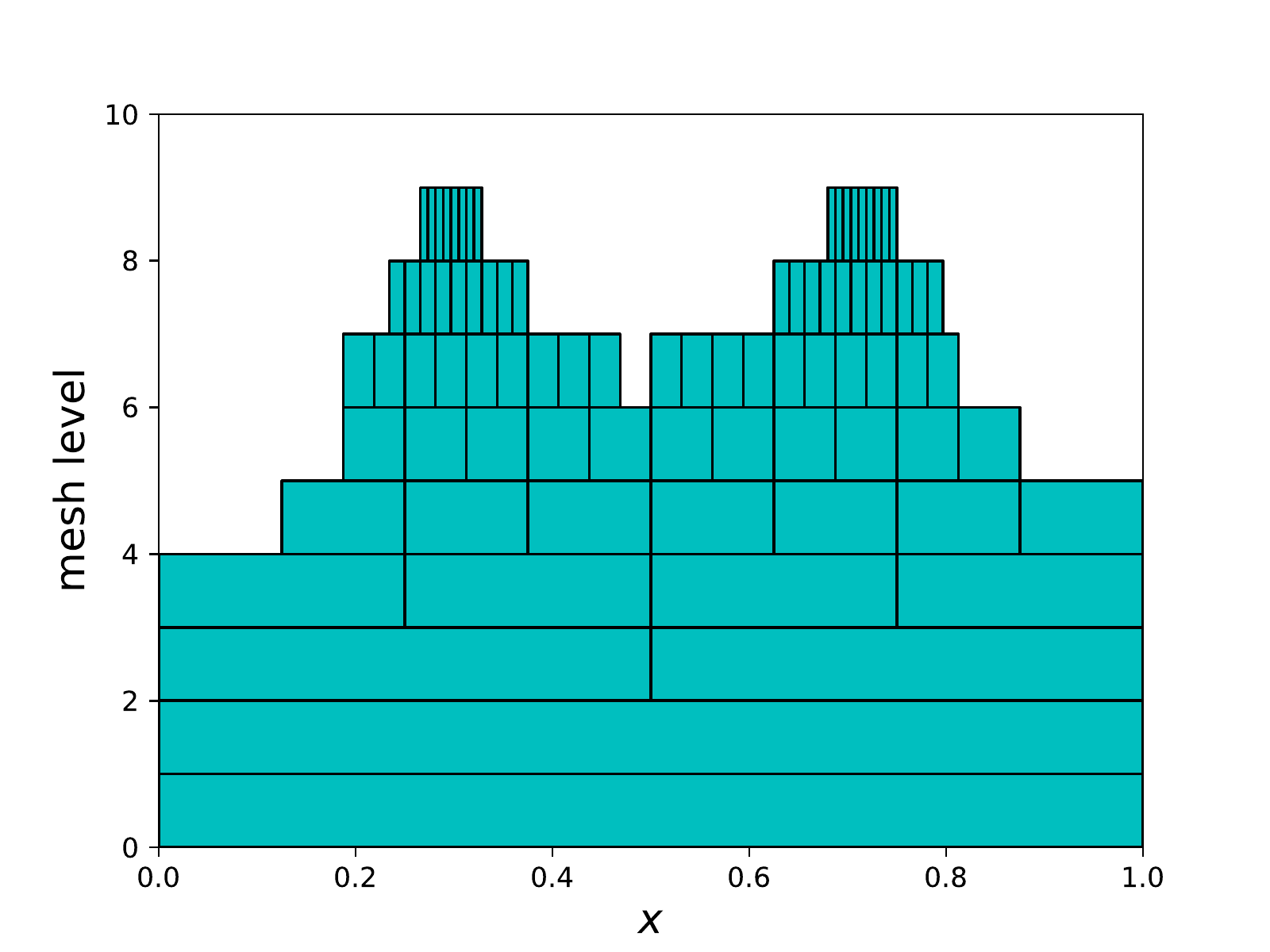}
		\end{minipage}
	}	
	\caption{Example \ref{exam:soliton-1D}: 1D NLS equation, double soliton. Left: numerical solutions; right: active elements. $t=0$, $2.5$ and $5$. $N=8$, $\epsilon=10^{-4}$, $\eta=10^{-5}$. }
	\label{fig:double-soliton-1D}
\end{figure}

The numerical solutions and the active elements for double solitons \eqref{eq:1D-double-init} are shown in Figure \ref{fig:double-soliton-1D}. The two waves propagate in opposite directions and collide at $t=2.5$. After that, the two waves separate. Such behaviors are accurately captured by our numerical simulations. Moreover, our numerical solution does not generate symmetric active elements, which is due to the fact that the ultra-weak DG in full grid does not preserve the symmetry exactly.

\begin{exam}\label{exam:bound_state}
    
In this example, we consider the bound state solution of the equation \cite{xu2005schrodinger}
\begin{equation}
i u_t + \frac{1}{M^2}  u_{xx} +   \beta  \abs{u}^2 u = 0,
\end{equation}
with initial condition
\begin{equation}
u(x,0) = \sech X
\end{equation}
where $X= M(x-0.5), M=30$.

When $\beta = 2L^2$, it will produce a bound state of $L$ solitons.
The theoretical solution for a bound state of solitons is known \cite{miles1981envelope}.
If $L \geq 3$, small narrow structures will develop in the solution which require high mesh resolution to capture. Clearly, using a uniform mesh is far from being optimal due to such a highly localized structure. 
We present the numerical solutions and active elements of the bound state of solitons with $L=3,4,5$ in Figures \ref{fig:bound_state_L3}-\ref{fig:bound_state_L5}. The multiscale structure of the solutions is accurately captured by our adaptive method.

\begin{figure}
	\centering
	\subfigure[numerical solution at $t=0$]{
		\begin{minipage}[b]{0.46\textwidth}
			\includegraphics[width=1\textwidth]{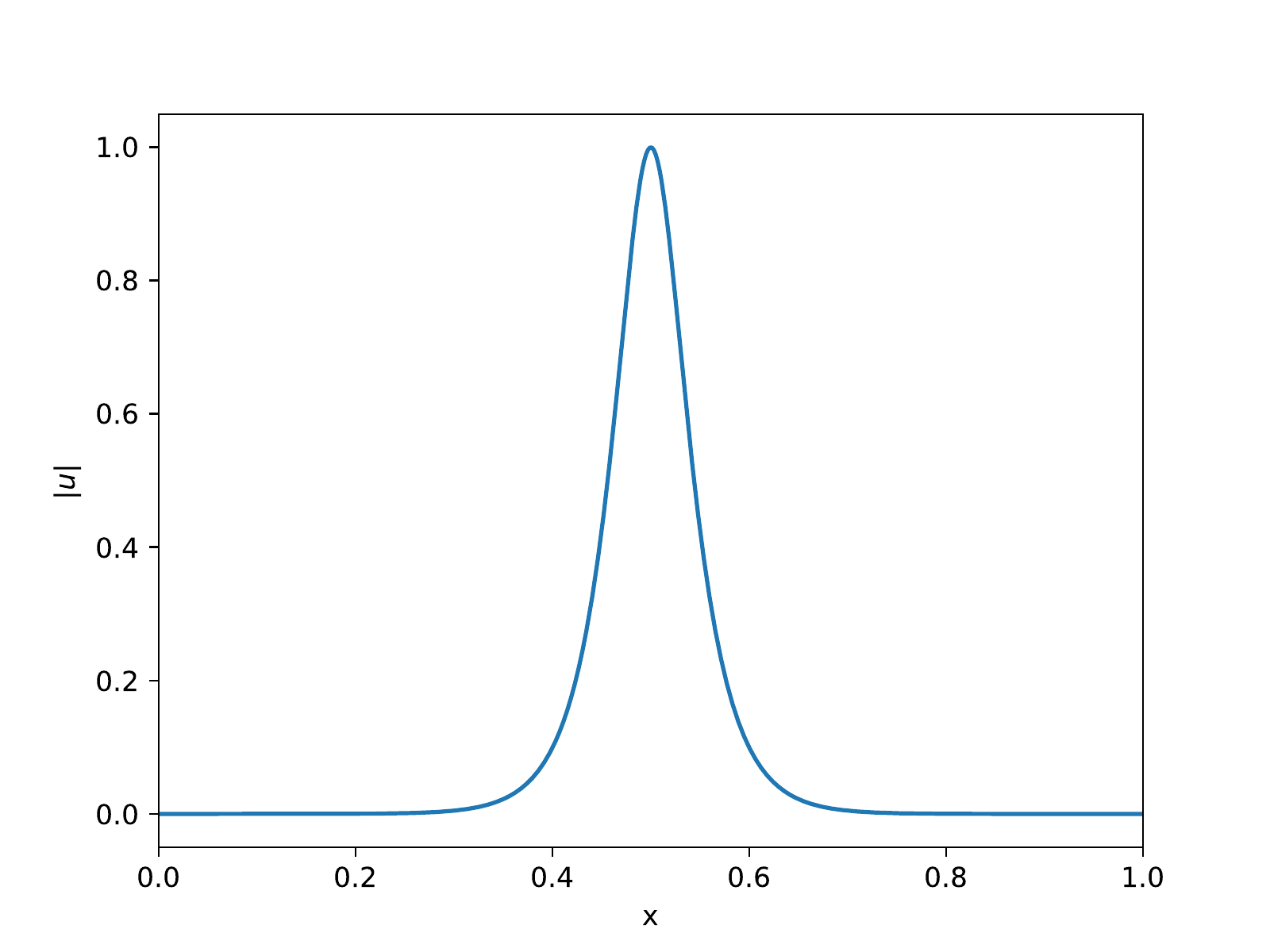}
		\end{minipage}
	}
	\subfigure[active elements at $t=0$]{
		\begin{minipage}[b]{0.46\textwidth}    
			\includegraphics[width=1\textwidth]{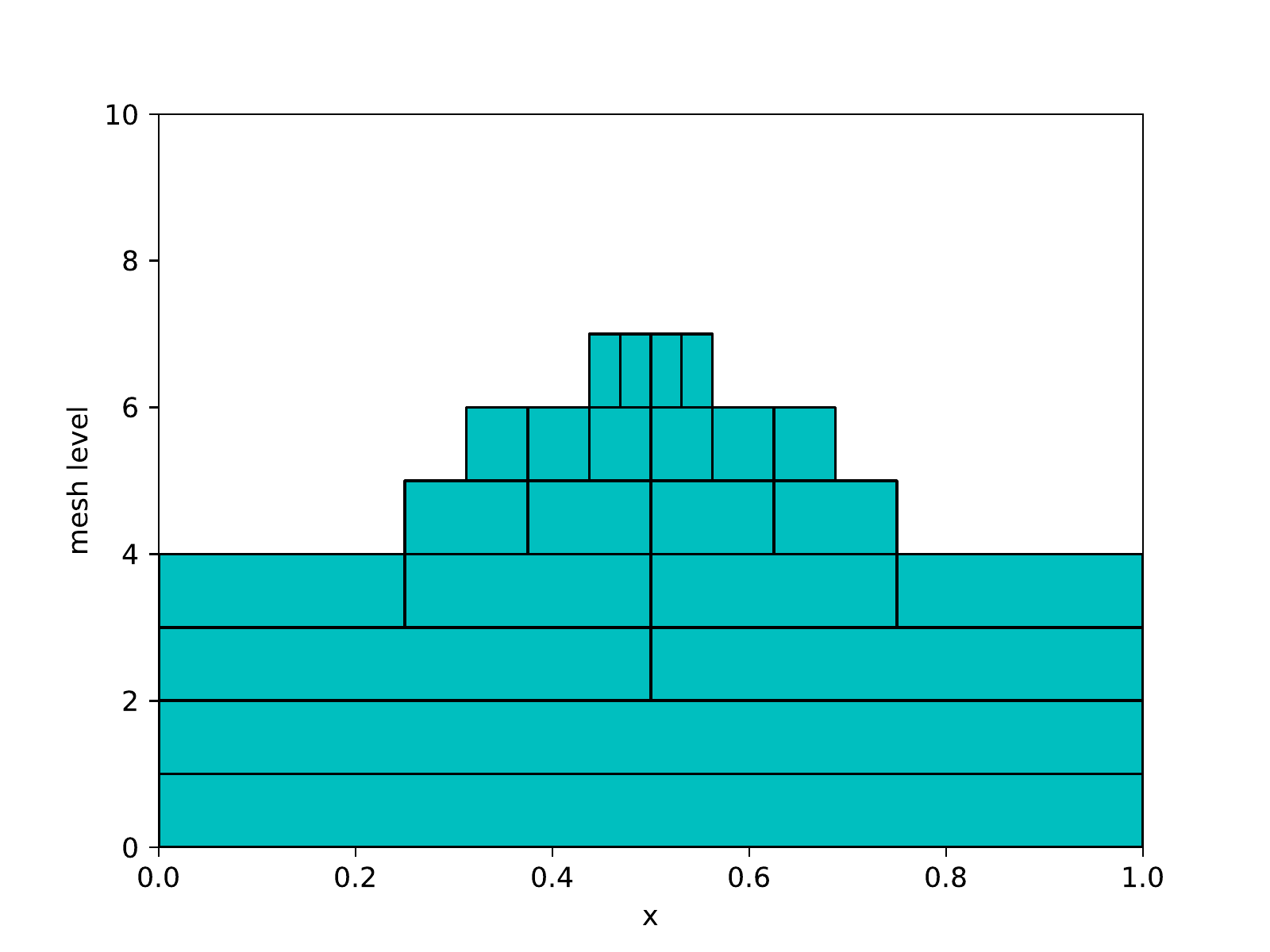}
		\end{minipage}
	}
	\bigskip
	\subfigure[numerical solution at $t=0.4$]{
		\begin{minipage}[b]{0.46\textwidth}
			\includegraphics[width=1\textwidth]{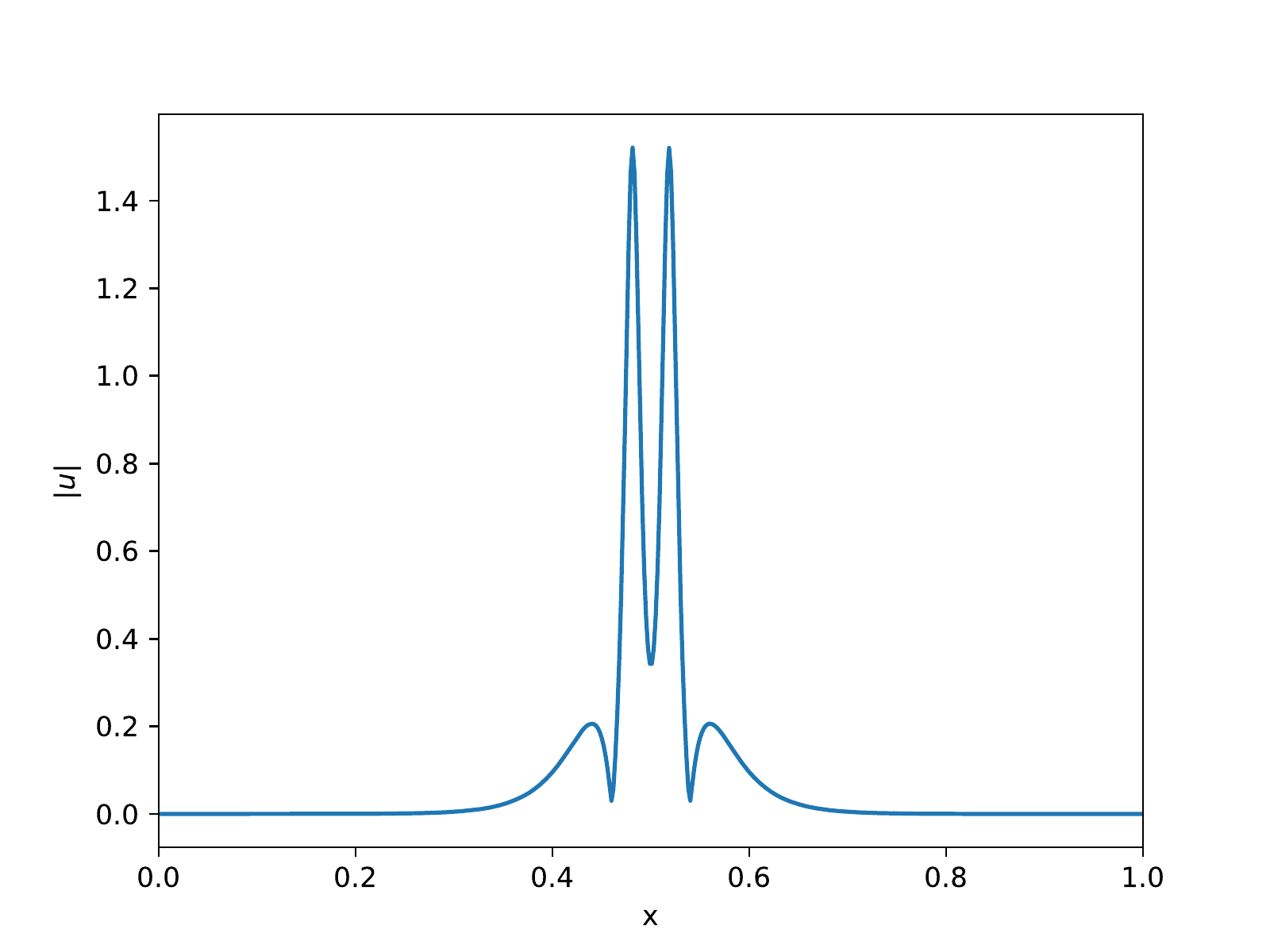}
		\end{minipage}
	}
	\subfigure[active elements at $t=0.4$]{
		\begin{minipage}[b]{0.46\textwidth}    
			\includegraphics[width=1\textwidth]{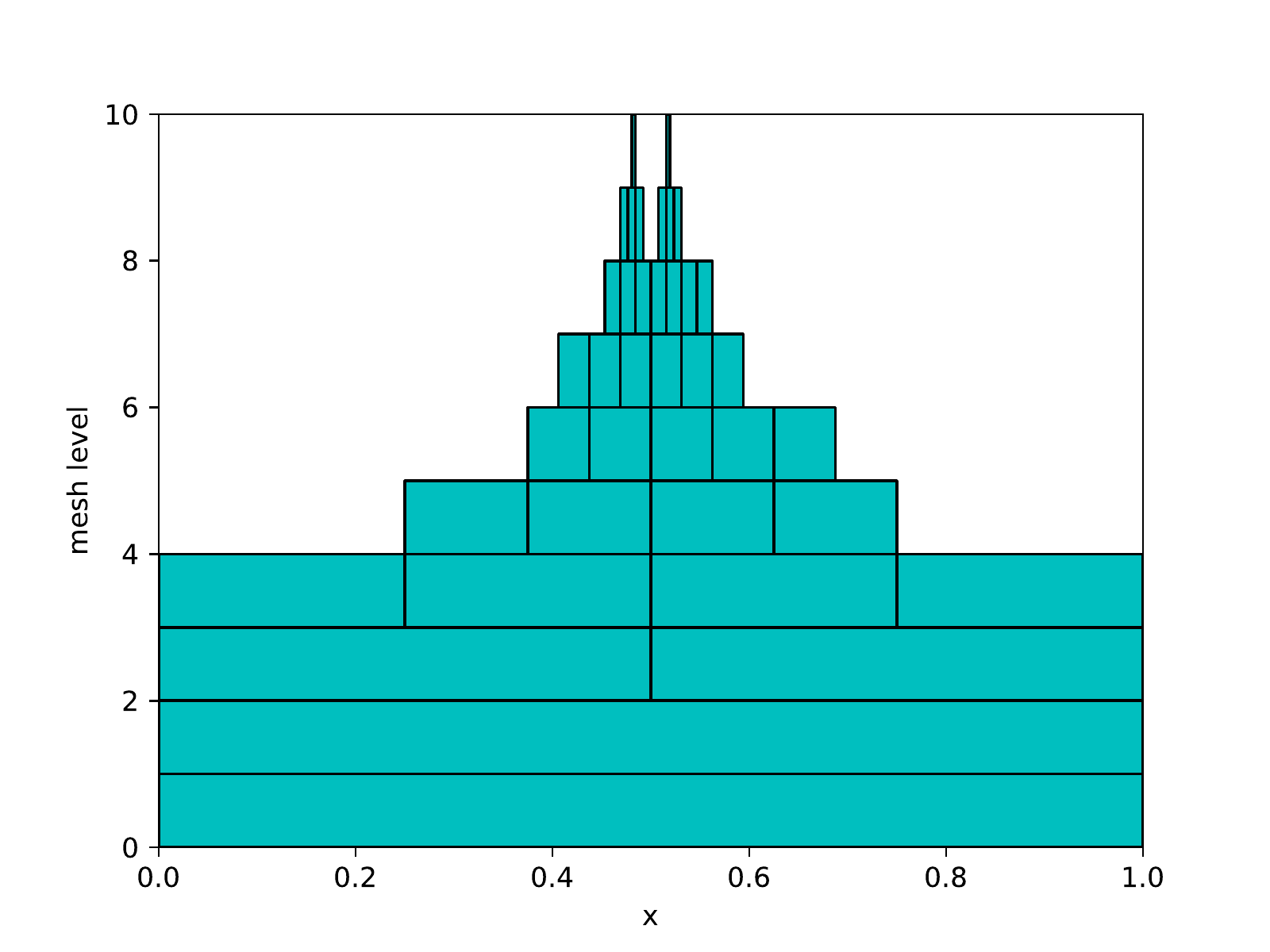}
		\end{minipage}
	}
	\bigskip
	\subfigure[numerical solution at $t=0.6$]{
		\begin{minipage}[b]{0.46\textwidth}
			\includegraphics[width=1\textwidth]{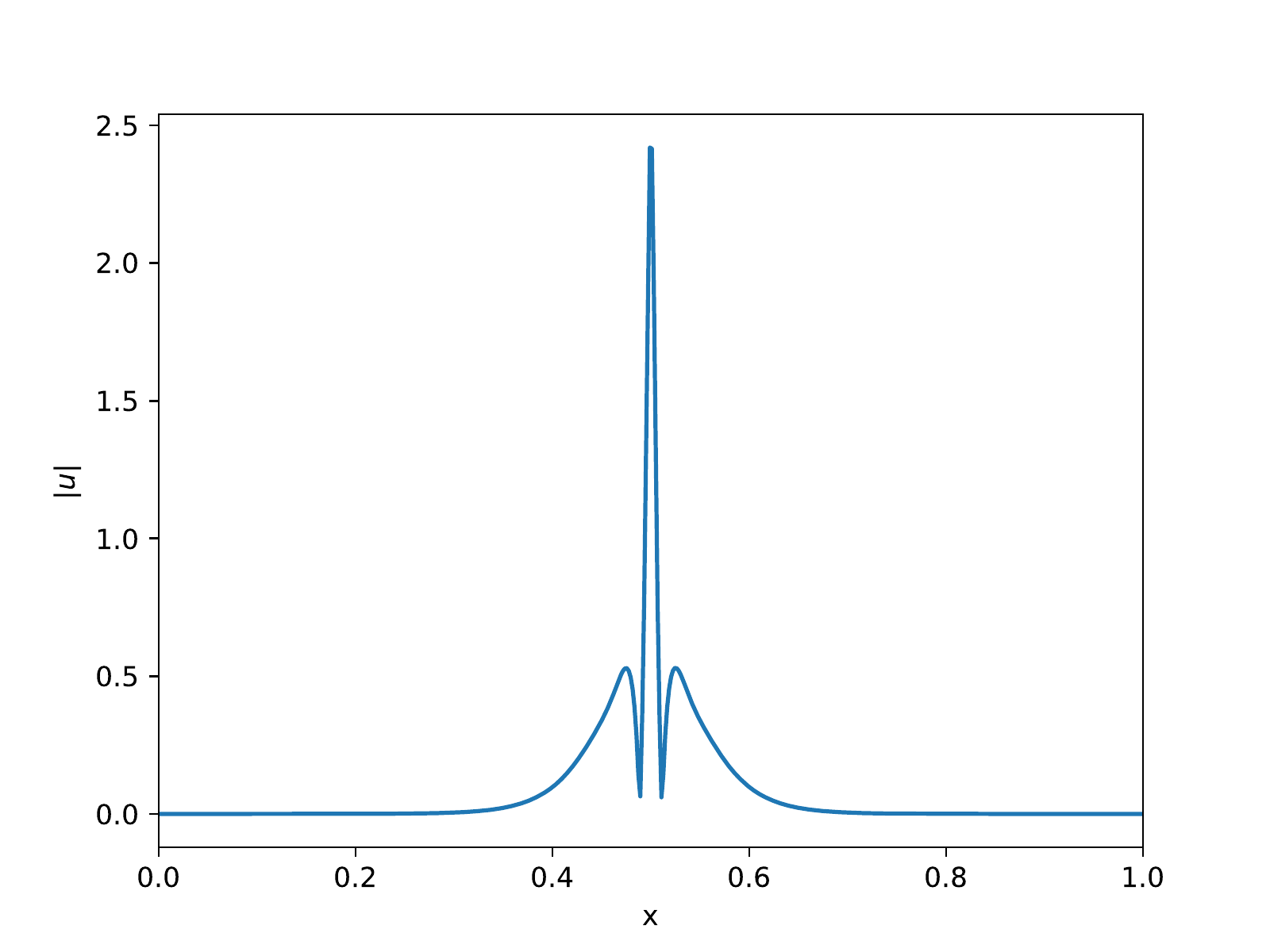}
		\end{minipage}
	}
	\subfigure[active elements at $t=0.6$]{
		\begin{minipage}[b]{0.46\textwidth}    
			\includegraphics[width=1\textwidth]{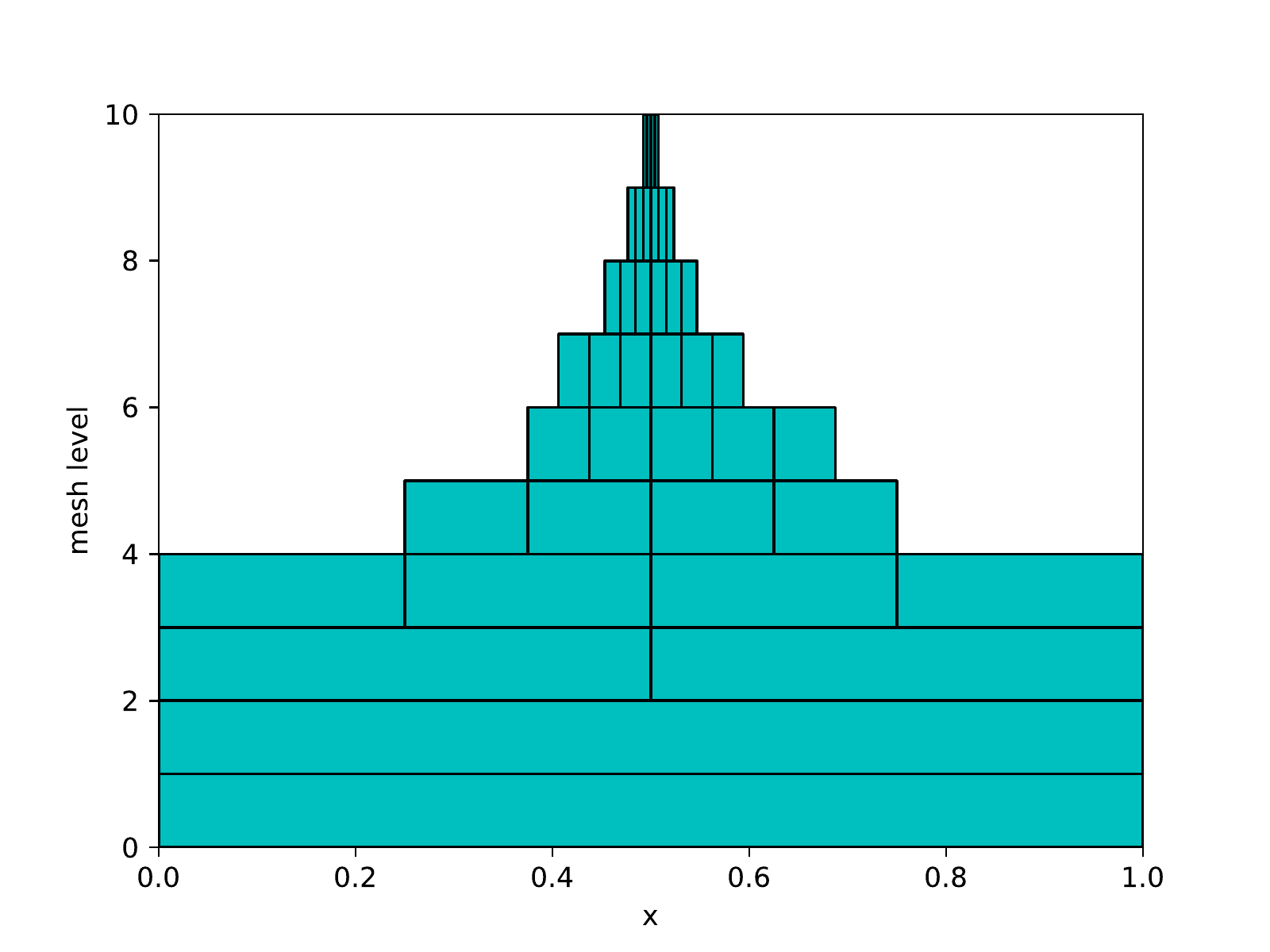}
		\end{minipage}
	}    
	\caption{Example \ref{exam:bound_state}: Bound state solution of solitons with $L=3$.  Left: numerical solutions; right: active elements. $t=0$, 0.4 and 0.6. $N=9$, $k=3$, $\epsilon=10^{-4}$ and $\eta= 10^{-5}$.}
	\label{fig:bound_state_L3}
\end{figure}

\begin{figure}
	\centering
	\subfigure[numerical solution at $t=0$]{
		\begin{minipage}[b]{0.46\textwidth}
			\includegraphics[width=1\textwidth]{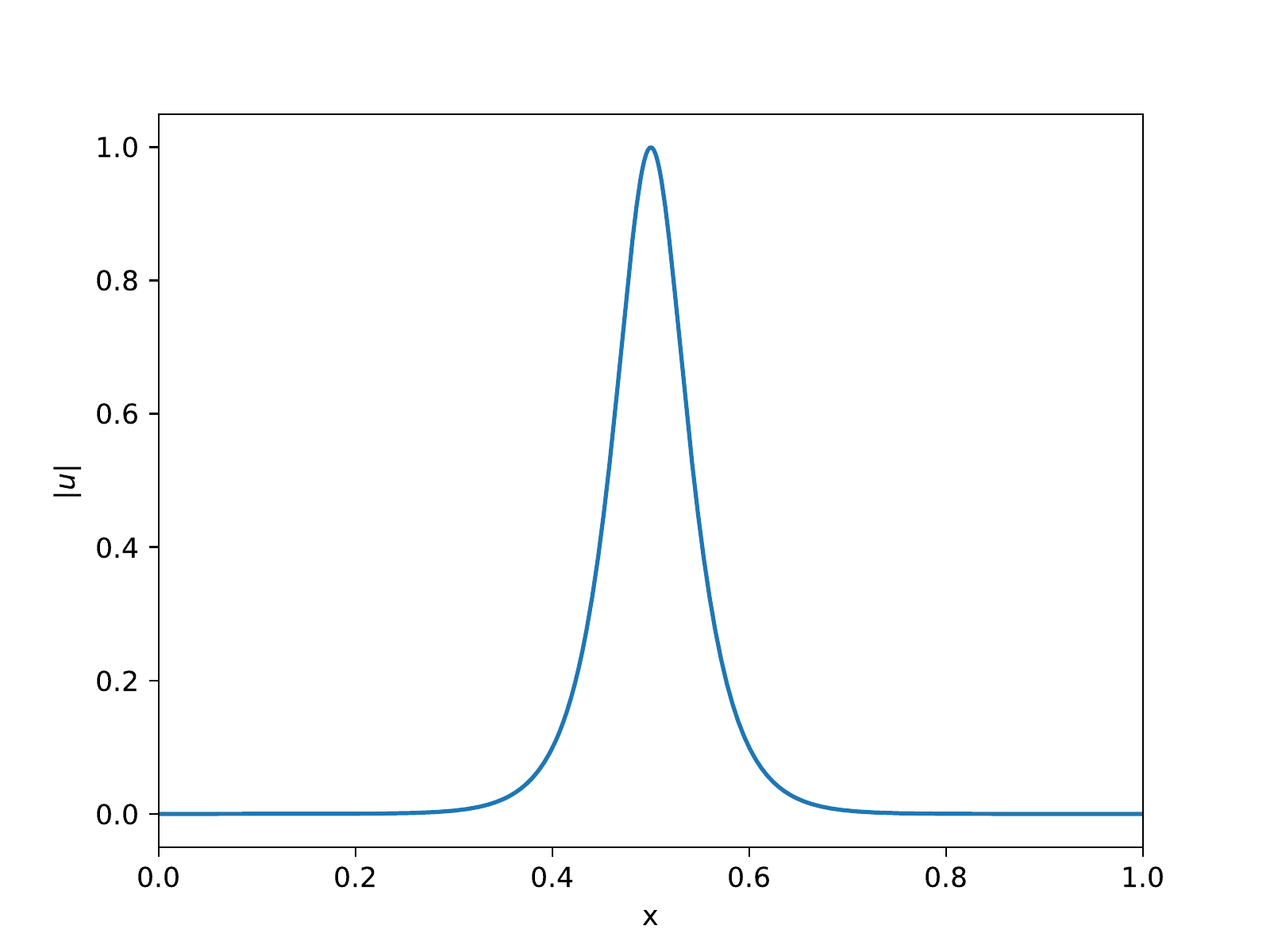}
		\end{minipage}
	}
	\subfigure[active elements at $t=0$]{
		\begin{minipage}[b]{0.46\textwidth}    
			\includegraphics[width=1\textwidth]{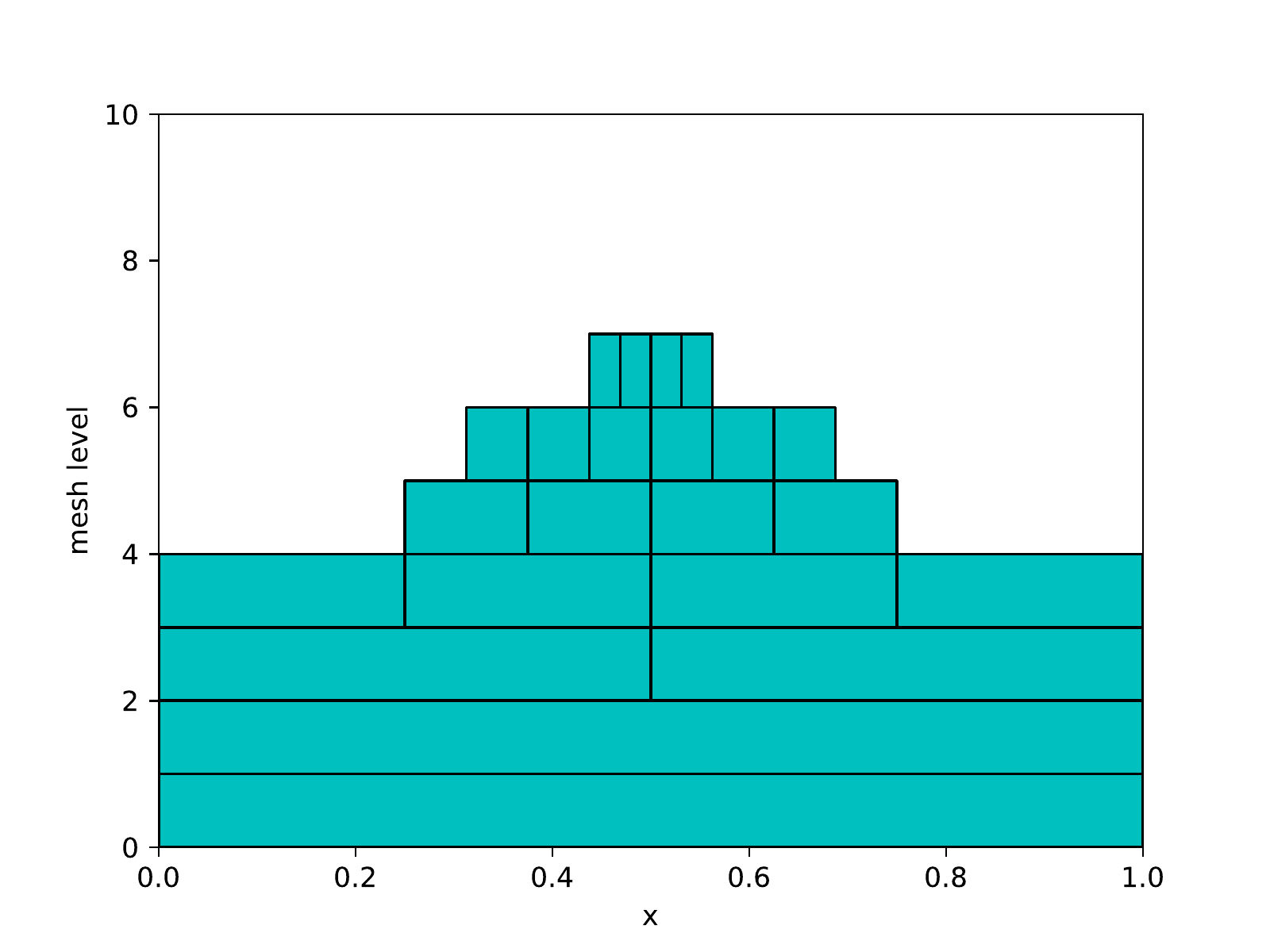}
		\end{minipage}
	}
	\bigskip
	\subfigure[numerical solution at $t=0.4$]{
		\begin{minipage}[b]{0.46\textwidth}
			\includegraphics[width=1\textwidth]{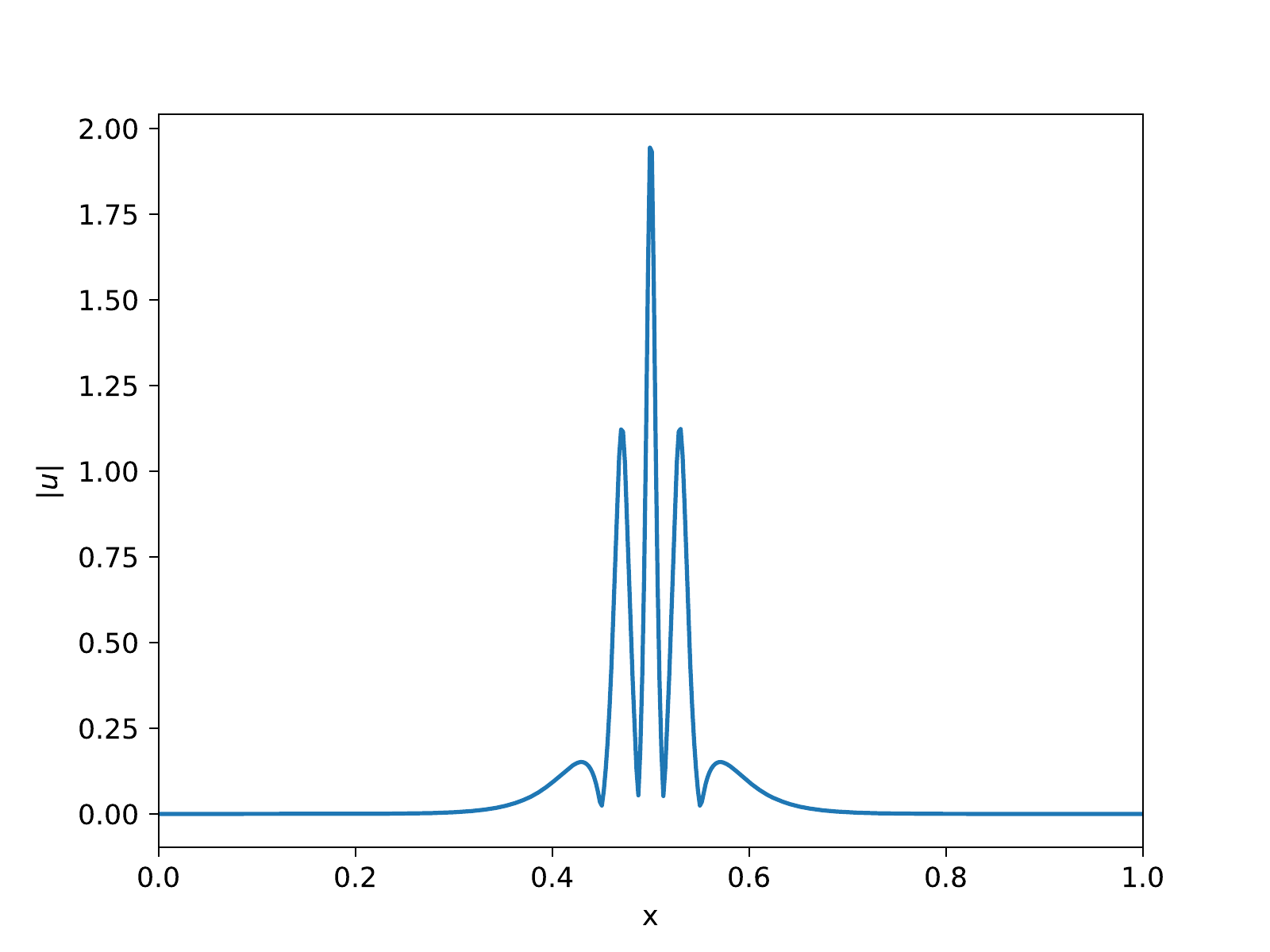}
		\end{minipage}
	}
	\subfigure[active elements at $t=0.4$]{
		\begin{minipage}[b]{0.46\textwidth}    
			\includegraphics[width=1\textwidth]{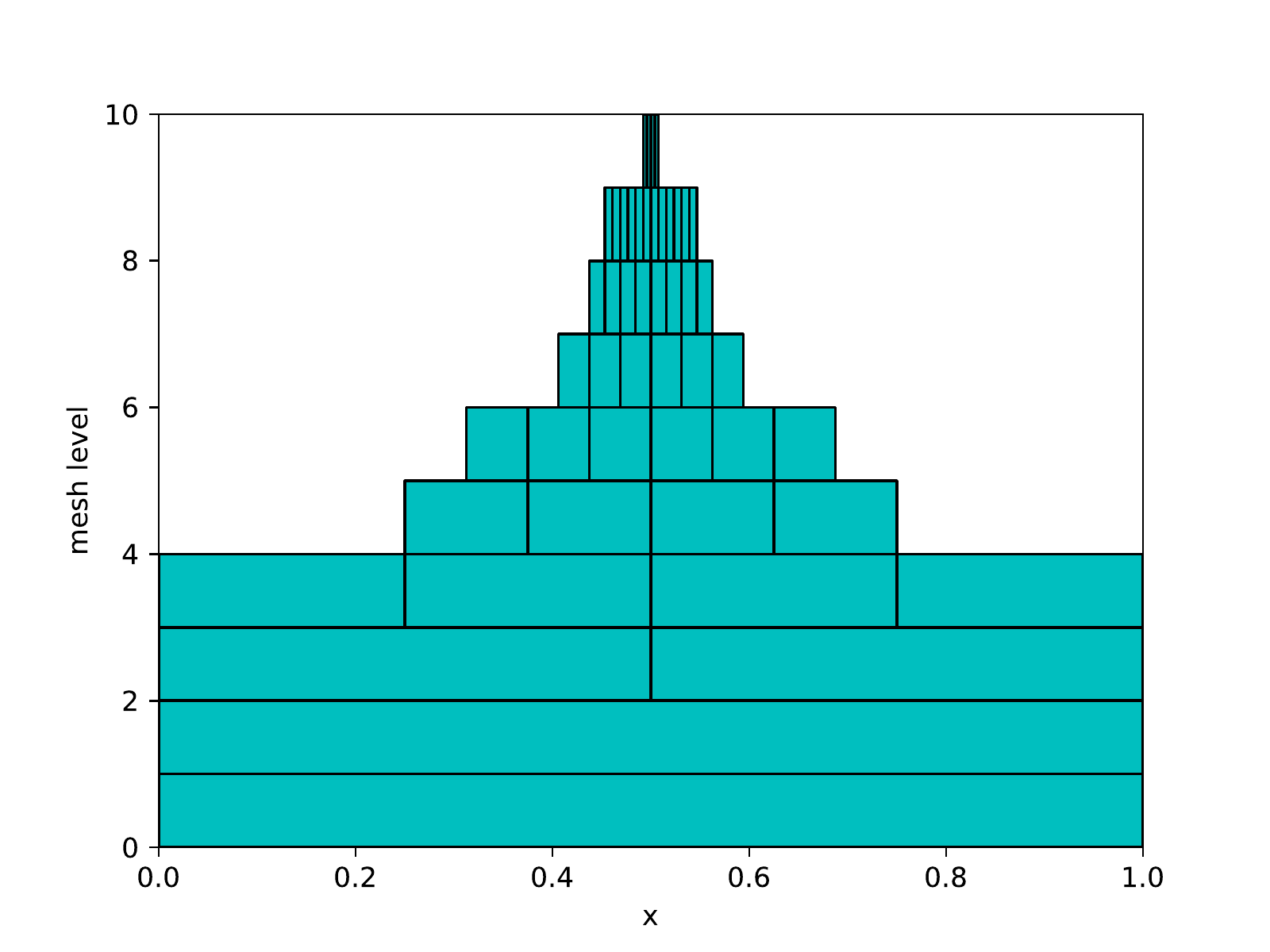}
		\end{minipage}
	}
	\bigskip
	\subfigure[numerical solution at $t=0.6$]{
		\begin{minipage}[b]{0.46\textwidth}
			\includegraphics[width=1\textwidth]{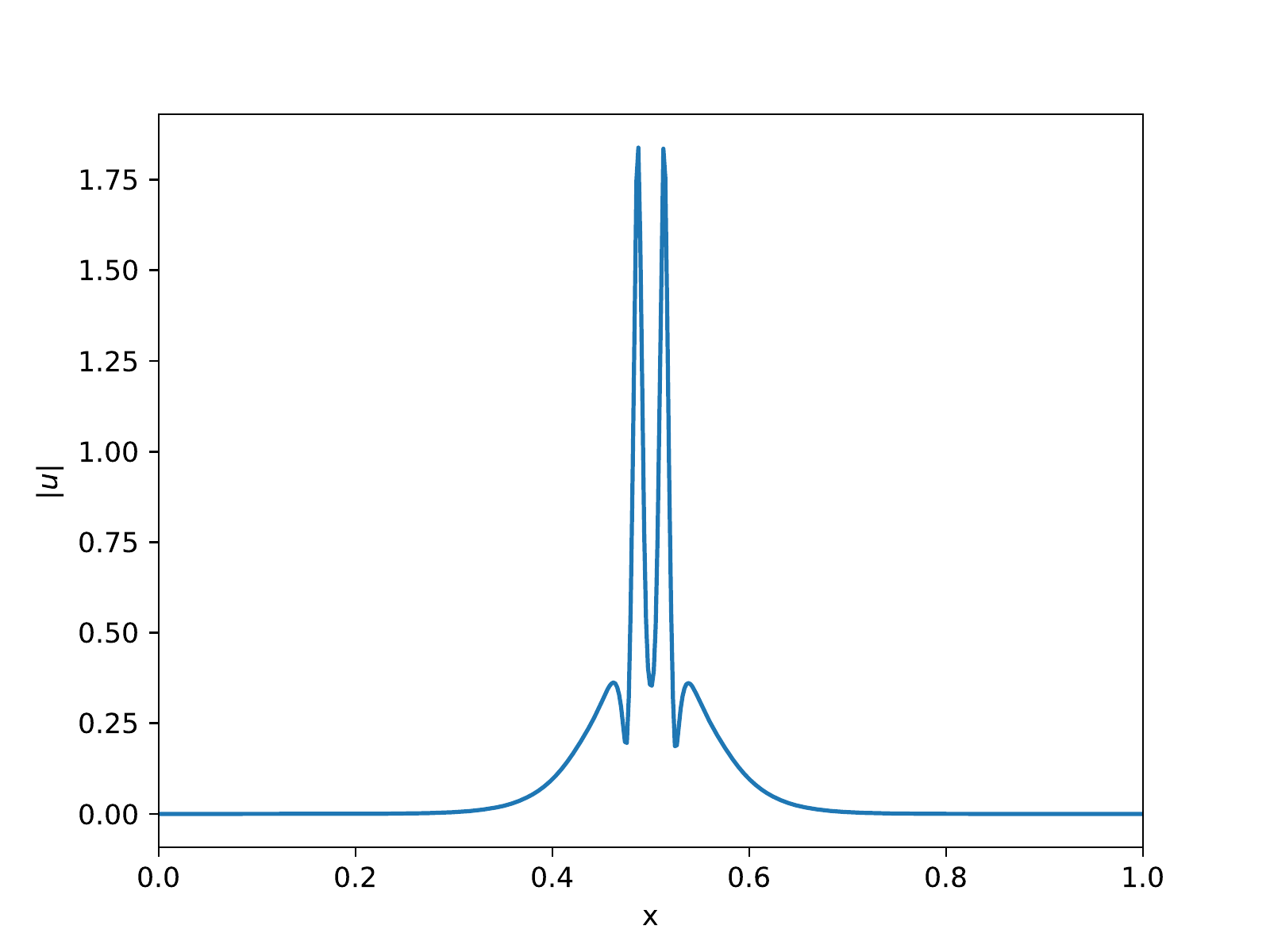}
		\end{minipage}
	}
	\subfigure[active elements at $t=0.6$]{
		\begin{minipage}[b]{0.46\textwidth}    
			\includegraphics[width=1\textwidth]{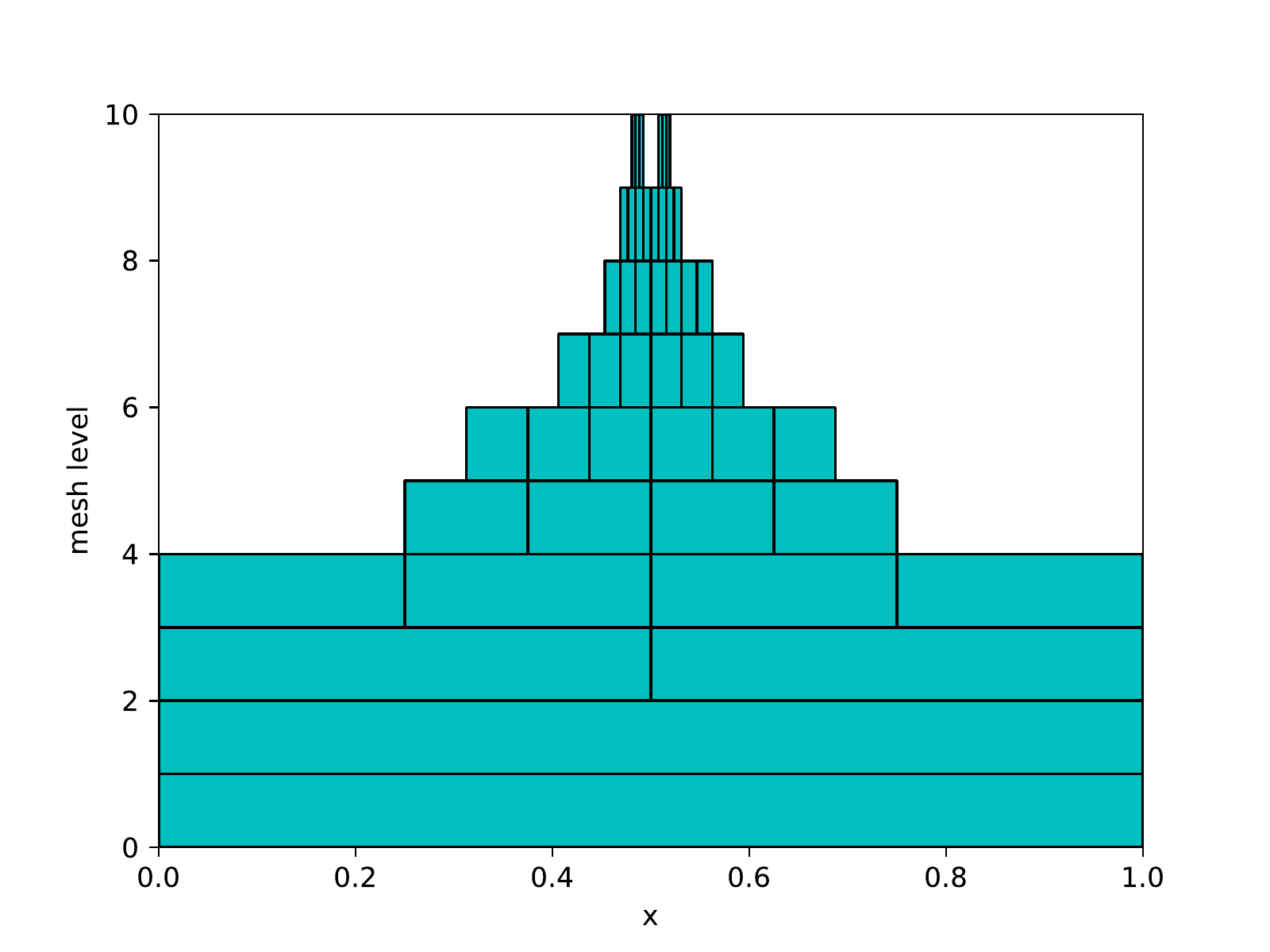}
		\end{minipage}
	}    
	\caption{Example \ref{exam:bound_state}: Bound state solution of solitons with $L=4$.  Left: numerical solutions; right: active elements. $t=0$, 0.4 and 0.6. $N=9, k=3$, $\epsilon=10^{-4}$ and $\eta=10^{-5}$.}
	\label{fig:bound_state_L4}
\end{figure}

\begin{figure}
	\centering
	\subfigure[numerical solution at $t=0$]{
		\begin{minipage}[b]{0.46\textwidth}
			\includegraphics[width=1\textwidth]{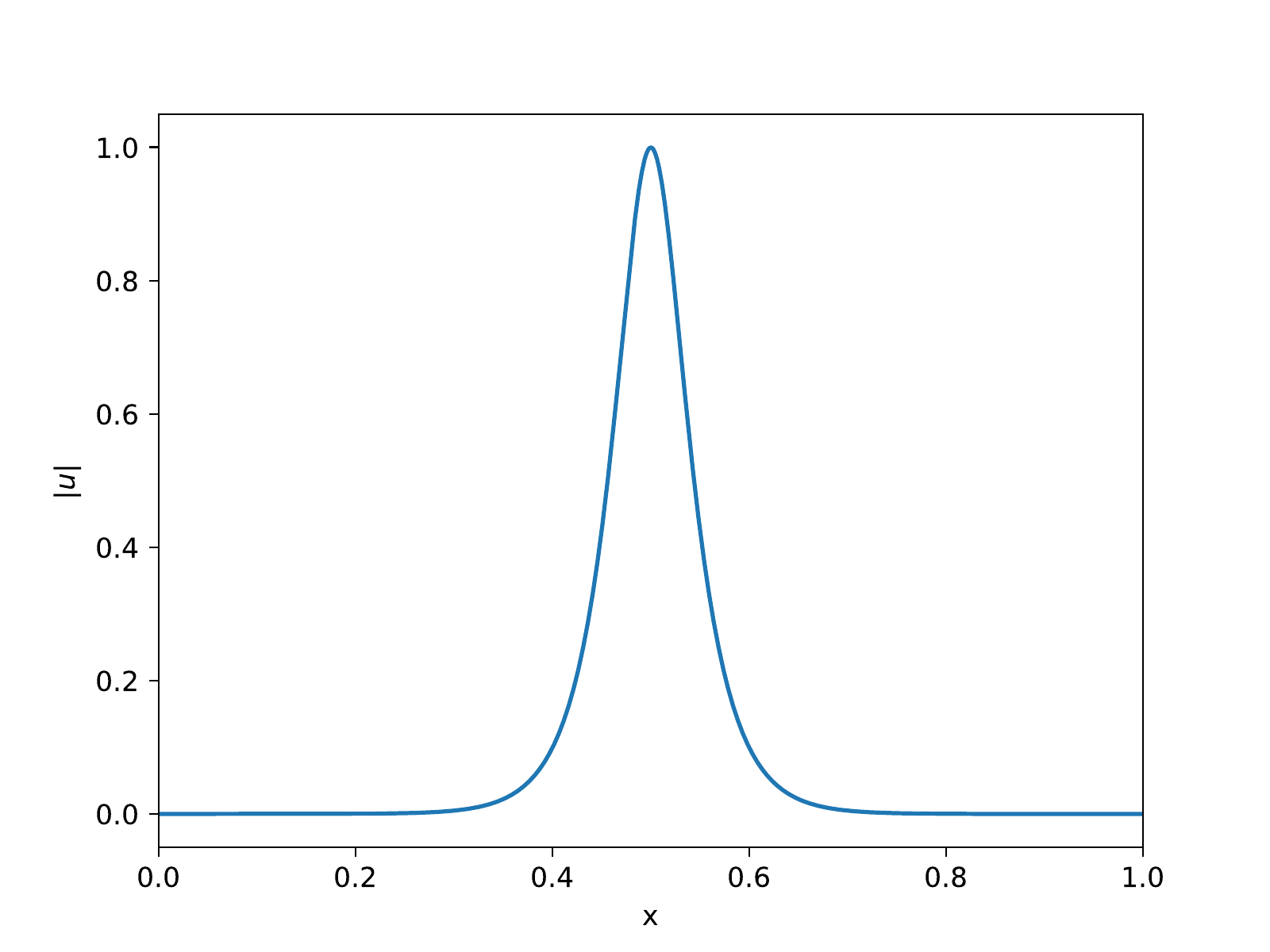}
		\end{minipage}
	}
	\subfigure[active elements at $t=0$]{
		\begin{minipage}[b]{0.46\textwidth}    
			\includegraphics[width=1\textwidth]{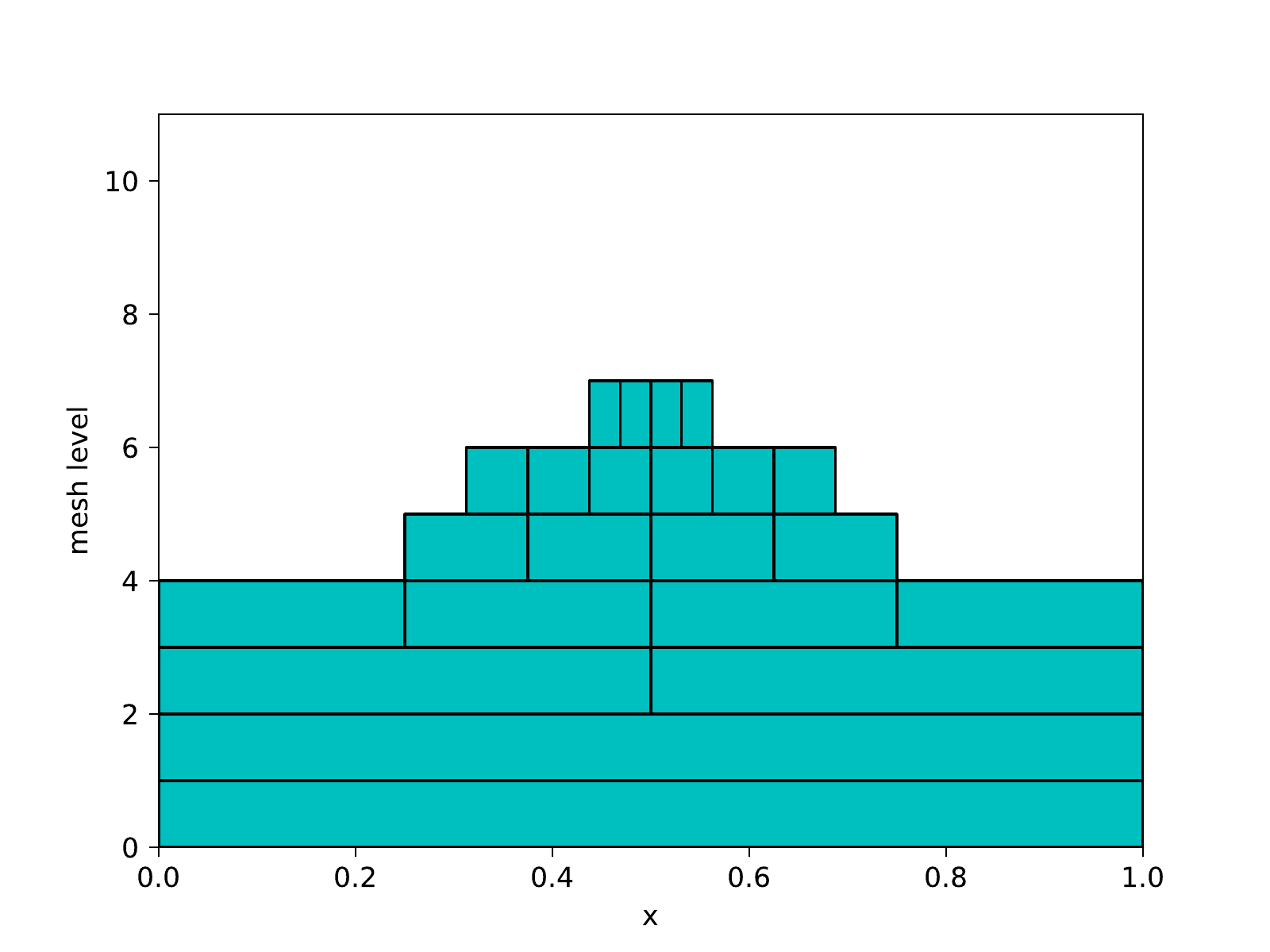}
		\end{minipage}
	}
	\bigskip
	\subfigure[numerical solution at $t=0.4$]{
		\begin{minipage}[b]{0.46\textwidth}
			\includegraphics[width=1\textwidth]{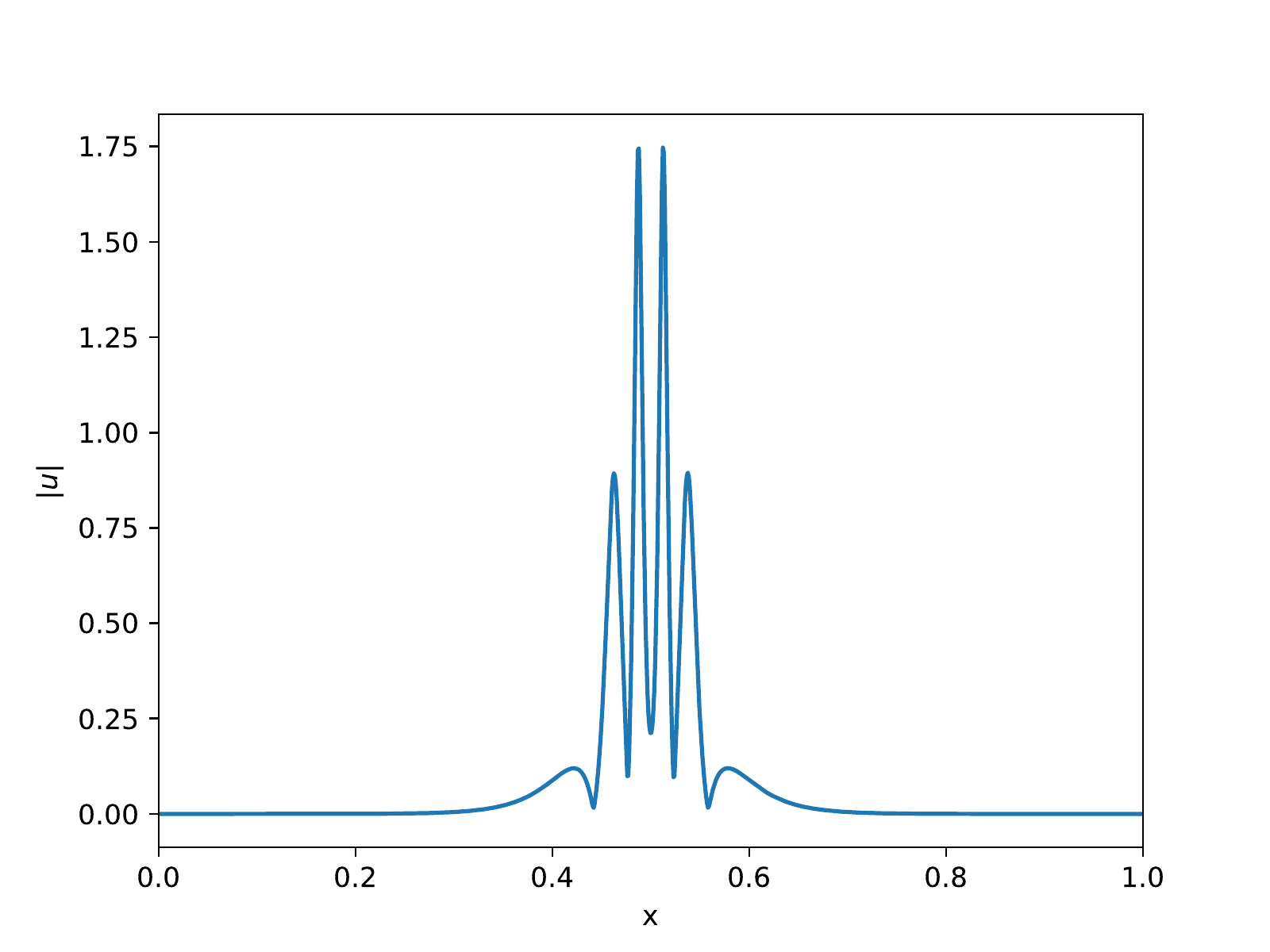}
		\end{minipage}
	}
	\subfigure[active elements at $t=0.4$]{
		\begin{minipage}[b]{0.46\textwidth}    
			\includegraphics[width=1\textwidth]{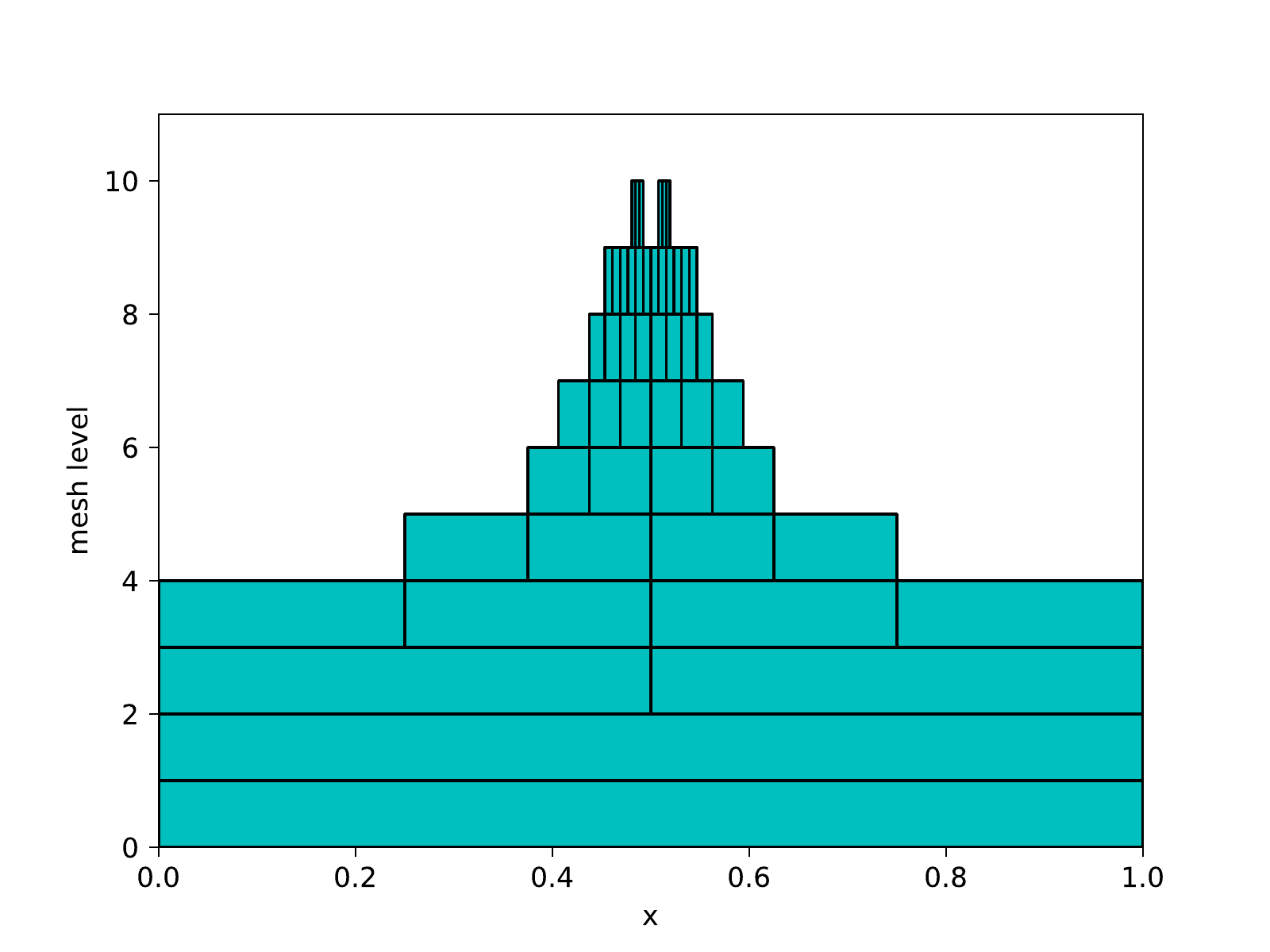}
		\end{minipage}
	}
	\bigskip
	\subfigure[numerical solution at $t=0.6$]{
		\begin{minipage}[b]{0.46\textwidth}
			\includegraphics[width=1\textwidth]{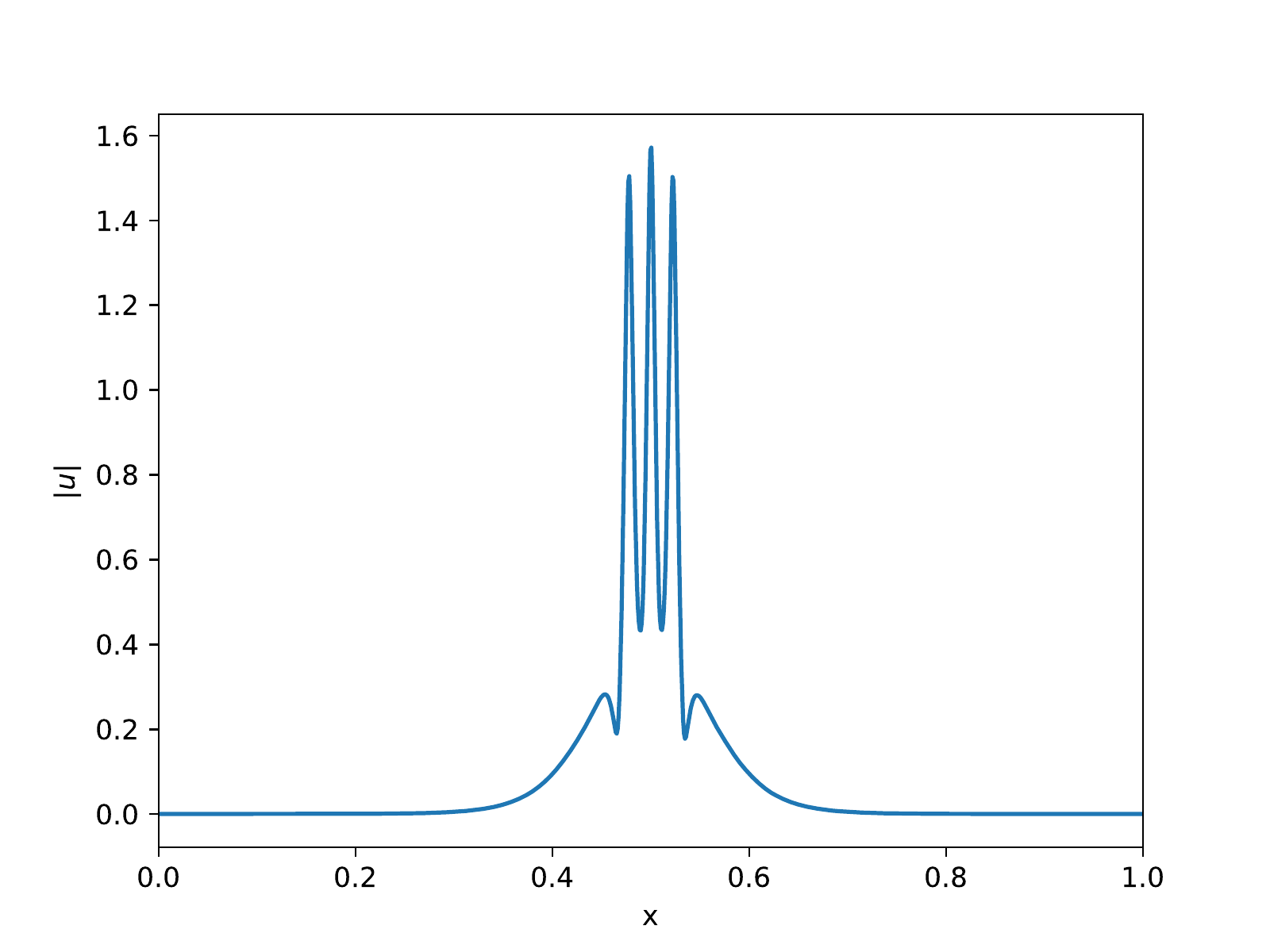}
		\end{minipage}
	}
	\subfigure[active elements at $t=0.6$]{
		\begin{minipage}[b]{0.46\textwidth}    
			\includegraphics[width=1\textwidth]{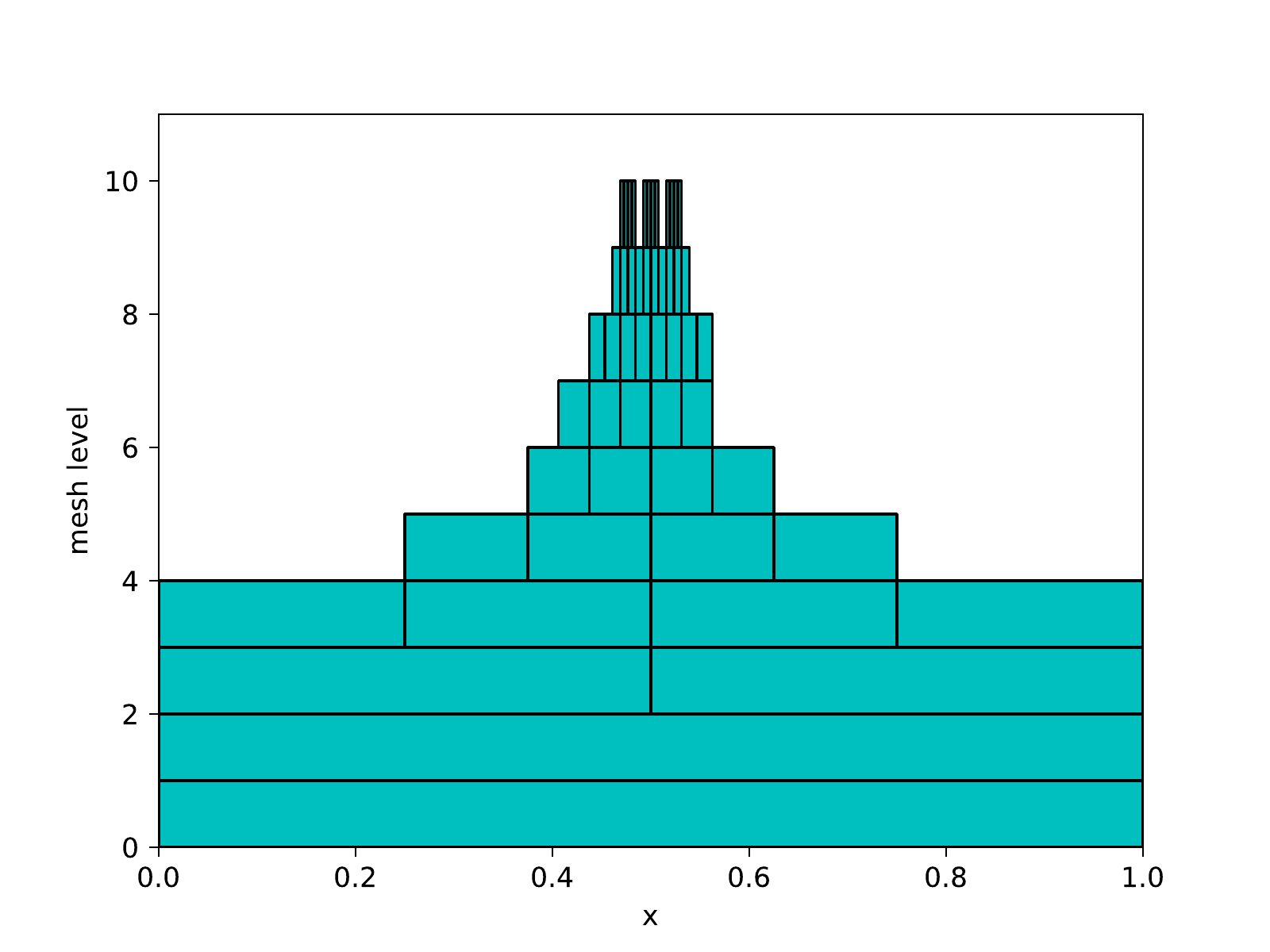}
		\end{minipage}
	}    
	\caption{Example \ref{exam:bound_state}: Bound state solution of solitons with $L=5$.  Left: numerical solutions; right: active elements. $t=0$, 0.4 and 0.6. $N=10, k=3$, $\epsilon=10^{-4}$ and $\eta=10^{-5}$.}
	\label{fig:bound_state_L5}
\end{figure}

\end{exam}

\subsection{Coupled NLS equation in 1D}

% \begin{exam}[Accuracy test for the coupled NLS equation]\label{exam:2-coupled-single}
\begin{exam}\label{exam:2-coupled-single}
We show an accuracy test for the coupled NLS equation \cite{xu2005schrodinger}
\begin{align}
\label{eq:2-coupled}
\left\{ \begin{array}{ll} 
i u_t + i  \frac{\alpha}{M} u_x + \frac{1}{2M^2} u_{xx} + (\abs{u}^2 + \beta \abs{v}^2)u = 0, \\
i v_t - i \frac{\alpha}{M} v_x + \frac{1}{2M^2} v_{xx}  + (\beta \abs{u}^2 + \abs{v}^2)v = 0,
\end{array}\right. 
\end{align}
with the soliton solution
\begin{align}
\left\{ \begin{array}{ll} 
u(x,t) = \sqrt{\frac{2a}{1+\beta} } \sech \left(\sqrt{2a}(X-c t)\right) \exp \left(i\left( (c-\alpha )X - \left( \frac{c^2 - \alpha^2}{2} - a\right) t \right)\right), \\
v(x,t) = \sqrt{\frac{2a}{1+\beta} } \sech \left(\sqrt{2a}(X-c t)\right) \exp \left(i\left( (c+\alpha )X - \left( \frac{c^2 - \alpha^2}{2} - a\right) t \right)\right),
\end{array}\right.
\end{align}
where $ c=1, a=1, \alpha = \frac{1}{2}, \beta = \frac{2}{3}$ and $X= M(x-0.5), M=50$. Periodic boundary condition is applied in $[0,1]$. The solutions are computed up to $t=1$. We take $\Delta t= \frac{0.1 M \Delta x}{\alpha}$, the maximum mesh level $N=10$, and $\eta = \epsilon/10$. The accuracy results are shown in table \ref{tb:2-coupled-single}. We can observe that approximation with higher polynomial degree outperforms that with lower one.
Note that the method  has saturated when $\epsilon = 10^{-4}$ for $k=1$,  therefore the error does not decay too much.

\begin{table}[!hbp]
	\centering
	\caption{Example \ref{exam:2-coupled-single}, accuracy test for the coupled NLS equation, $d=1$. adaptive. $t=1$.}
	\label{tb:2-coupled-single}
	\begin{tabular}{c|c|c|c|c|c|c|c|c}
		\hline
		 \multirow{2}[4]{*}{} & \multirow{2}[4]{*}{$\epsilon$} & \multirow{2}[4]{*}{DoF} & \multicolumn{3}{c|}{Real part of $u$} & \multicolumn{3}{c}{Imaginary part of $u$} \\
		\cline{4-9} 
		&  &  & L$^2$-error & $R_{\textrm{DoF}}$ & $R_{\epsilon}$ & L$^2$-error & $R_{\textrm{DoF}}$ & $R_{\epsilon}$\\
		\hline
		\multirow{4}{3em}{$k=1$}
		& 1e-01 & 28 & 4.42e-02 & - & - & 5.49e-02 & - & -  \\ 
		& 1e-02 & 74 & 9.25e-03 & 1.61 & 0.68 & 1.43e-02 & 1.39 & 0.59 \\ 
		& 1e-03 & 152 & 1.30e-03 & 2.72 & 0.85 & 2.04e-03 & 2.70 & 0.85\\ 
		& 1e-04 & 304 & 6.22e-04 & 1.07 & 0.32 & 1.05e-03 & 0.96 & 0.29\\ 
		\hline
		\multirow{4}{3em}{$k=2$}
		& 1e-01 & 36 & 9.02e-03 & - & - & 1.12e-02 & - & - \\
		& 1e-02 & 54 & 1.38e-03 & 4.62 & 0.81 & 1.59e-03 & 4.82 & 0.85 \\
		& 1e-03 & 105 & 1.39e-04 & 3.45 & 1.00 & 1.68e-04 & 3.38 & 0.98 \\
		& 1e-04 & 186 & 2.05e-05 & 3.35 & 0.83 & 1.99e-05 & 3.73 & 0.93 \\
		\hline
		\multirow{4}{3em}{$k=3$}
		& 1e-01 & 44 & 1.26e-02 & - & - & 1.69e-02 & - & - \\
		& 1e-02 & 60 & 6.82e-04 & 9.40 & 1.27 & 1.21e-03 & 8.50 & 1.15\\
		& 1e-03 & 84 & 8.22e-05 & 6.29 & 0.92 & 1.28e-04 & 6.68 & 0.98 \\
		& 1e-04 & 136 & 1.16e-05 & 4.06 & 0.85 & 1.75e-05 & 4.13 & 0.86 \\
		\hline
		\hline
		\multirow{2}[4]{*}{} & \multirow{2}[4]{*}{$\epsilon$} & \multirow{2}[4]{*}{DoF} & \multicolumn{3}{c|}{Real part of $v$} & \multicolumn{3}{c}{Imaginary part of $v$} \\
		\cline{4-9} 
		&  &  & L$^2$-error & $R_{\textrm{DoF}}$ & $R_{\epsilon}$ & L$^2$-error & $R_{\textrm{DoF}}$ & $R_{\epsilon}$\\
		\hline
		\multirow{4}{3em}{$k=1$}
		& 1e-01 & 28 & 1.13e-01 & - & - & 8.66e-02 & - & -  \\ 
		& 1e-02 & 74 & 2.75e-02 & 1.46 & 0.62 & 2.57e-02 & 1.25 & 0.53 \\ 
		& 1e-03 & 152 & 4.60e-03 & 2.49 & 0.78 & 4.34e-03 & 2.47 & 0.77\\ 
		& 1e-04 & 304 & 1.97e-03 & 1.22 & 0.37 & 1.79e-03 & 1.28 & 0.38\\ 
		\hline
		\multirow{4}{3em}{$k=2$}
		& 1e-01 & 36 & 2.07e-02 & - & - & 2.12e-02 & - & - \\
		& 1e-02 & 54 & 3.56e-03 & 4.34 & 0.76 & 3.39e-03 & 4.52 & 0.80 \\
		& 1e-03 & 105 & 3.27e-04 & 3.59 & 1.04 & 3.76e-04 & 3.31 & 0.95 \\
		& 1e-04 & 186 & 4.82e-05 & 3.35 & 0.83& 5.81e-05 & 3.27 & 0.81 \\
		\hline
		\multirow{4}{3em}{$k=3$}
		& 1e-01 & 44 & 2.81e-02 & - & - & 2.50e-02 & - & - \\
		& 1e-02 & 60 & 1.44e-03 & 9.58 & 1.29 & 1.43e-03 & 9.24 & 1.24\\
		& 1e-03 & 84 & 1.49e-04 & 6.75 & 0.99 & 1.68e-04 & 6.36 & 0.93\\
		& 1e-04 & 136 & 2.16e-05 & 4.01 & 0.84 & 3.43e-05 & 3.30 & 0.69 \\
		\hline
	\end{tabular}
\end{table}

\end{exam}

\begin{exam}\label{exam:2-coupled-interaction}

In this example, we consider the solitary wave propagation and the soliton interaction for the coupled NLS equation \eqref{eq:2-coupled} following \cite{xu2005schrodinger}. In this example, $\Delta t$ is taken as $\Delta t= \frac{0.1 M \Delta x}{\alpha}$.

We first take the initial condition for soliton propagation
\begin{align}
\left\{ \begin{array}{ll} 
u(x,0) = \sqrt{\frac{2a}{1+\beta} } \sech \left(\sqrt{2a}X\right) \exp \left(i (c-\alpha )X \right), \\
v(x,0) = \sqrt{\frac{2a}{1+\beta} } \sech \left(\sqrt{2a}X\right) \exp \left(i (c+\alpha )X \right),
\end{array}\right.
\end{align}
with the same parameters in Example \ref{exam:2-coupled-single} except $X= M(x-0.2), M=100$. 
Periodic boundary condition is used in $[0,1]$. 
%We take $N=9, k=3$, $\epsilon=10^{-4}$ and $\eta=5 \times 10^{-5}$.
The numerical solutions and active elements at $t=0, 20$ and 50 are presented in Figure \ref{fig:2-coupled-single}.
The plots of $|u|$ and $|v|$ are similar, thus we only show the results of $|u|$ here. 
%Observe that the soliton travels to the right at velocity $1/M = 0.01$. 

\begin{figure}
	\centering
	\subfigure[numerical solution at $t=0$]{
		\begin{minipage}[b]{0.46\textwidth}
			\includegraphics[width=1\textwidth]{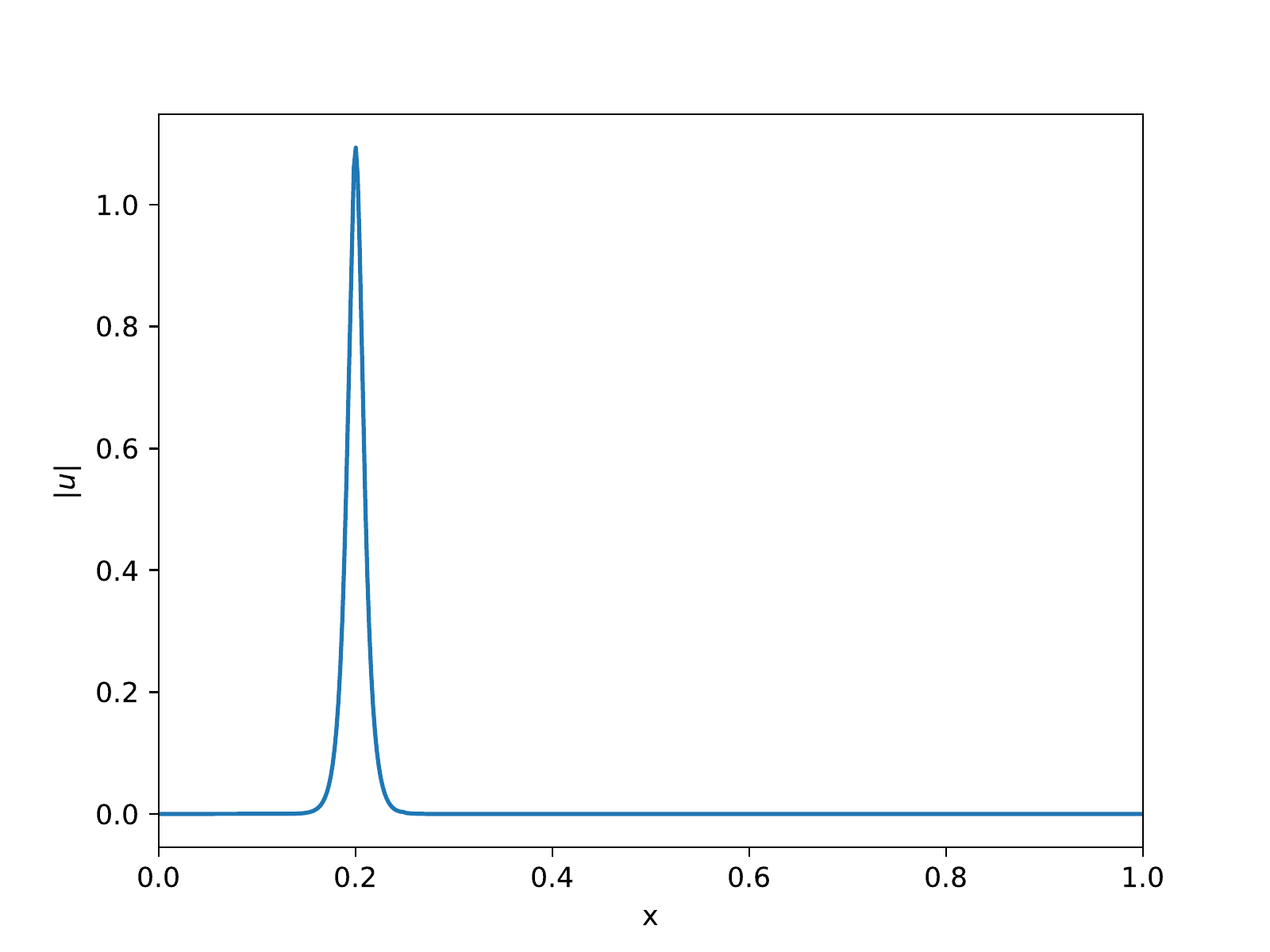}
		\end{minipage}
	}
	\subfigure[active elements at $t=0$]{
		\begin{minipage}[b]{0.46\textwidth}    
			\includegraphics[width=1\textwidth]{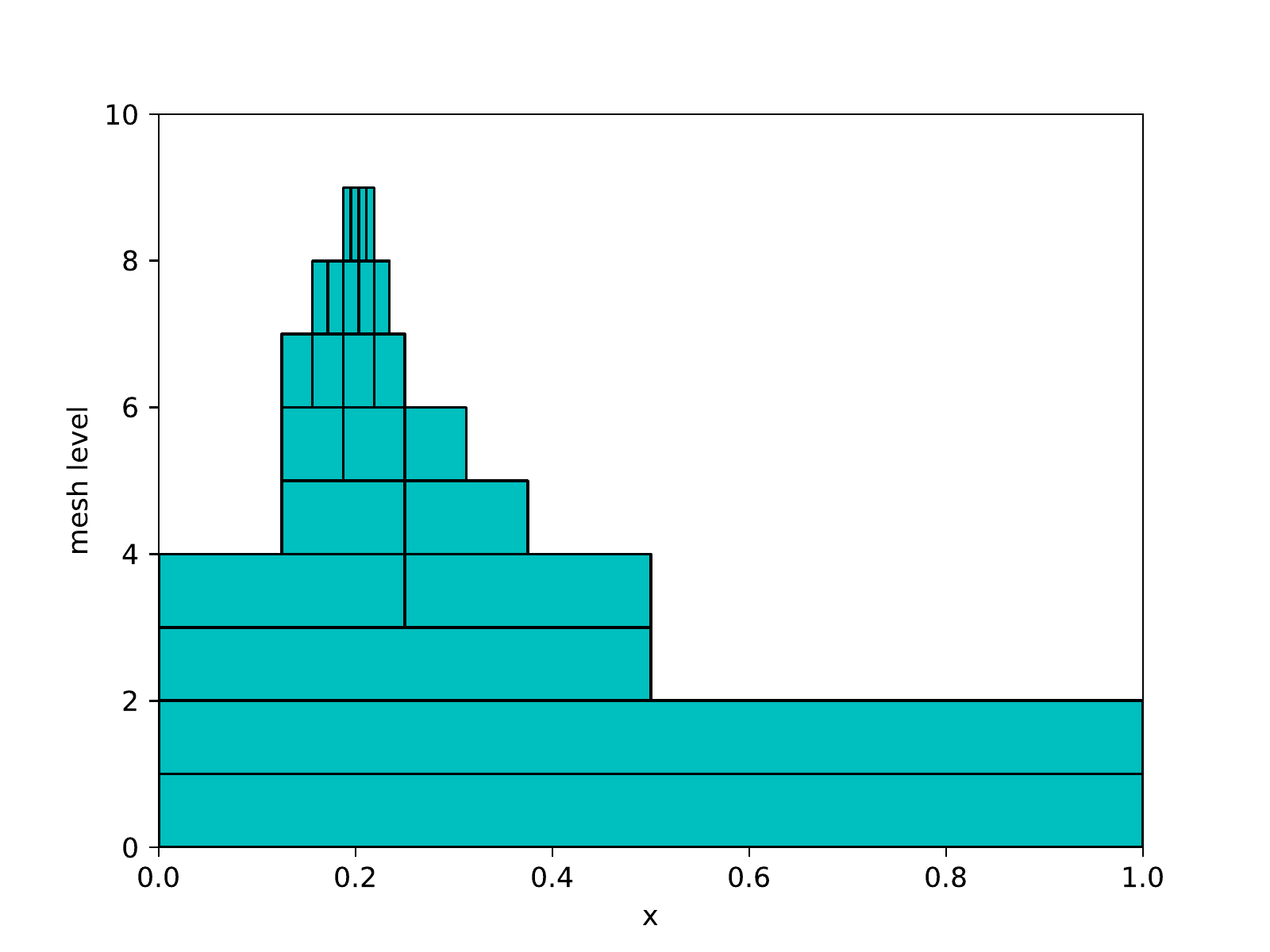}
		\end{minipage}
	}
	\bigskip
	\subfigure[numerical solution at $t=20$]{
		\begin{minipage}[b]{0.46\textwidth}
			\includegraphics[width=1\textwidth]{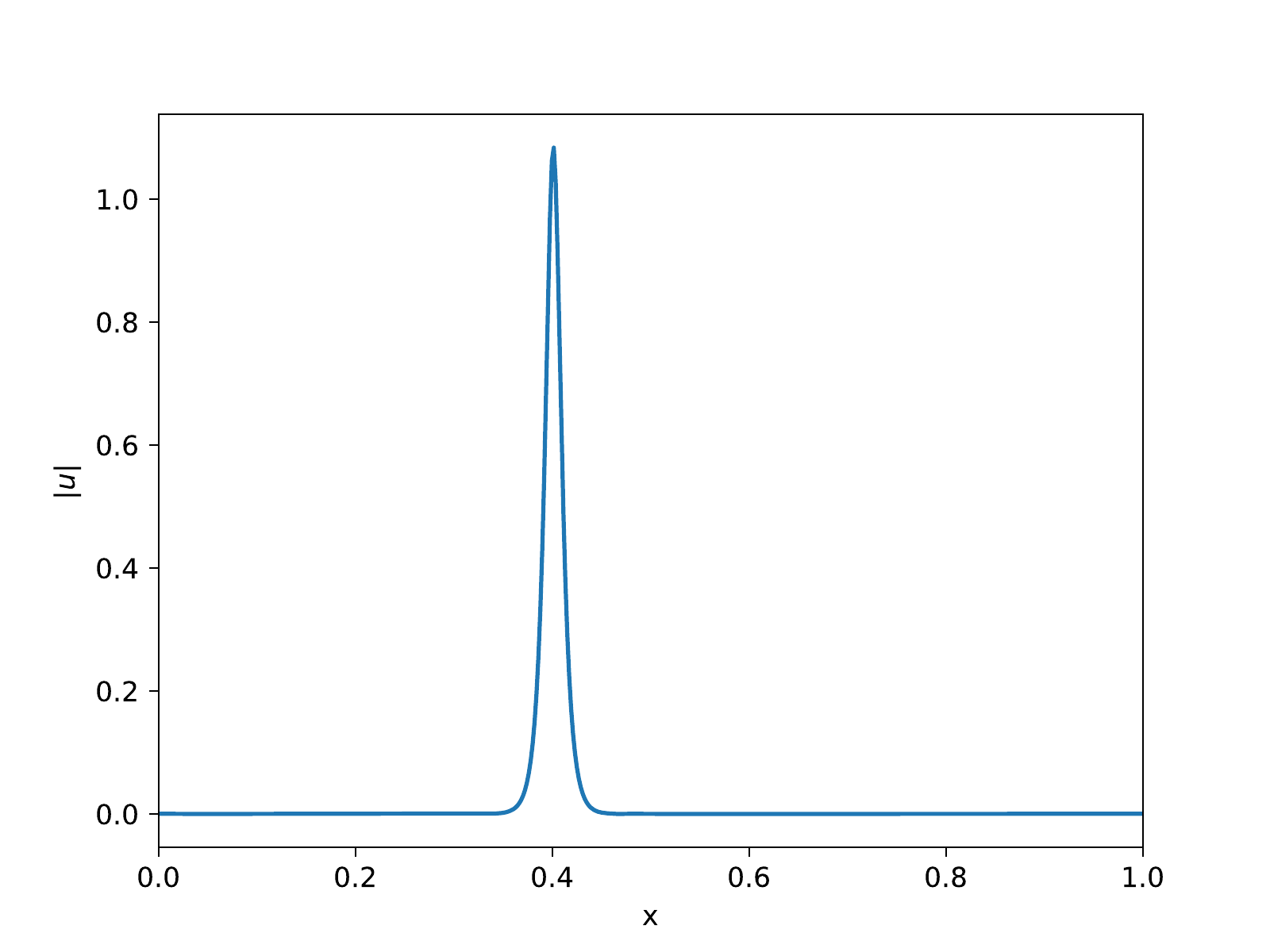}
		\end{minipage}
	}
	\subfigure[active elements at $t=20$]{
		\begin{minipage}[b]{0.46\textwidth}    
			\includegraphics[width=1\textwidth]{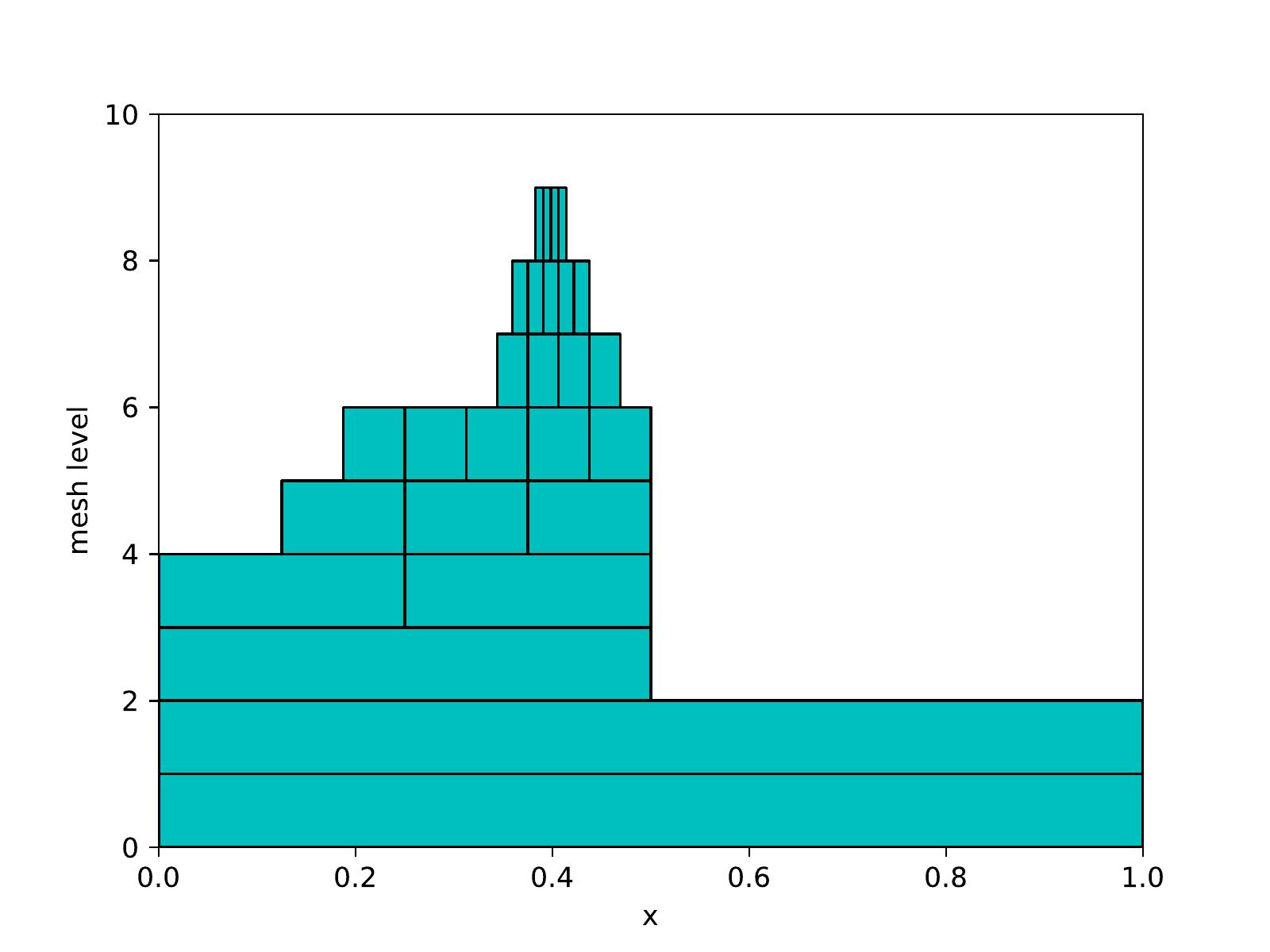}
		\end{minipage}
	}
    \bigskip
	\subfigure[numerical solution at $t=50$]{
		\begin{minipage}[b]{0.46\textwidth}
			\includegraphics[width=1\textwidth]{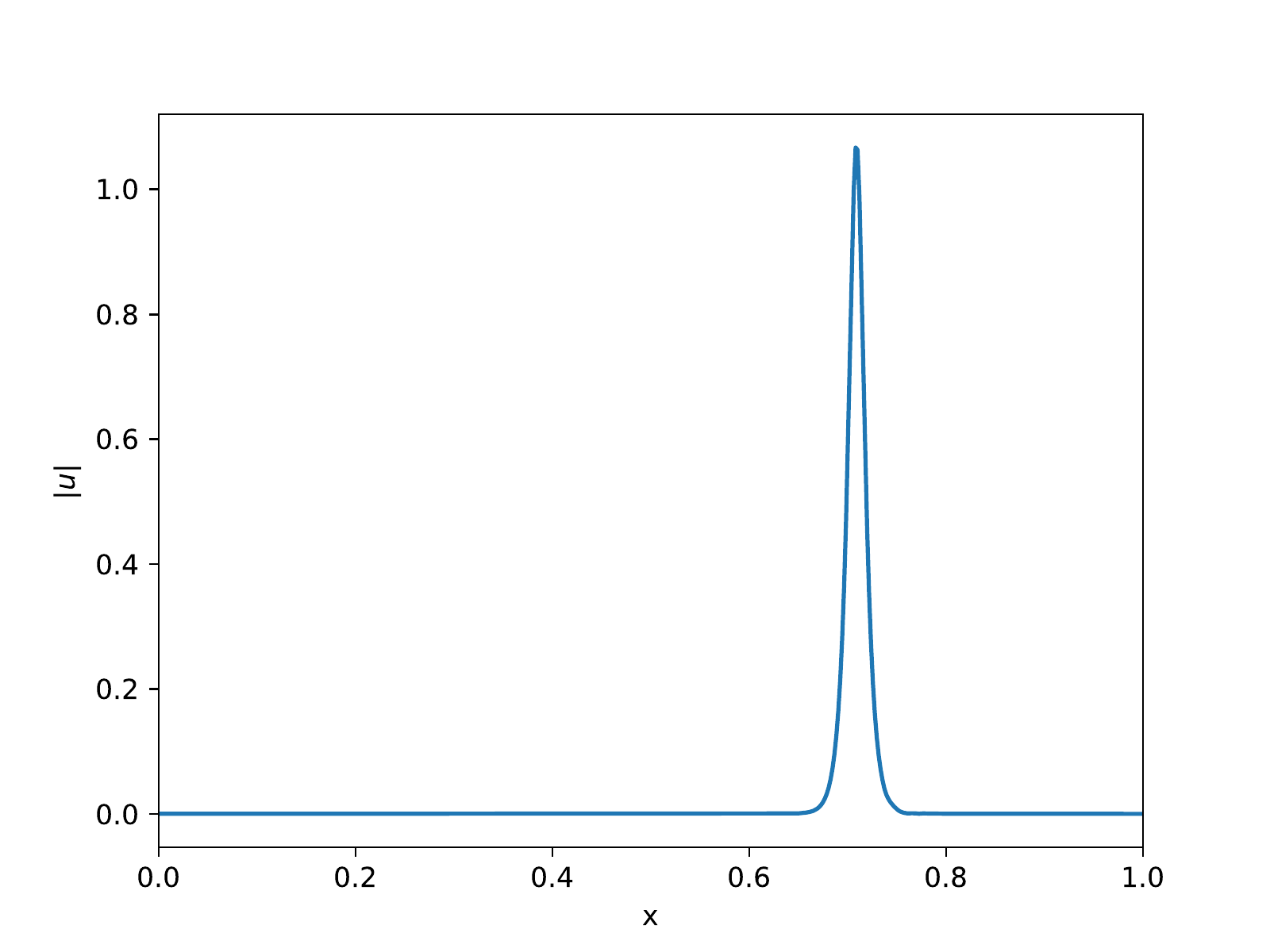}
		\end{minipage}
	}
	\subfigure[active elements at $t=50$]{
		\begin{minipage}[b]{0.46\textwidth}    
			\includegraphics[width=1\textwidth]{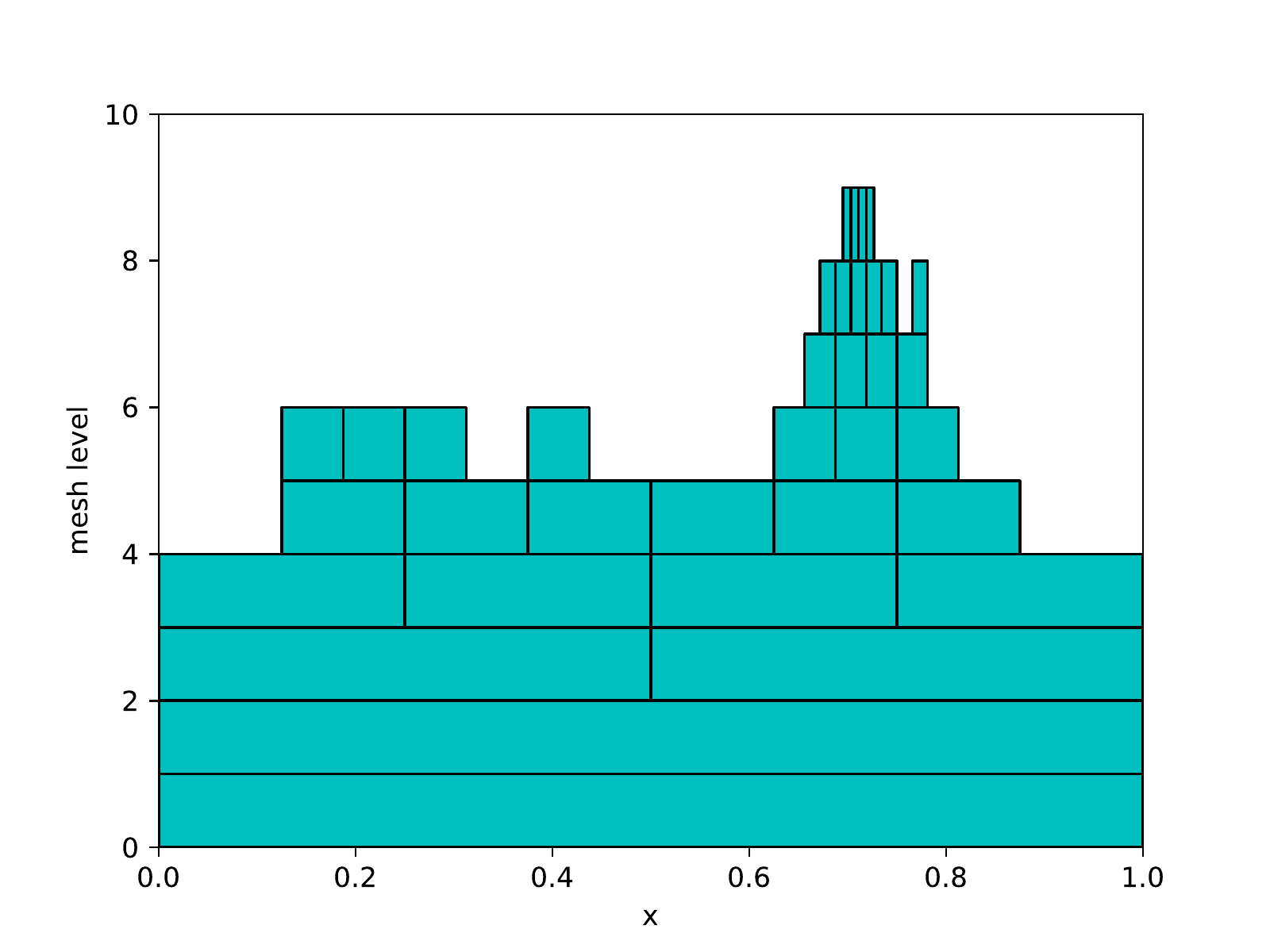}
		\end{minipage}
	}    
	\caption{Example \ref{exam:2-coupled-single}: coupled NLS equation, single  soliton.  Left: numerical solutions; right: active elements. $t=0$, 20 and 50. $N=9, k=3$, $\epsilon=10^{-4}$ and $\eta=5 \times 10^{-5}$.}
	\label{fig:2-coupled-single}
\end{figure}

For interaction of two solitons, we use the following initial condition
\begin{align}
\left\{ \begin{array}{ll} 
u(x,0) = \sum_{j=1}^{2} \sqrt{\frac{2a_j}{1+\beta} } \sech \left(\sqrt{2a_j}X_j\right) \exp \left(i (c_j-\alpha )X_j \right), \\
v(x,0) = \sum_{j=1}^{2} \sqrt{\frac{2a_j}{1+\beta} } \sech \left(\sqrt{2a_j}X_j\right) \exp \left(i (c_j+\alpha )X_j \right),
\end{array}\right.
\end{align}
where $c_1=1, c_2 = 0.1,  a_1=1, a_2 = 0.5, \alpha = \frac{1}{2}, \beta = \frac{2}{3}$ and $X_j= M(x-0.2-x_j), M=100, x_1=0, x_2=0.25$. Periodic boundary condition is used in $[0,1]$. 
%We take $N=9, k=3$, $\epsilon=10^{-4}$ and $\eta=4 \times 10^{-5}$.
The numerical solutions of $|u|$  and active elements at $t=0, 20$ and 50 are presented in Figure \ref{fig:2-coupled-double}.
%The plots of $|u|$ and $|v|$ are similar for this problem, thus we only show the results of $|u|$ here. %We observe from the figure that the two solitons travel to the right at velocities 1 and 0.1 and they interact between $t=20$ and $t=30$. 
The interaction is elastic and the solitons restore their original shapes. 

Next, we consider interaction of three solitons with initial condition
\begin{align}
\left\{ \begin{array}{ll} 
u(x,0) = \sum_{j=1}^{3} \sqrt{\frac{2a_j}{1+\beta} } \sech \left(\sqrt{2a_j}X_j\right) \exp \left(i (c_j-\alpha )X_j \right), \\
v(x,0) = \sum_{j=1}^{3} \sqrt{\frac{2a_j}{1+\beta} } \sech \left(\sqrt{2a_j}X_j\right) \exp \left(i (c_j+\alpha )X_j \right),
\end{array}\right.
\end{align}
where $c_1=1, c_2 = 0.1, c_3 = -1, a_1=1.2, a_2 = 0.72, a_3 = 0.36, \alpha = \frac{1}{2}, \beta = \frac{2}{3}$ and $X_j= M(x-0.2-x_j), M=100, x_1=0, x_2=0.25, x_3 = 0.5$. Periodic boundary condition is used in $[0,1]$. 
%We take $N=9, k=3$, $\epsilon=10^{-4}$ and $\eta=4 \times 10^{-5}$.
The numerical solutions of $|u|$ and active elements at $t=0, 20$ and 50 are presented in Figure \ref{fig:2-coupled-triple}.
%Notice that the three solitons travel at velocities 1, 0.1 and -1 along x direction. They interact between  $t=20$ and $t=30$ and restore their original shapes.
Notice that the three solitons restore their original shapes after interaction.

\begin{figure}
	\centering
	\subfigure[numerical solution at $t=0$]{
		\begin{minipage}[b]{0.46\textwidth}
			\includegraphics[width=1\textwidth]{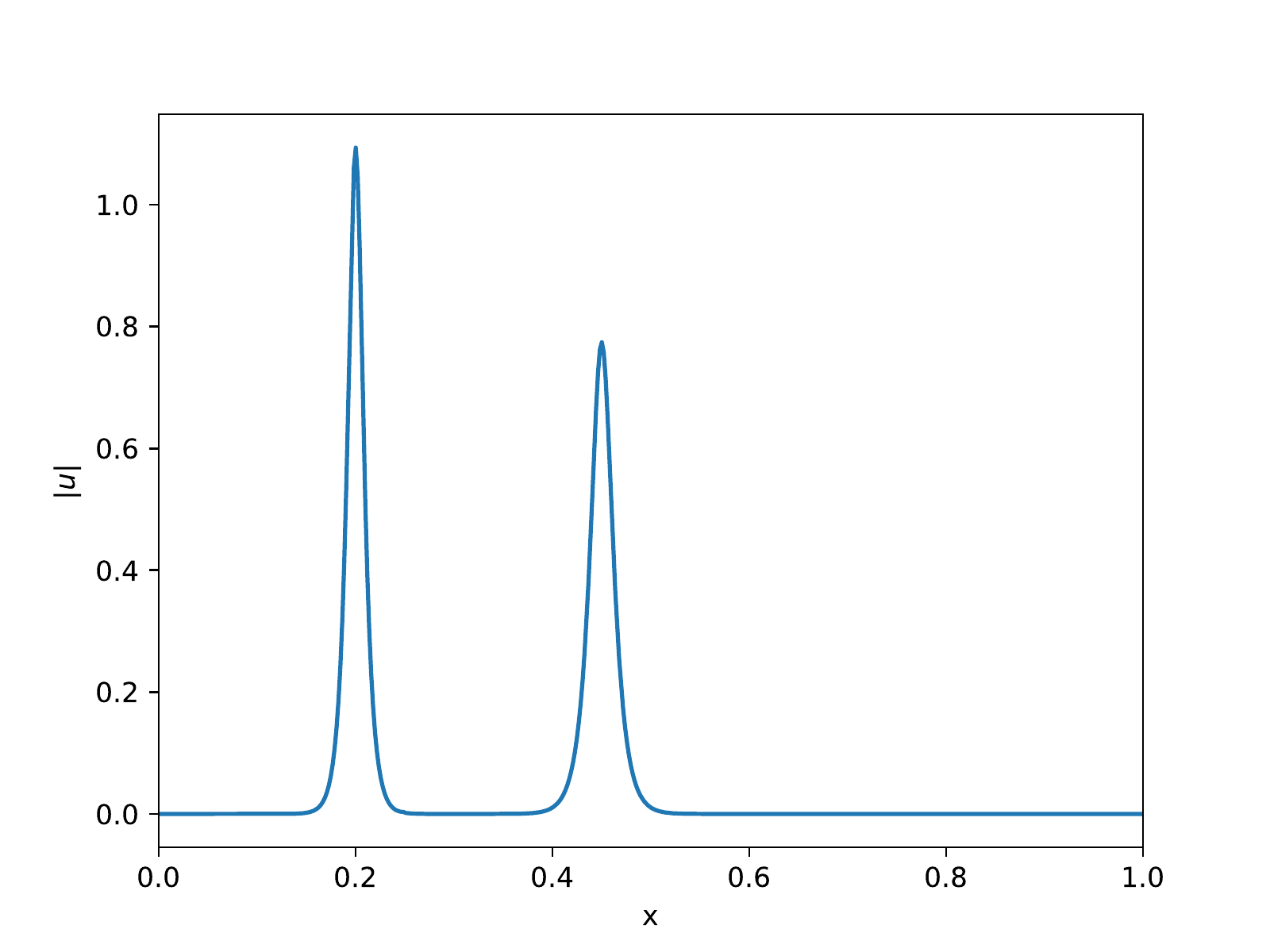}
		\end{minipage}
	}
	\subfigure[active elements at $t=0$]{
		\begin{minipage}[b]{0.46\textwidth}    
			\includegraphics[width=1\textwidth]{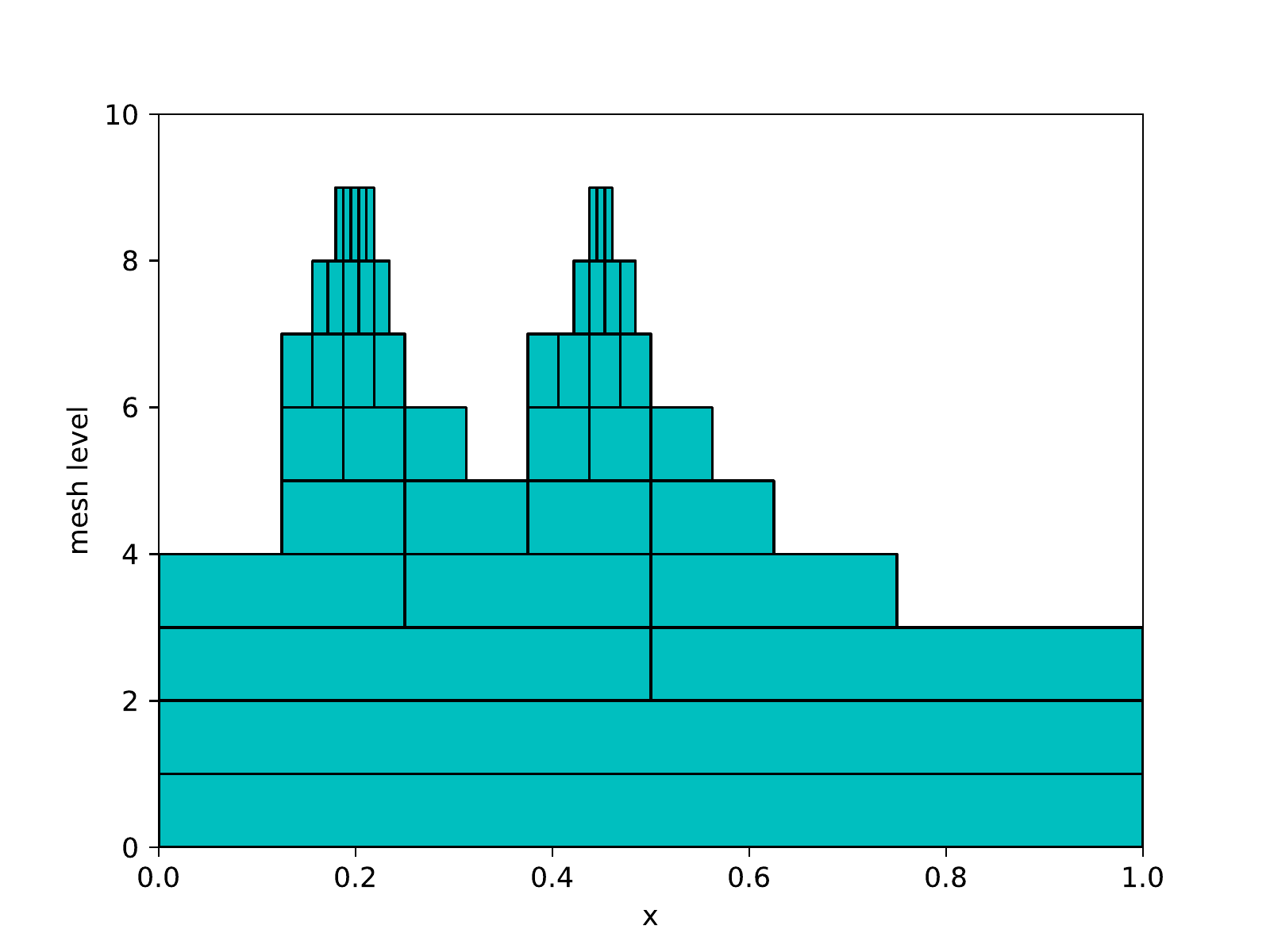}
		\end{minipage}
	}
	\bigskip
	\subfigure[numerical solution at $t=20$]{
		\begin{minipage}[b]{0.46\textwidth}
			\includegraphics[width=1\textwidth]{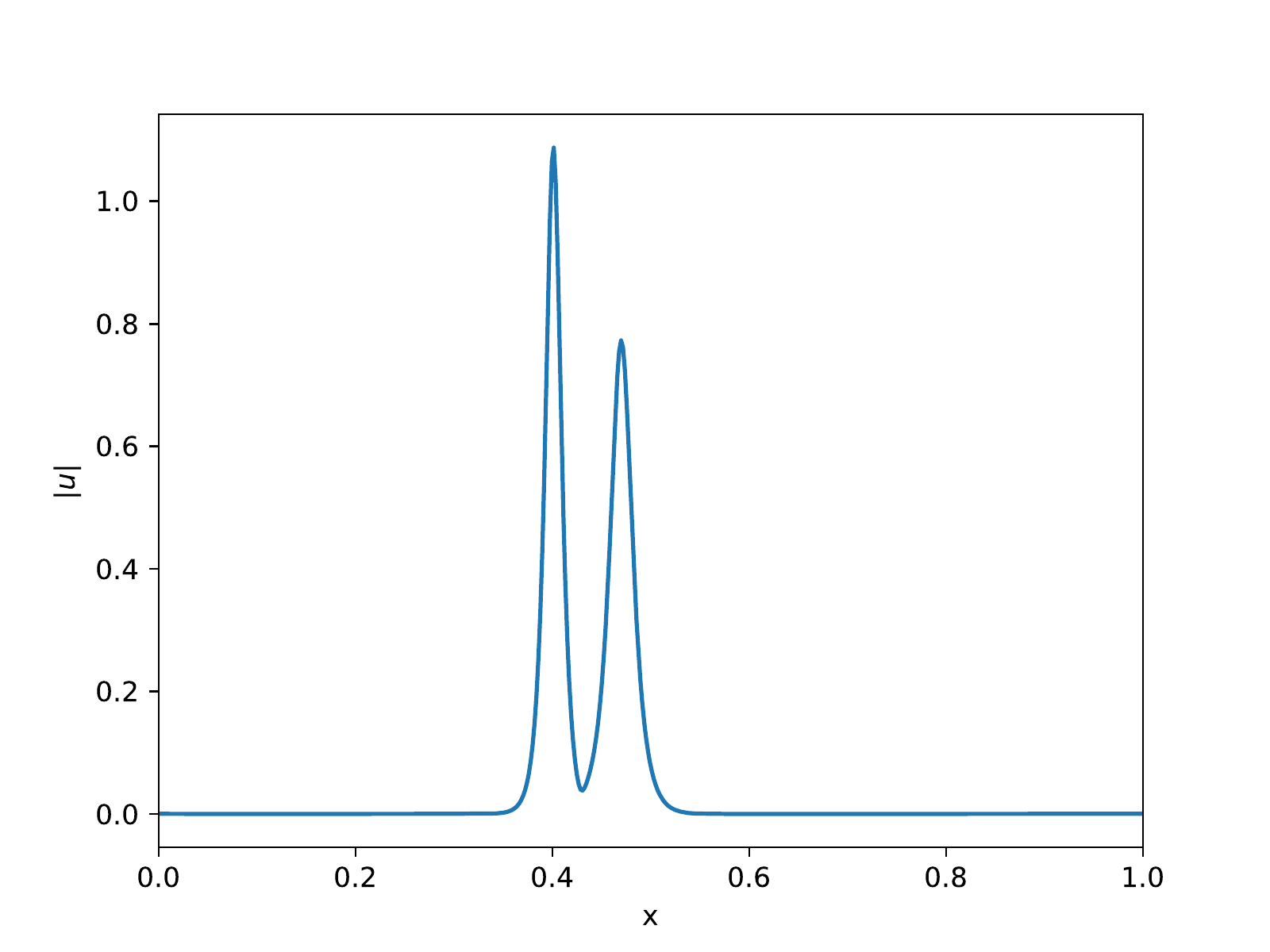}
		\end{minipage}
	}
	\subfigure[active elements at $t=20$]{
		\begin{minipage}[b]{0.46\textwidth}    
			\includegraphics[width=1\textwidth]{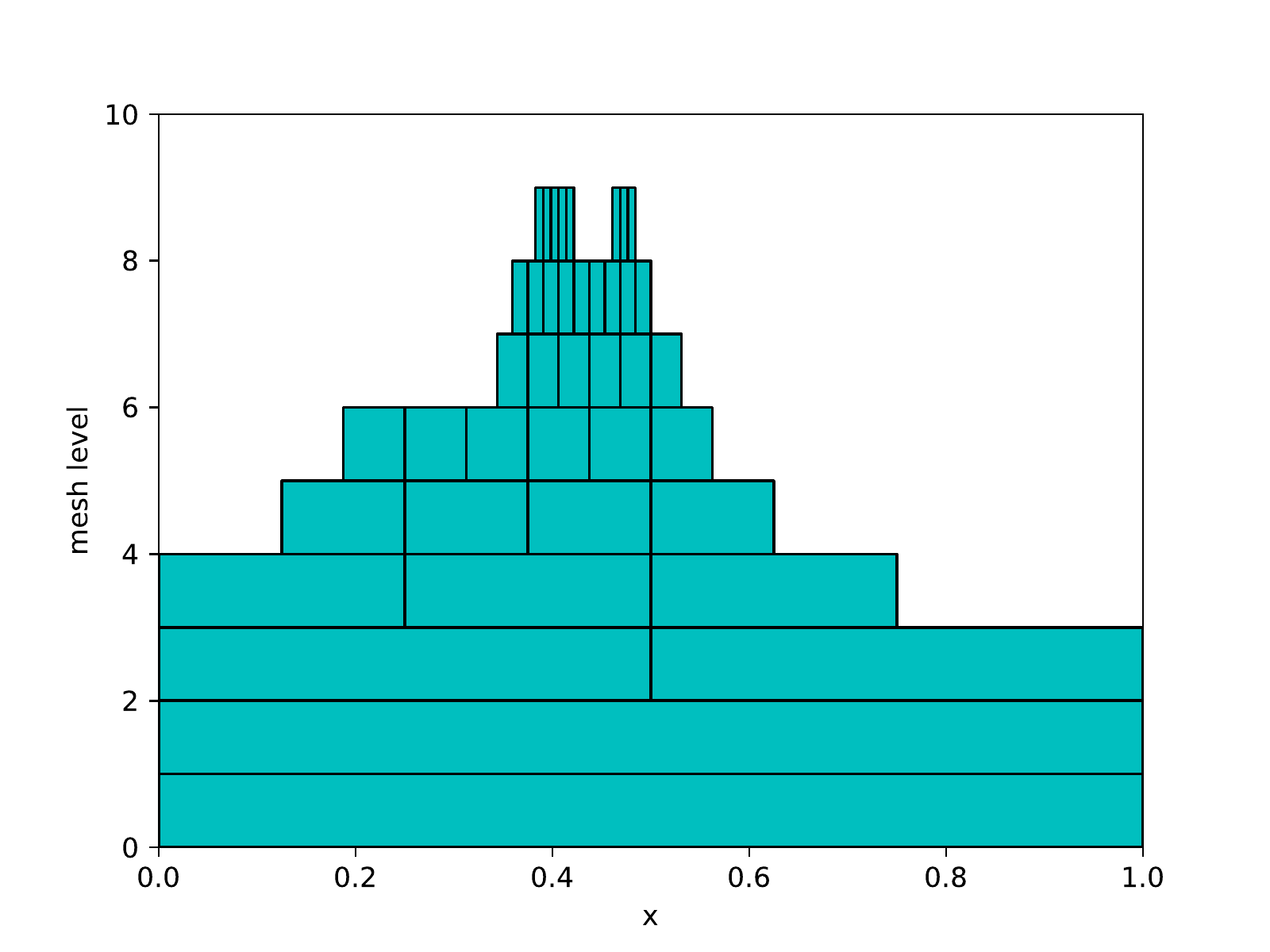}
		\end{minipage}
	}
	\bigskip
	\subfigure[numerical solution at $t=50$]{
		\begin{minipage}[b]{0.46\textwidth}
			\includegraphics[width=1\textwidth]{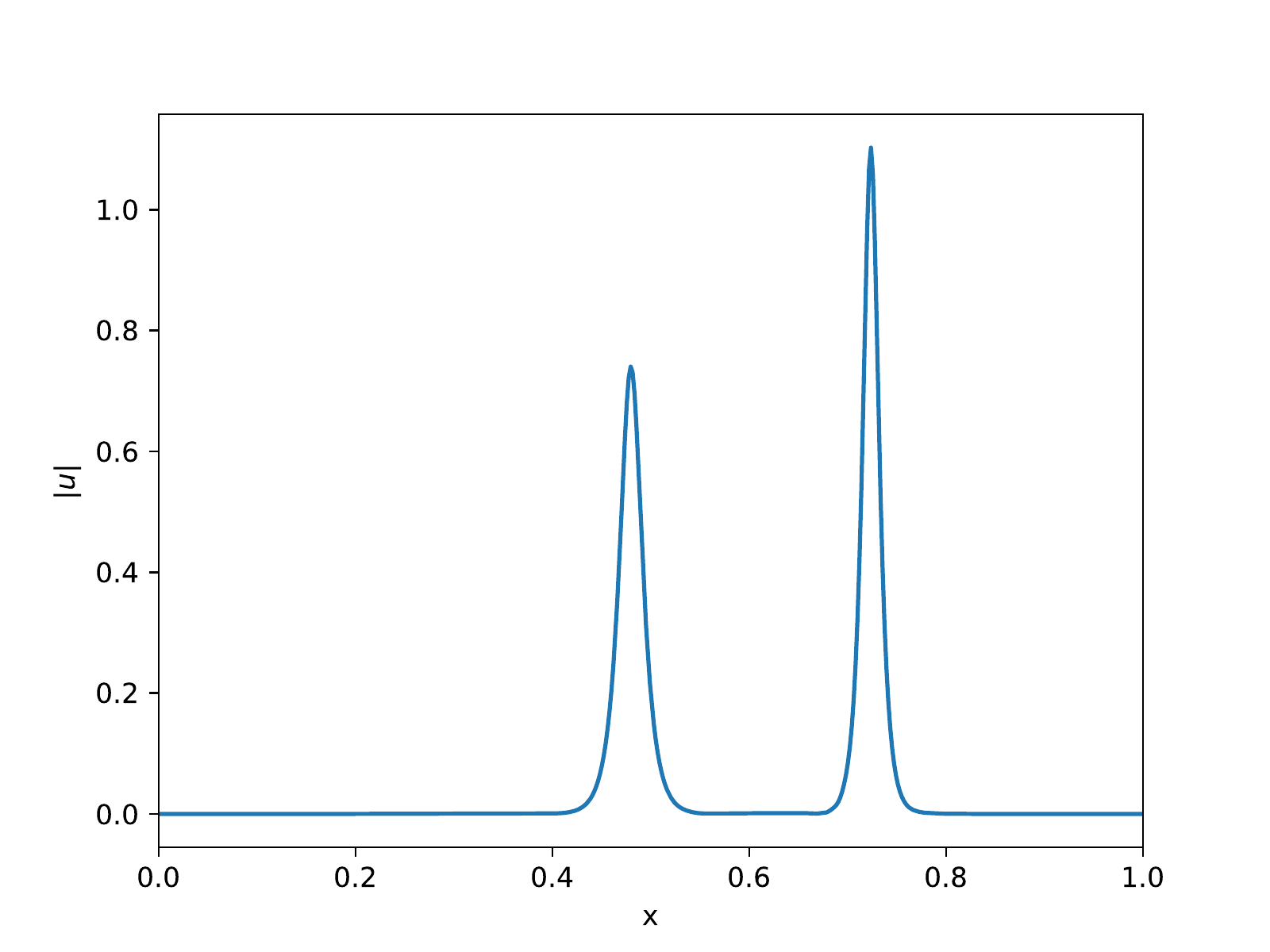}
		\end{minipage}
	}
	\subfigure[active elements at $t=50$]{
		\begin{minipage}[b]{0.46\textwidth}    
			\includegraphics[width=1\textwidth]{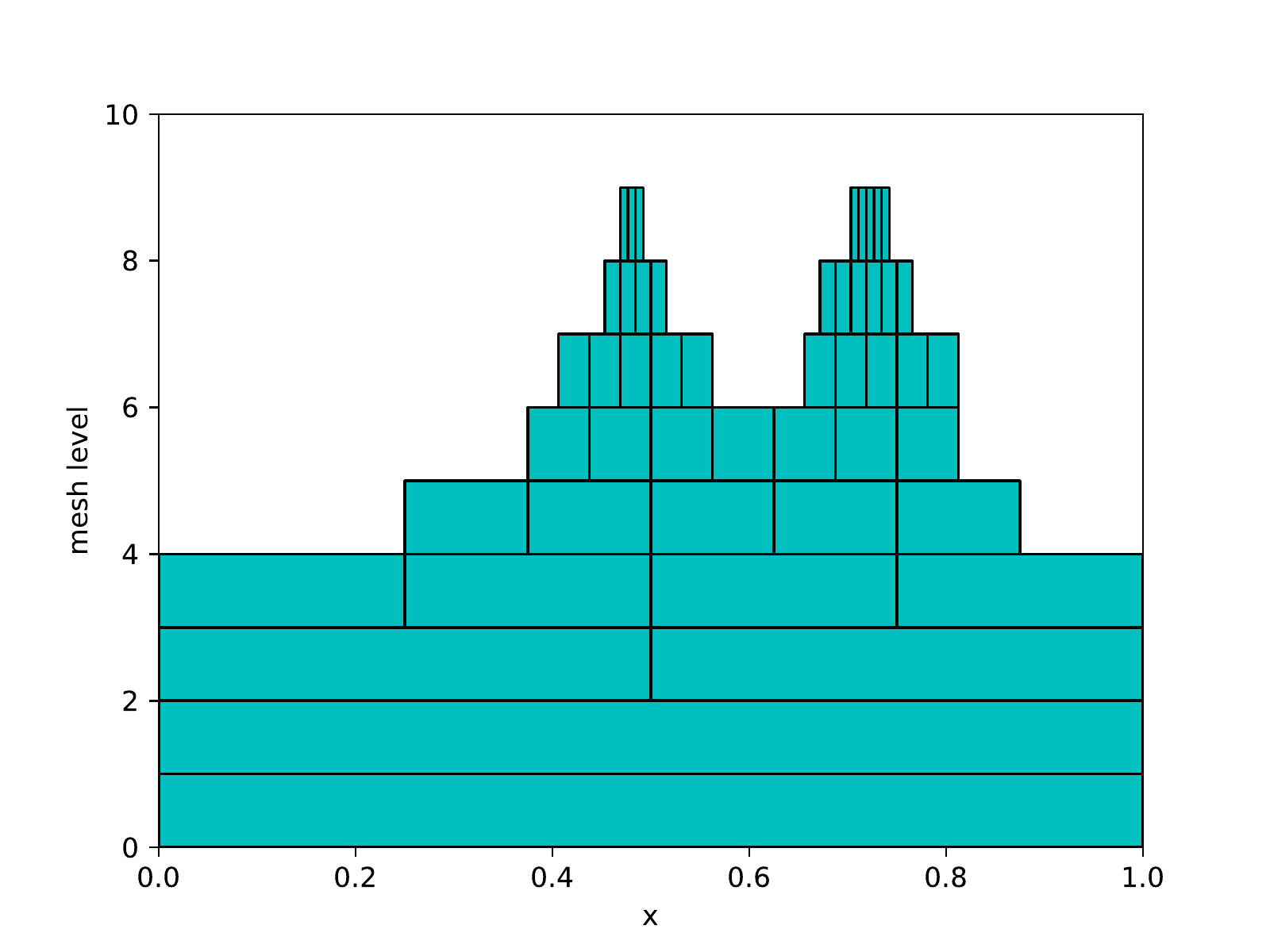}
		\end{minipage}
	}    
	\caption{Example \ref{exam:2-coupled-interaction}: coupled NLS equation, double  solitons.  Left: numerical solutions; right: active elements. $t=0$, 20 and 50. $N=9, k=3$, $\epsilon=10^{-4}$ and $\eta=4 \times 10^{-5}$. }
	\label{fig:2-coupled-double}
\end{figure}

\begin{figure}
	\centering
	\subfigure[numerical solution at $t=0$]{
		\begin{minipage}[b]{0.46\textwidth}
			\includegraphics[width=1\textwidth]{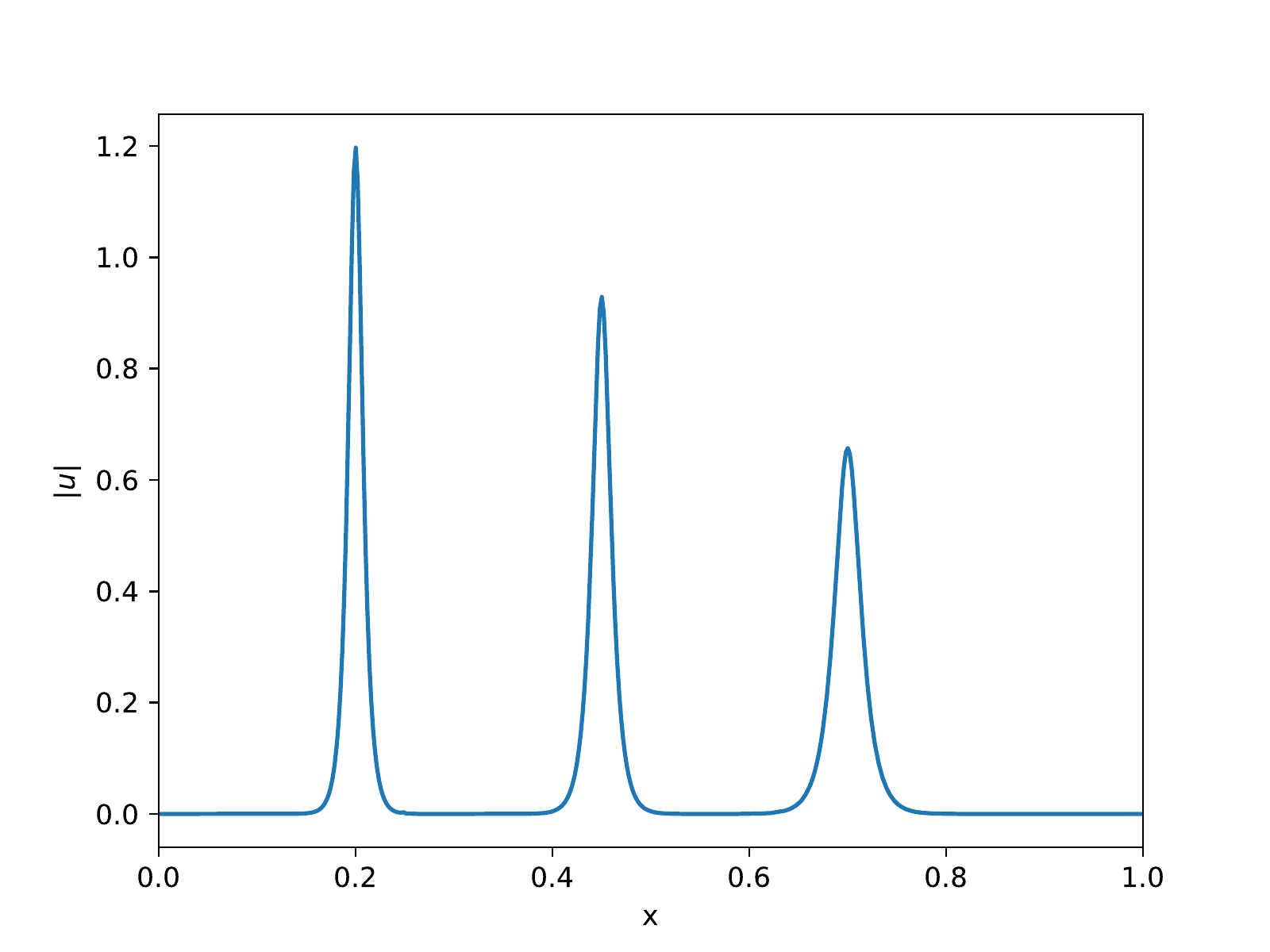}
		\end{minipage}
	}
	\subfigure[active elements at $t=0$]{
		\begin{minipage}[b]{0.46\textwidth}    
			\includegraphics[width=1\textwidth]{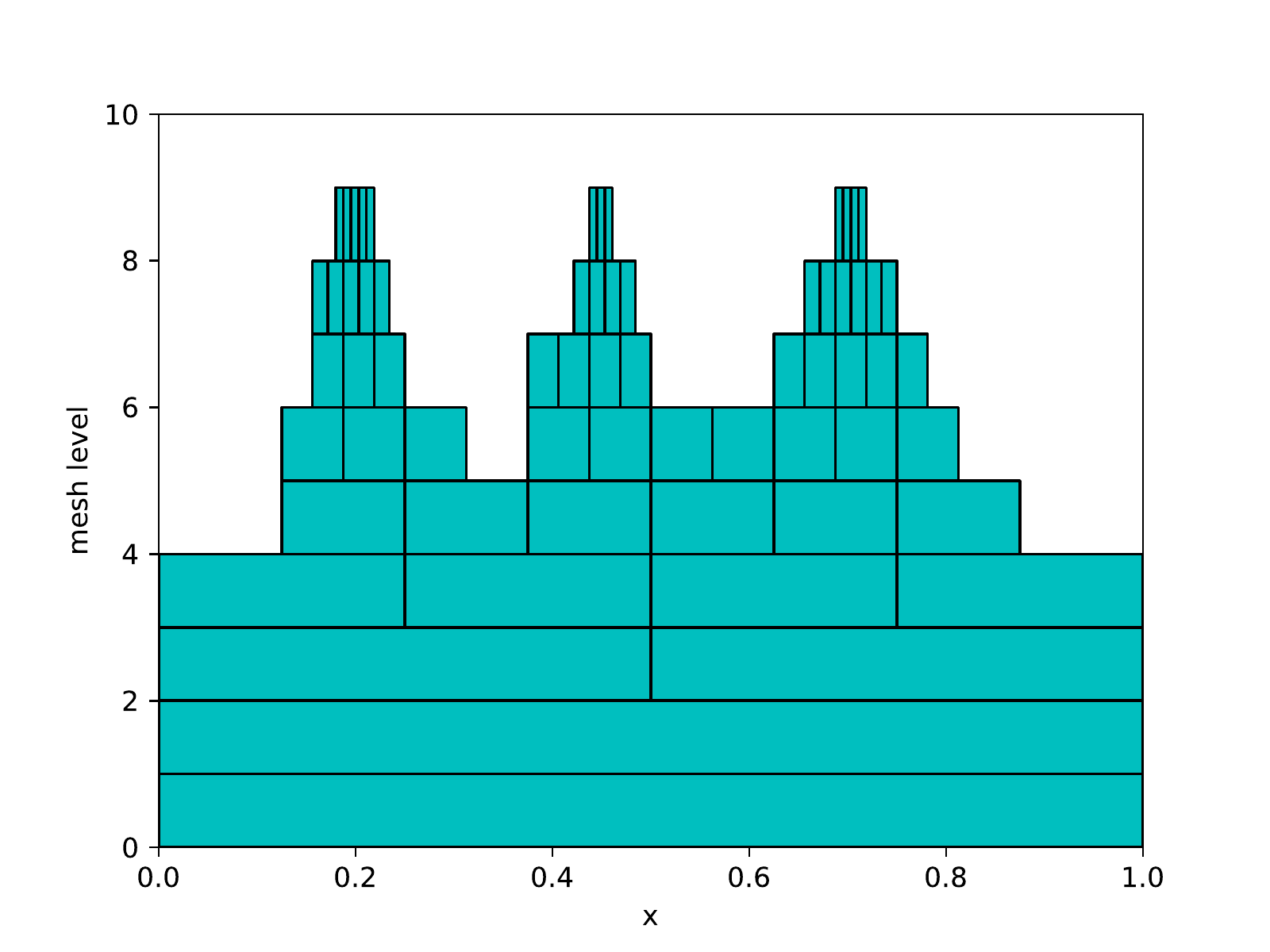}
		\end{minipage}
	}
	\bigskip
	\subfigure[numerical solution at $t=20$]{
		\begin{minipage}[b]{0.46\textwidth}
			\includegraphics[width=1\textwidth]{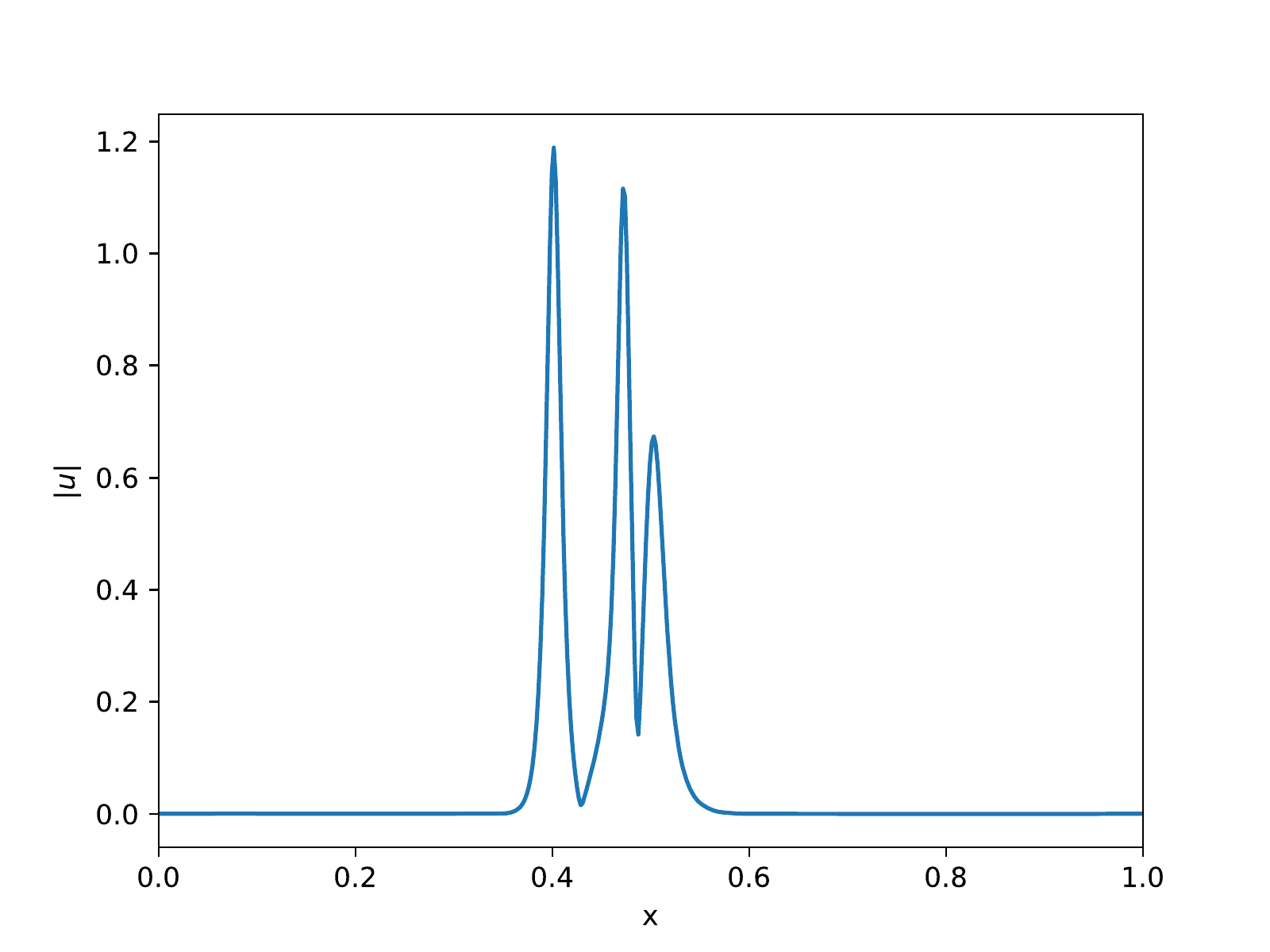}
		\end{minipage}
	}
	\subfigure[active elements at $t=20$]{
		\begin{minipage}[b]{0.46\textwidth}    
			\includegraphics[width=1\textwidth]{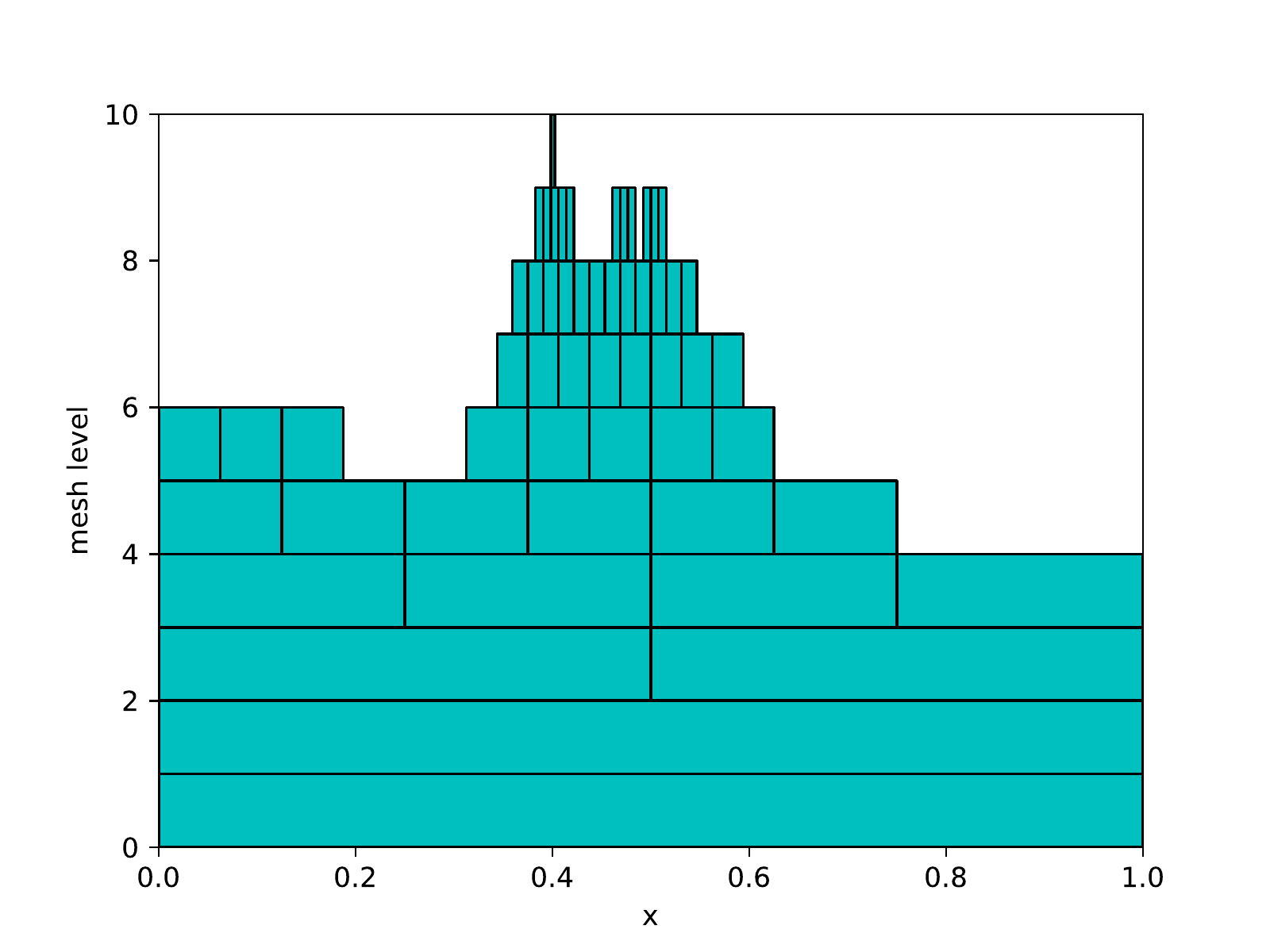}
		\end{minipage}
	}
	\bigskip
	\subfigure[numerical solution at $t=50$]{
		\begin{minipage}[b]{0.46\textwidth}
			\includegraphics[width=1\textwidth]{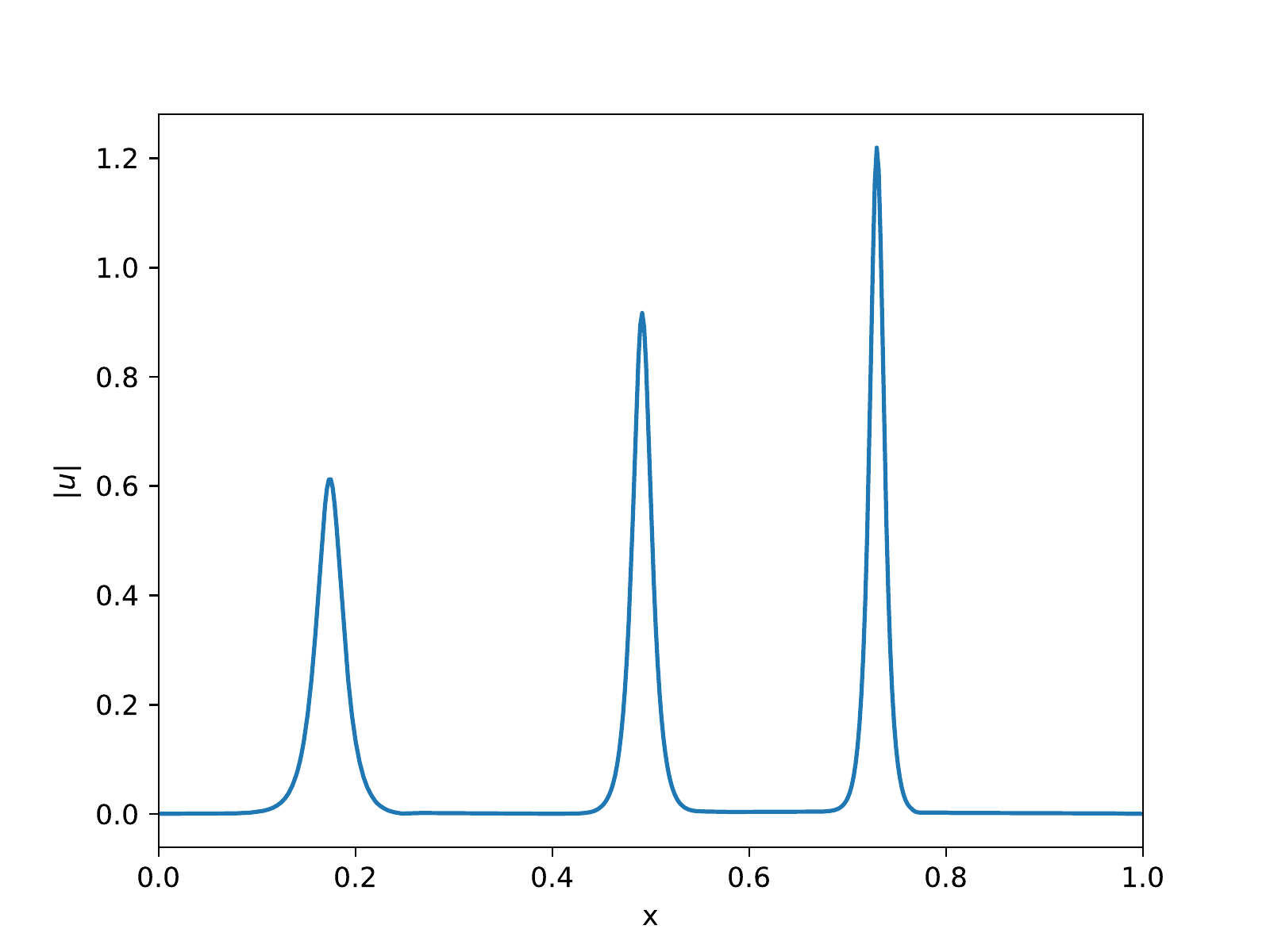}
		\end{minipage}
	}
	\subfigure[active elements at $t=50$]{
		\begin{minipage}[b]{0.46\textwidth}    
			\includegraphics[width=1\textwidth]{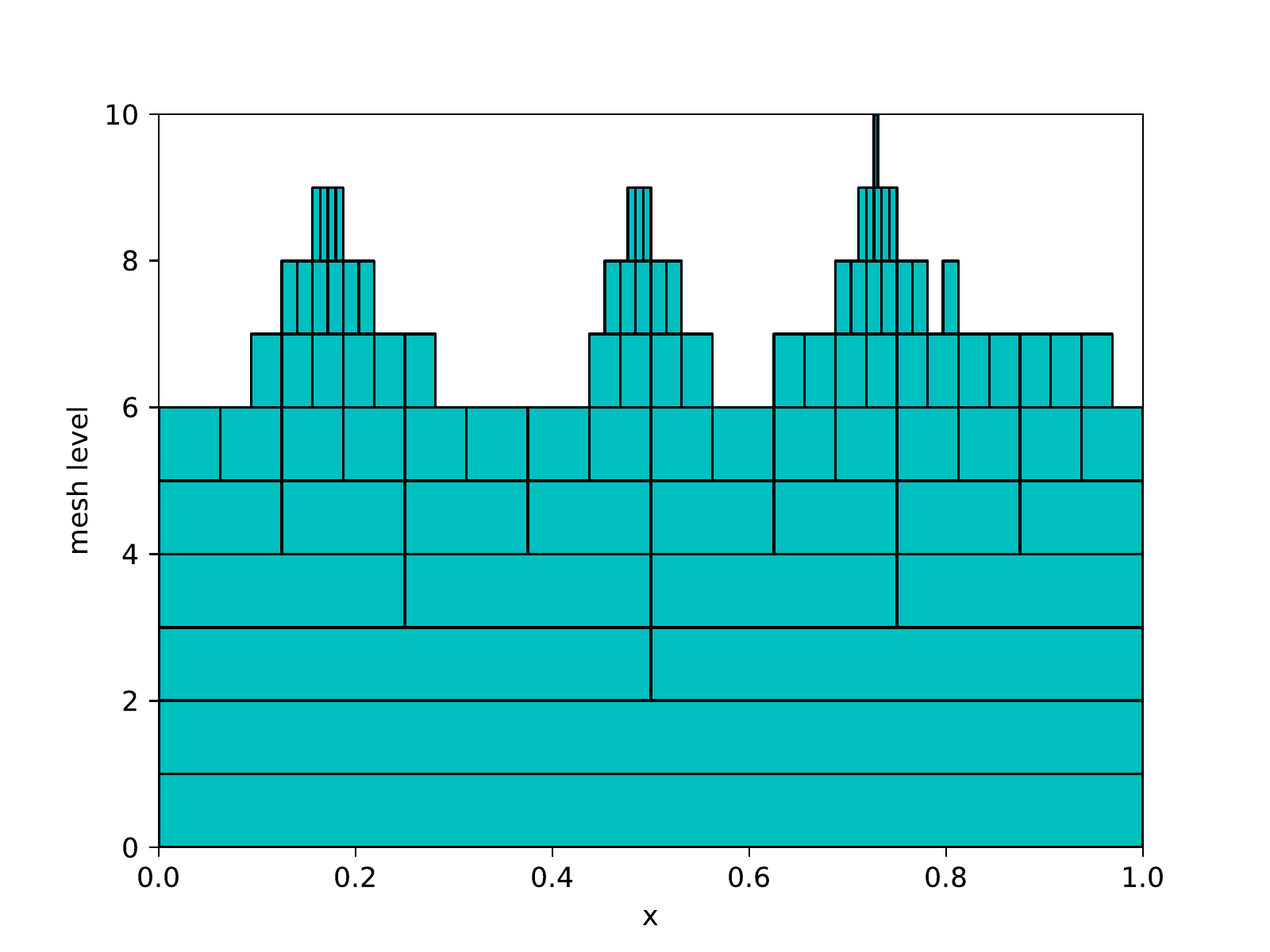}
		\end{minipage}
	}    
	\caption{Example \ref{exam:2-coupled-interaction}: coupled NLS equation, triple  solitons. Left: numerical solutions; right: active elements. $t=0$, 20 and 50. $N=9, k=3$, $\epsilon=10^{-4}$ and $\eta=4 \times 10^{-5}$.}
	\label{fig:2-coupled-triple}
\end{figure}

\end{exam}

\subsection{NLS equation in 2D}

\begin{exam}\label{exam:2D_singular}
In this example, we consider the singular solutions for the 2D NLS equation 
\begin{equation}
\label{2D_NLS}
i u_t + \frac{1}{M^2}  u_{xx} + \frac{1}{M^2}  u_{yy} +    \abs{u}^2 u = 0,
\end{equation}
with initial condition \cite{xu2005schrodinger}
\begin{equation}
u(x,0) = (1 + \sin X) (2 + \sin Y)
\end{equation}
where $X= Mx, Y=My, M=2 \pi$. Periodic boundary conditions are applied in $[0,1]^2$.
Strong evidence of a singularity in finite time is obtained. 
%We take $N=7, k=3$, $\epsilon=10^{-4}$ and $\eta= 10^{-5}$.
The plots of $|u|$  and active elements at $t=0$ and $t=0.108$ are shown in Figure \ref{fig:2D_singular}. %The contour plots are with 17 equally-spaced contour lines. 
From the results, we can observe that a singular is generated at $t=0.108$ and our method can capture the structure adaptively.

 \begin{figure}
 	\centering
 	\subfigure[surface of $|u|$ at $t=0$]{
 		\begin{minipage}[b]{0.46\textwidth}
 			\includegraphics[width=1\textwidth]{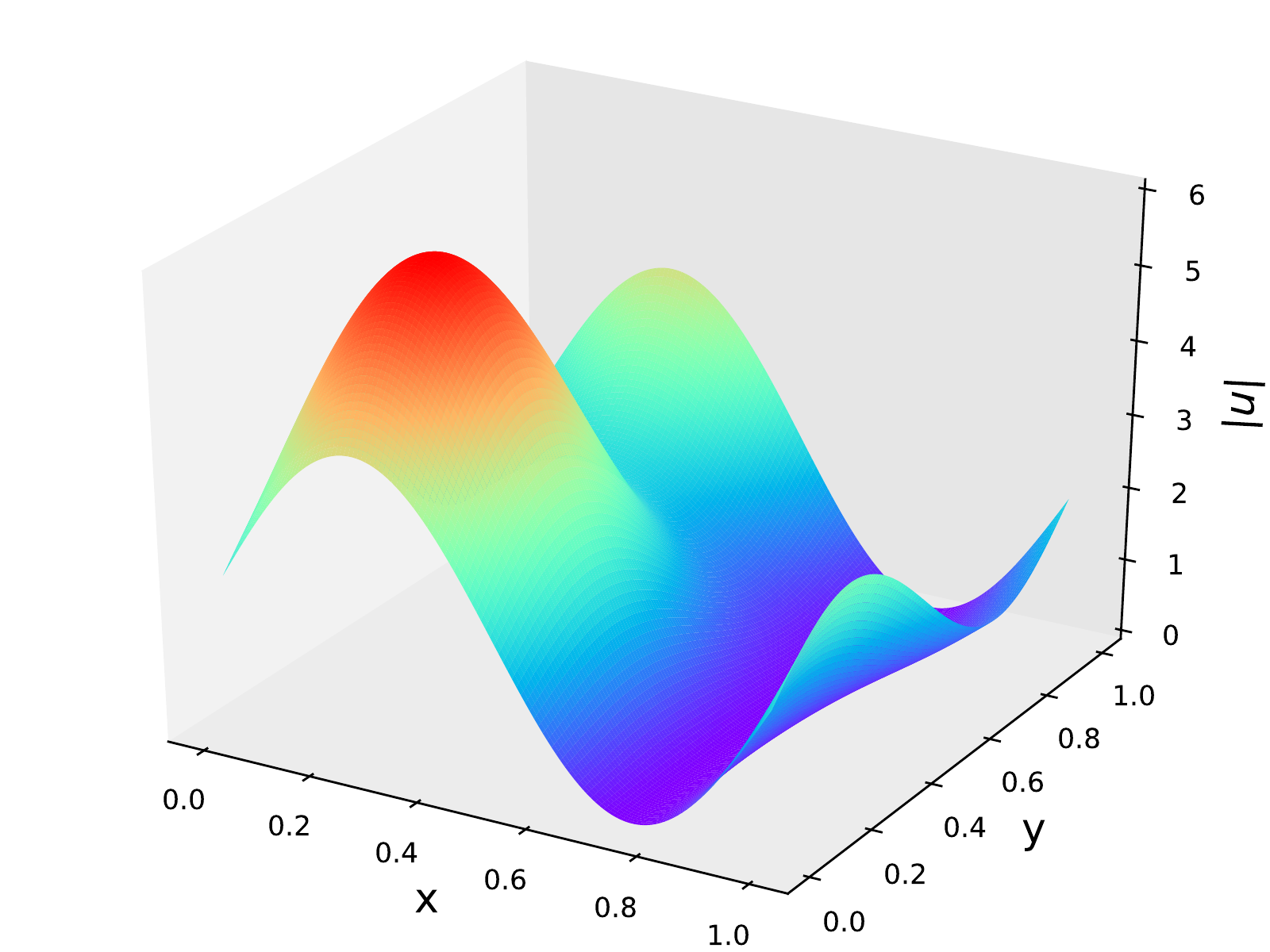}
 		\end{minipage}
 	}
 	\subfigure[active elements at $t=0$]{
 		\begin{minipage}[b]{0.46\textwidth}
 			\includegraphics[width=1\textwidth]{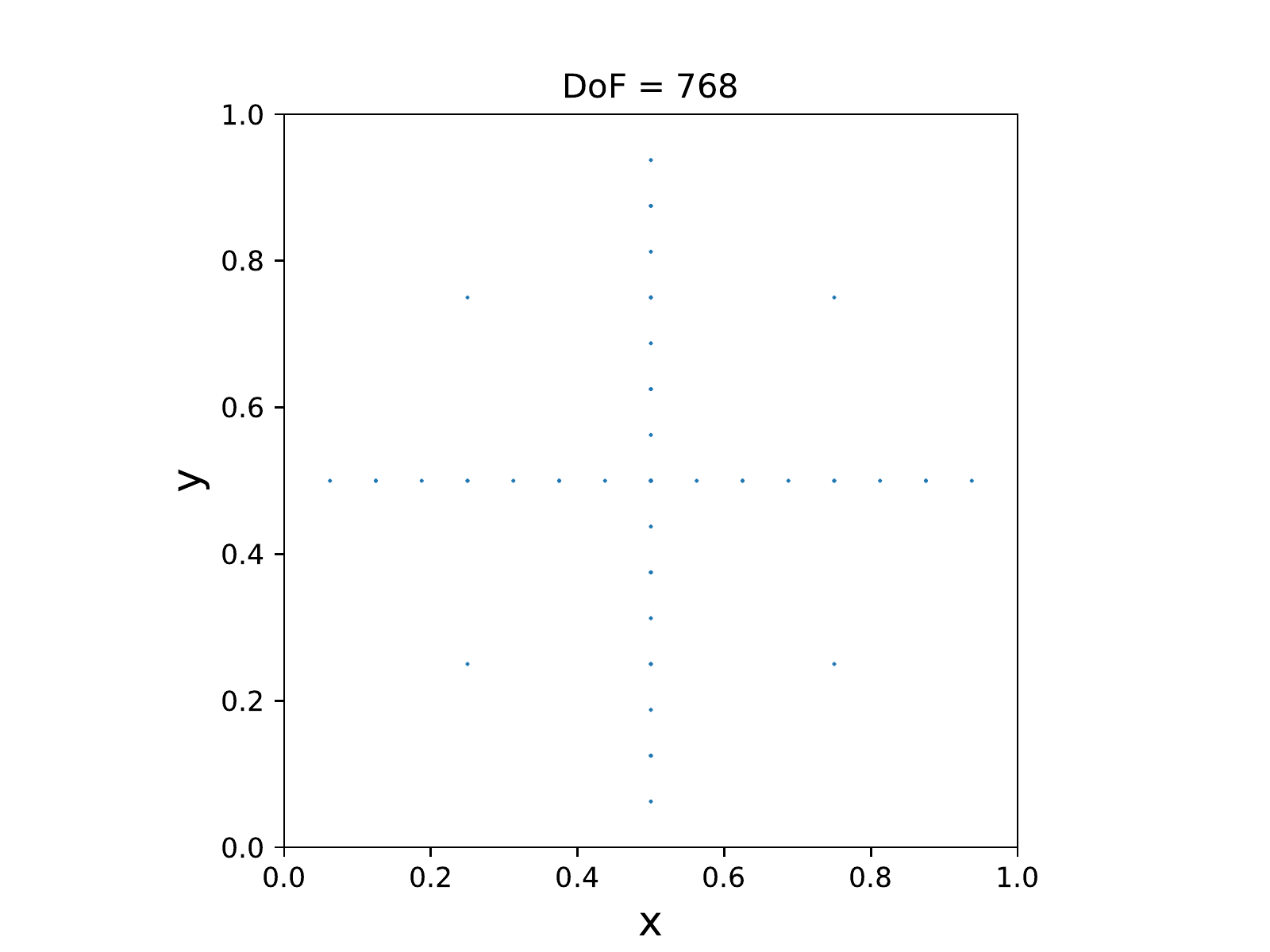}
 		\end{minipage}
 	}
 	\bigskip
% 	\subfigure[contour of $|u|$ at $t=0$]{
% 		\begin{minipage}[b]{0.46\textwidth}
% 			\includegraphics[width=1\textwidth]{fig/2D_singular/contour_init.eps}
% 		\end{minipage}
% 	}
% 	\subfigure[contour of $|u|$ at $t=0.108$]{
% 		\begin{minipage}[b]{0.46\textwidth}    
% 			\includegraphics[width=1\textwidth]{fig/2D_singular/contour_final.eps}
% 		\end{minipage}
% 	}
% 	\bigskip
 	 \subfigure[surface of $|u|$ at $t=0.108$]{
 		\begin{minipage}[b]{0.46\textwidth}    
 			\includegraphics[width=1\textwidth]{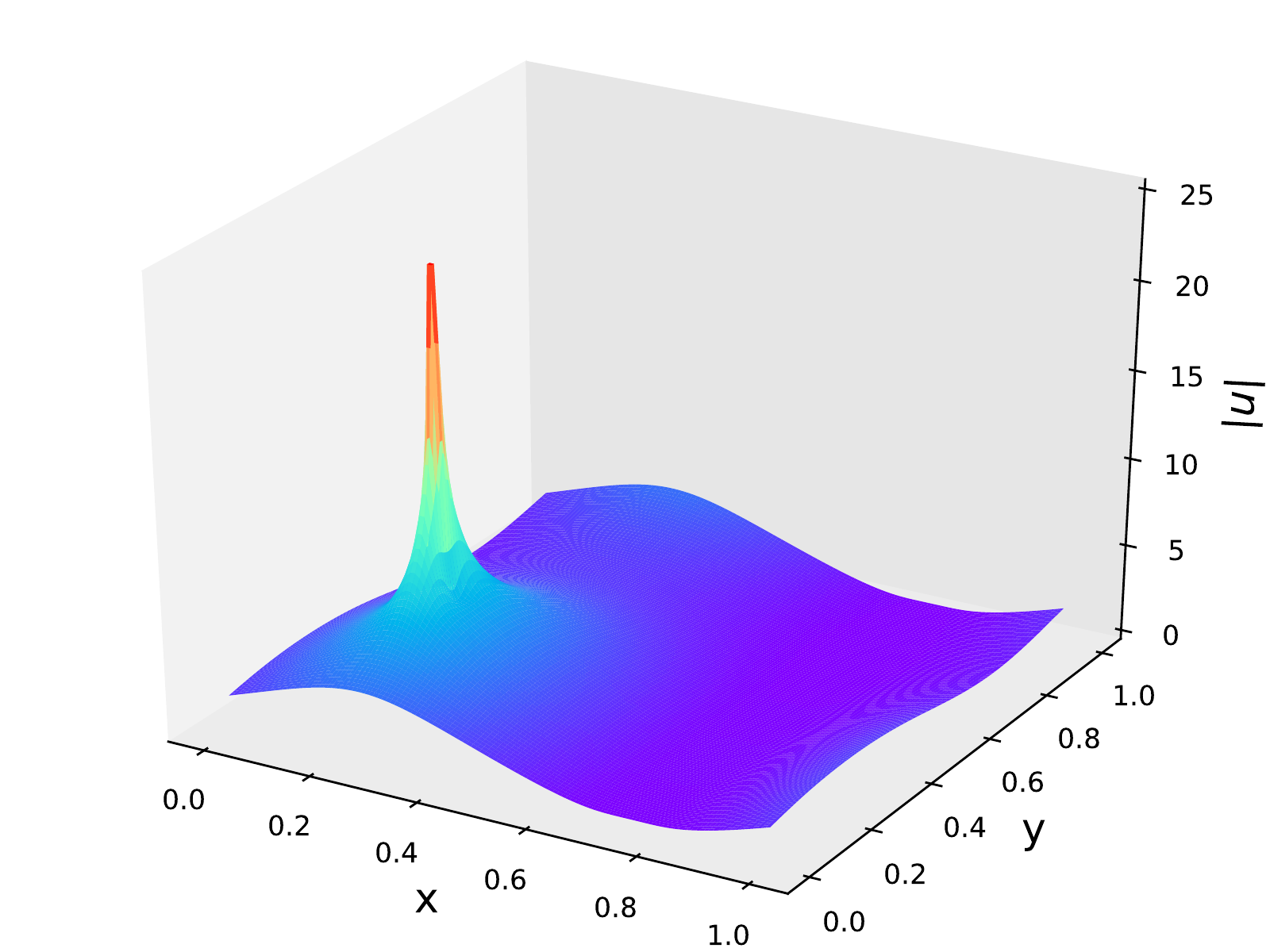}
 		\end{minipage}
 	}
 	\subfigure[active elements at $t=0.108$]{
 		\begin{minipage}[b]{0.46\textwidth}    
 			\includegraphics[width=1\textwidth]{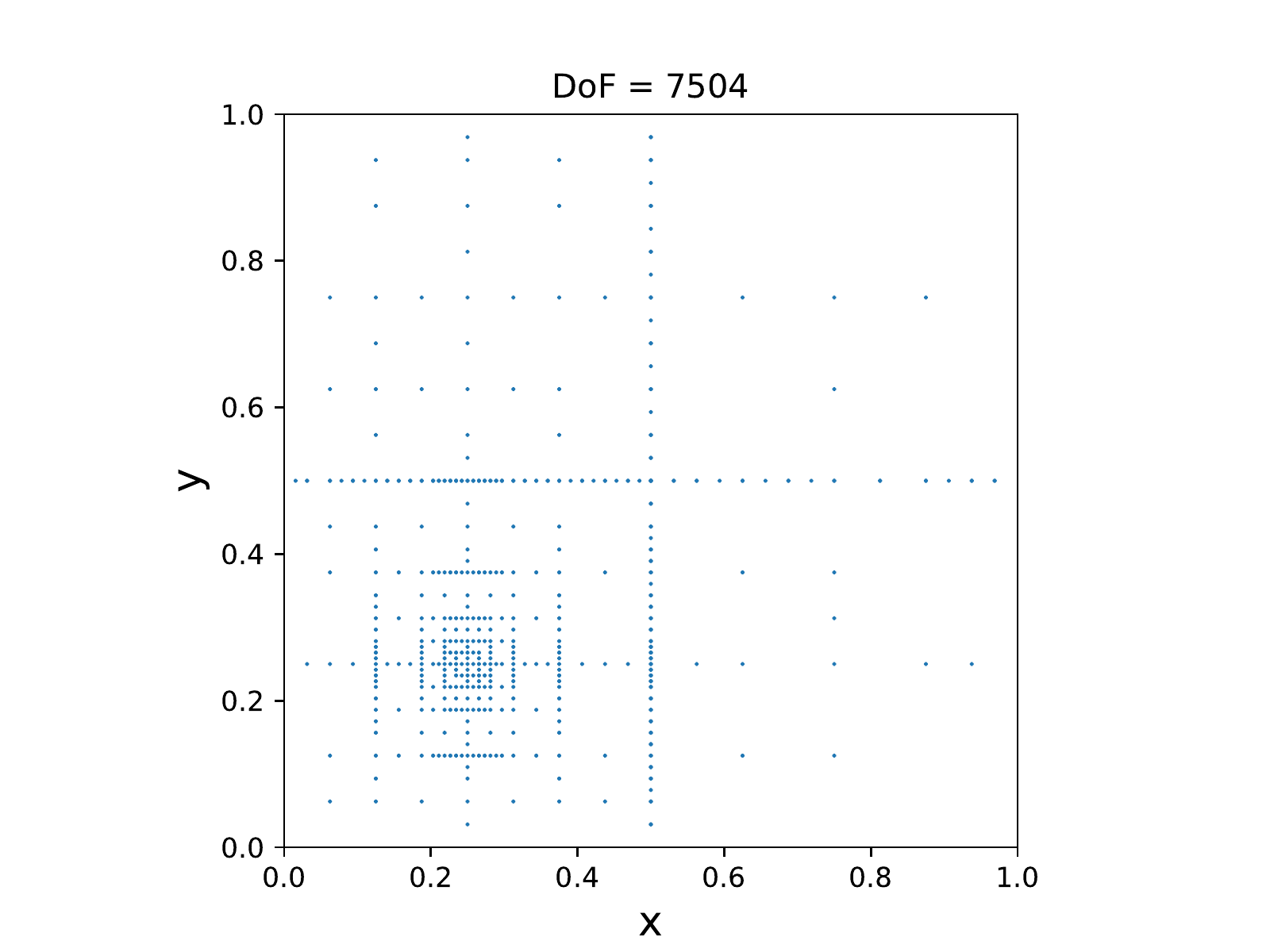}
 		\end{minipage}
 	}
 	\caption{Example \ref{exam:2D_singular}: singular solutions in 2D NLS equation. Left: numerical solutions; right: active elements. $t=0$ and 0.108.
 	$N=7, k=3$, $\epsilon=10^{-4}$ and $\eta=10^{-5}$.}    
% 	\caption{Example \ref{exam:2D_singular}: singular solutions in 2D NLS equation. Left: $t=0$; right: $t=0.108$. Top: Surface of $|u|$;  middle: contour of  $|u|$; bottom: active elements.
% 	$N=7, k=3$, $\epsilon=10^{-4}$ and $\eta=10^{-5}$.}
 	\label{fig:2D_singular}
 \end{figure}
 
\end{exam}

\begin{exam}\label{exam:2D_blow_up_1}
	
	In this example, we consider the  2D NLS equation \eqref{2D_NLS} with initial condition \cite{zhang2012compact}
	\begin{equation}
	u(x,0) = 2.0 + 0.01\sin(X+\frac{\pi }{4})\sin(Y+\frac{\pi }{4})
	\end{equation}
	with $X= M(x-0.5)$, $Y=M(y-0.5)$, $M=2 \pi$. Periodic boundary conditions are applied in $[0,1]^2$.
%	We take $N=7, k=3$, $\epsilon=10^{-4}$ and $\eta= 10^{-5}$.
	The plots of $|u|$  and active elements at $t=0$ and $t=1.5813$ are shown in Figure \ref{fig:2D_blow_1}. We can observe the blow-up phenomenon in $|u|$ at $t=1.5813$.

	\begin{figure}
		\centering
		\subfigure[surface of $|u|$ at $t=0$]{
			\begin{minipage}[b]{0.46\textwidth}
				\includegraphics[width=1\textwidth]{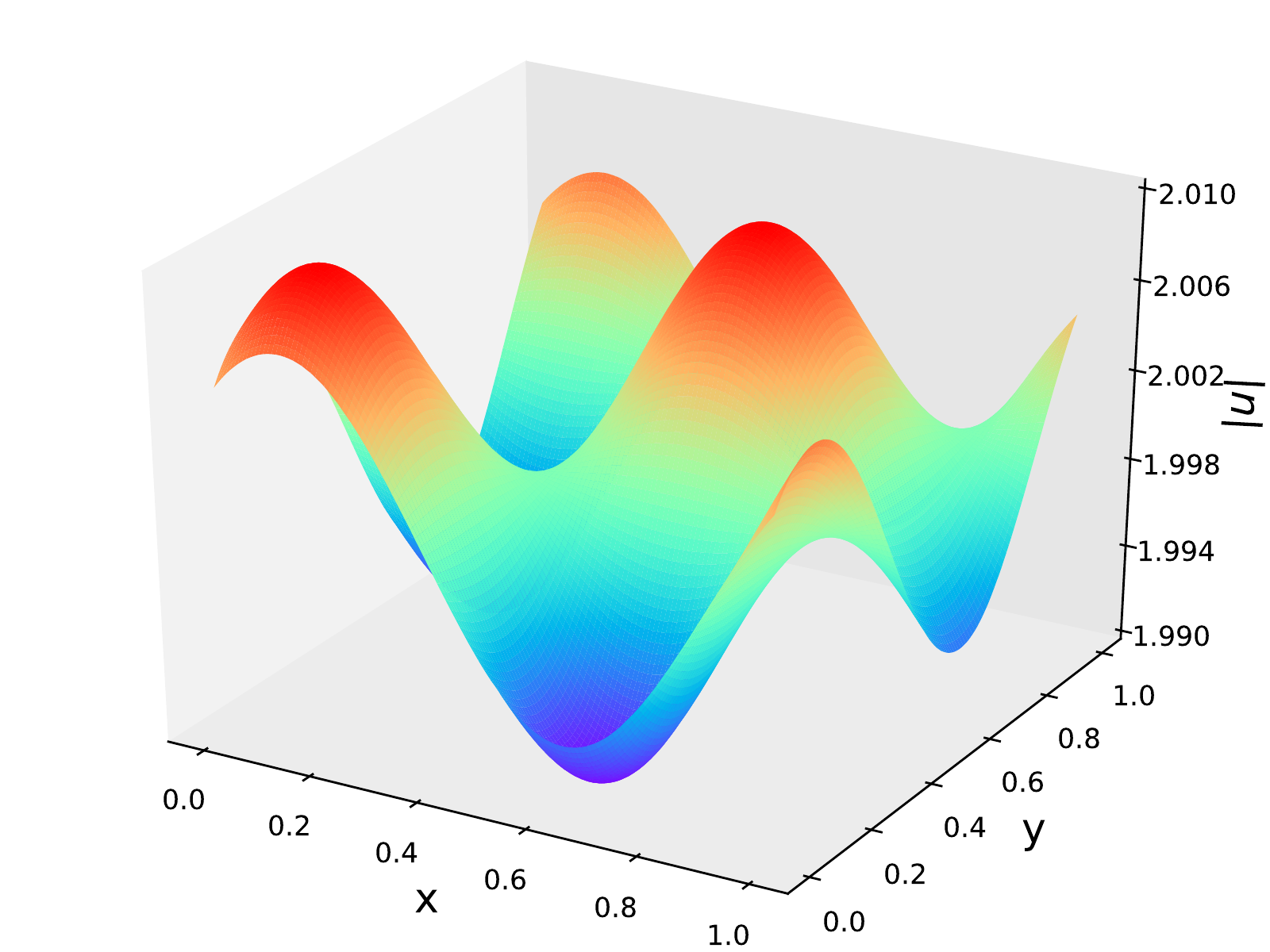}
			\end{minipage}
		}
		\subfigure[active elements at $t=0$]{
			\begin{minipage}[b]{0.46\textwidth}
				\includegraphics[width=1\textwidth]{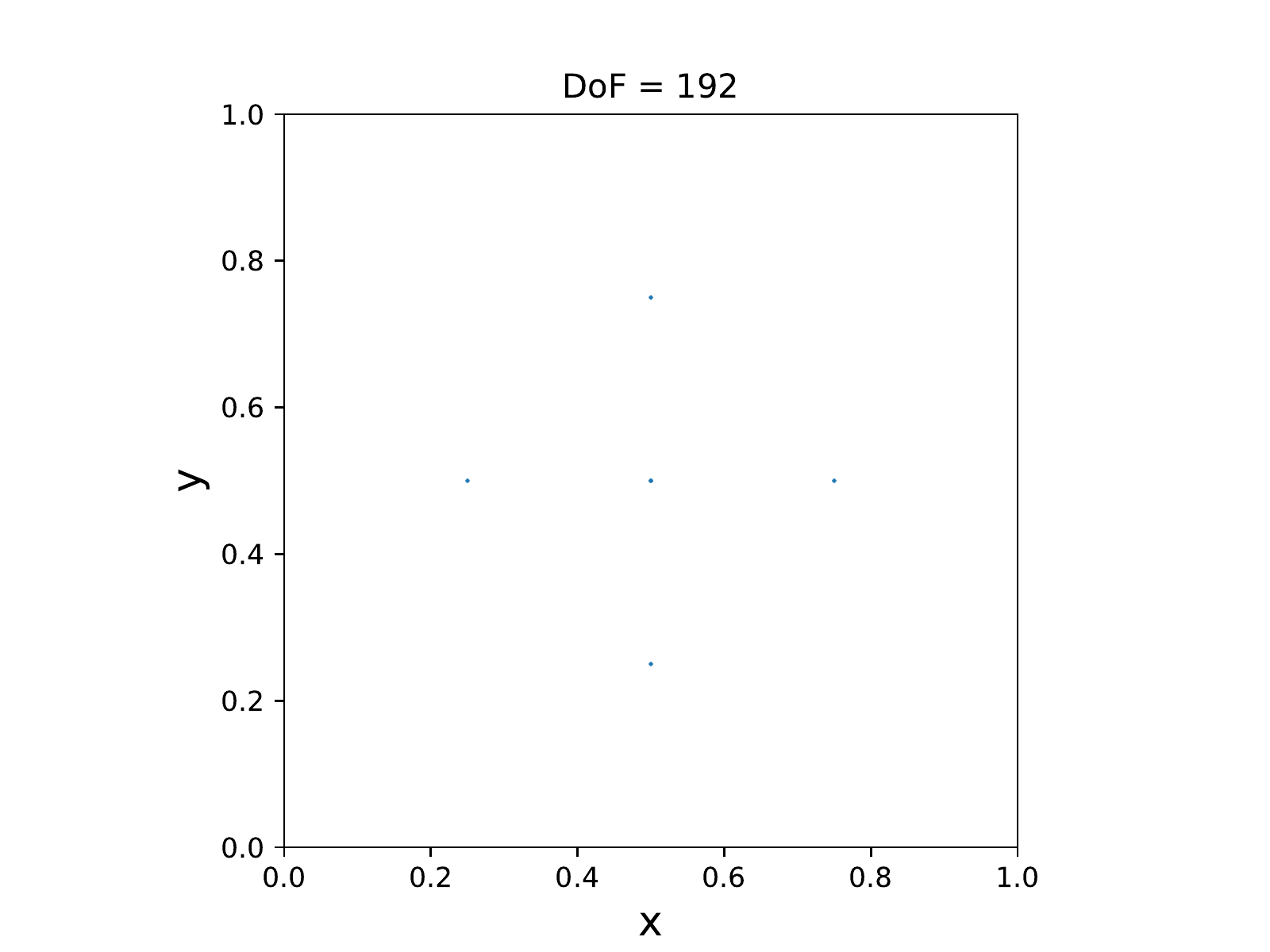}
			\end{minipage}
		}
		\bigskip
%		\subfigure[contour of $|u|$ at $t=0$]{
%			\begin{minipage}[b]{0.46\textwidth}
%				\includegraphics[width=1\textwidth]{fig/2D_blow/contour_init.eps}
%			\end{minipage}
%		}
%		\subfigure[contour of $|u|$ at $t=1.5813$]{
%			\begin{minipage}[b]{0.46\textwidth}    
%				\includegraphics[width=1\textwidth]{fig/2D_blow/contour_final.eps}
%			\end{minipage}
%		}
%		\bigskip
		\subfigure[surface of $|u|$ at $t=1.5813$]{
			\begin{minipage}[b]{0.46\textwidth}    
				\includegraphics[width=1\textwidth]{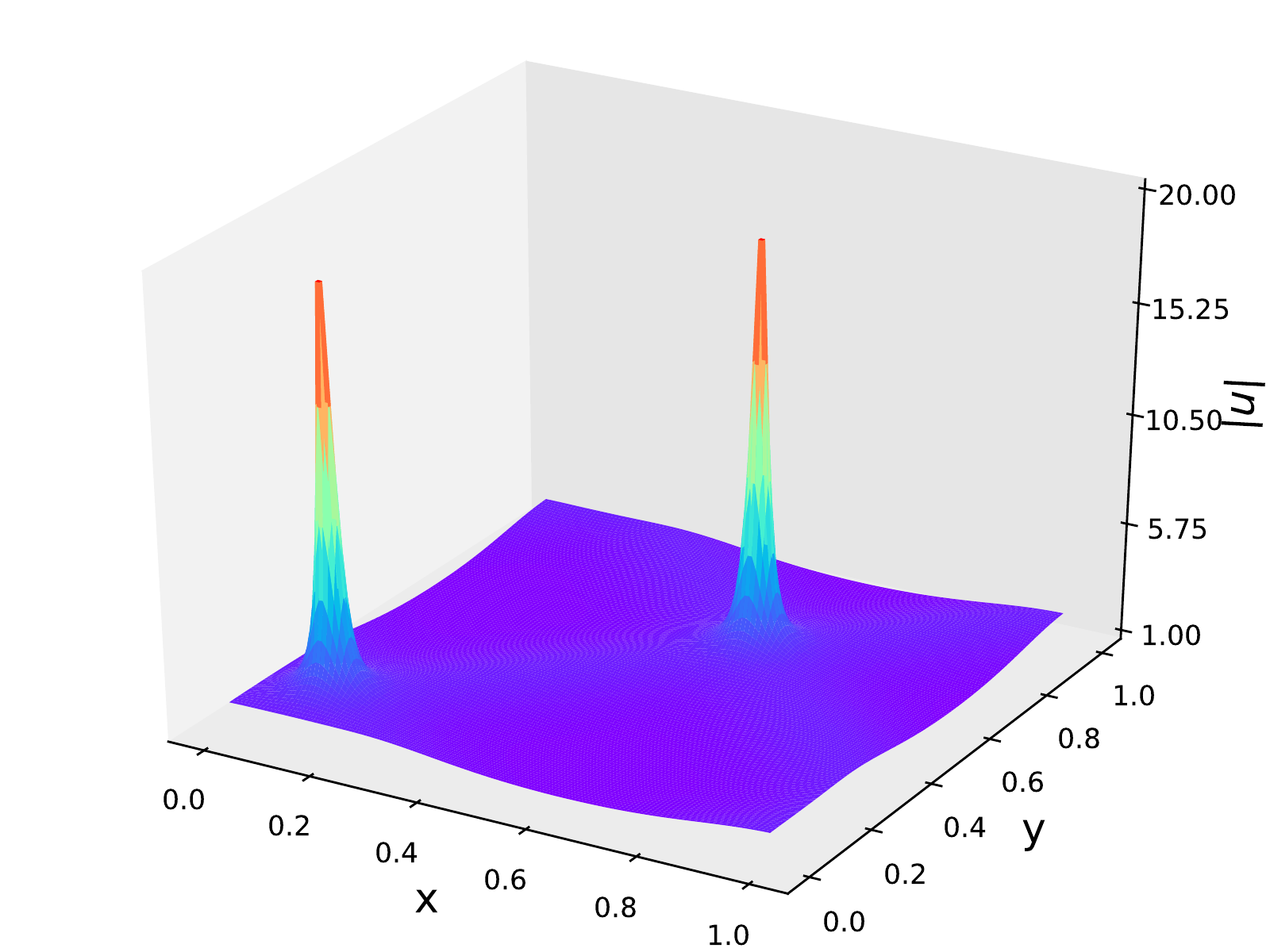}
			\end{minipage}
		}
		\subfigure[active elements at $t=1.5813$]{
			\begin{minipage}[b]{0.46\textwidth}    
				\includegraphics[width=1\textwidth]{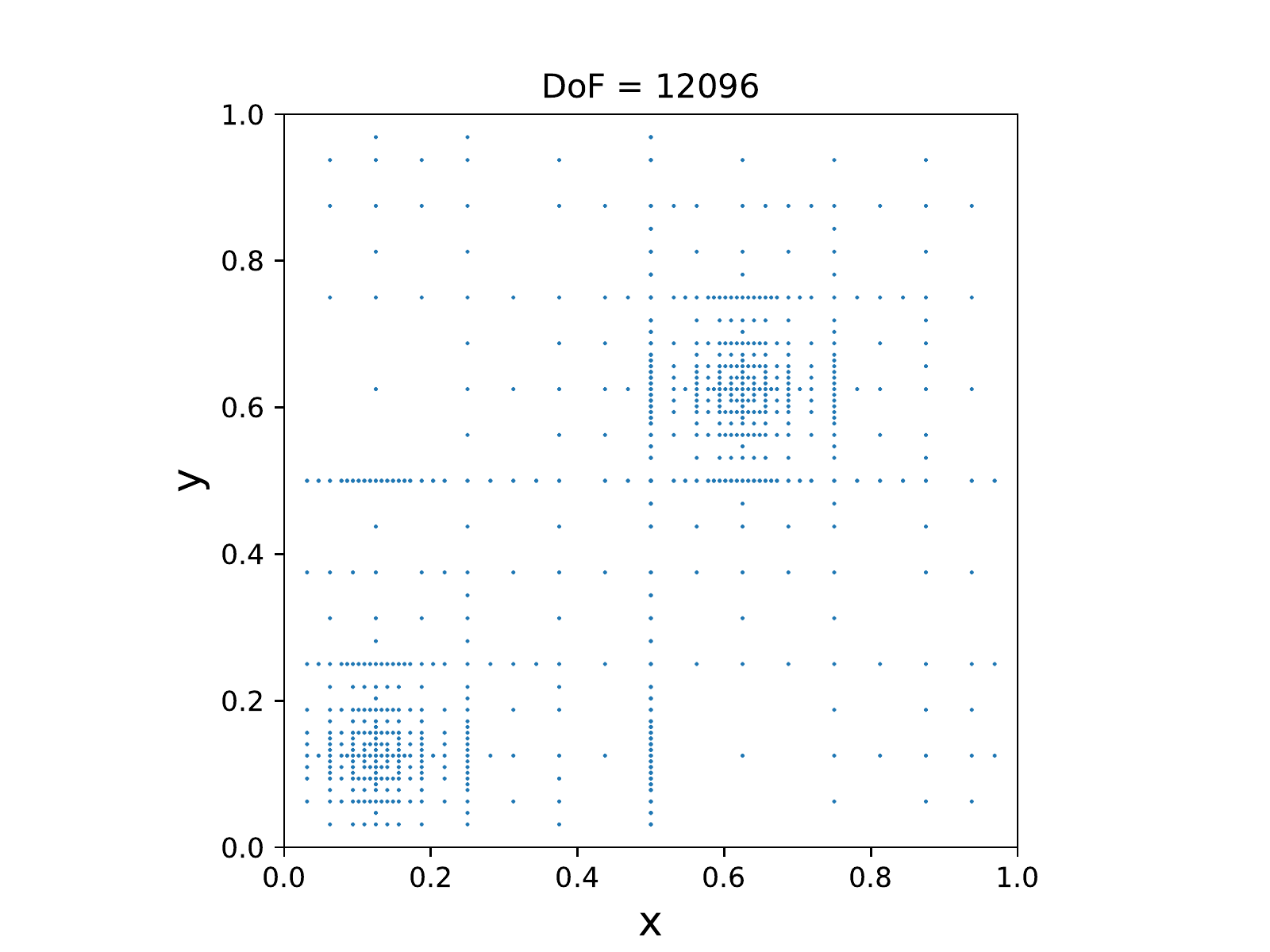}
			\end{minipage}
		} 
		\caption{Example \ref{exam:2D_blow_up_1}: blow up solution in 2D NLS equation. Left: numerical solutions; right: active elements. $t=0$ and 1.5813.  
			$N=7, k=3$, $\epsilon=10^{-4}$ and $\eta=10^{-5}$.}   
%		\caption{Example \ref{exam:2D_blow_up_1} for 2D NLS equation. Left: $t=0$; right: $t=1.5813$. Top: Surface of $|u|$;  middle: contour of  $|u|$; bottom: active elements.
%			$N=7, k=3$, $\epsilon=10^{-4}$ and $\eta=10^{-5}$.}
		\label{fig:2D_blow_1}
	\end{figure}

\end{exam}

\begin{exam}\label{exam:2D_blow_up_2}
	
	In this example, we consider the 2D NLS equation \eqref{2D_NLS} with initial condition \cite{zhang2017conservative}
	\begin{equation}
	u(x,0) = 6\sqrt{2}\exp(-X^2-Y^2)
	\end{equation}
	with $X= M(x-0.5)$, $Y=M(y-0.5)$ and $M=10$. Periodic boundary conditions are applied in $[0,1]^2$.
	The plots of $|u|$  and active elements at $t=0$ and $t=0.04$ are shown in Figure \ref{fig:2D_blow_2}. We observe from the results that the solution blows up in the center and our method can capture the blow up phenomenon.

	\begin{figure}
		\centering
		\subfigure[surface of $|u|$ at $t=0$]{
			\begin{minipage}[b]{0.46\textwidth}
				\includegraphics[width=1\textwidth]{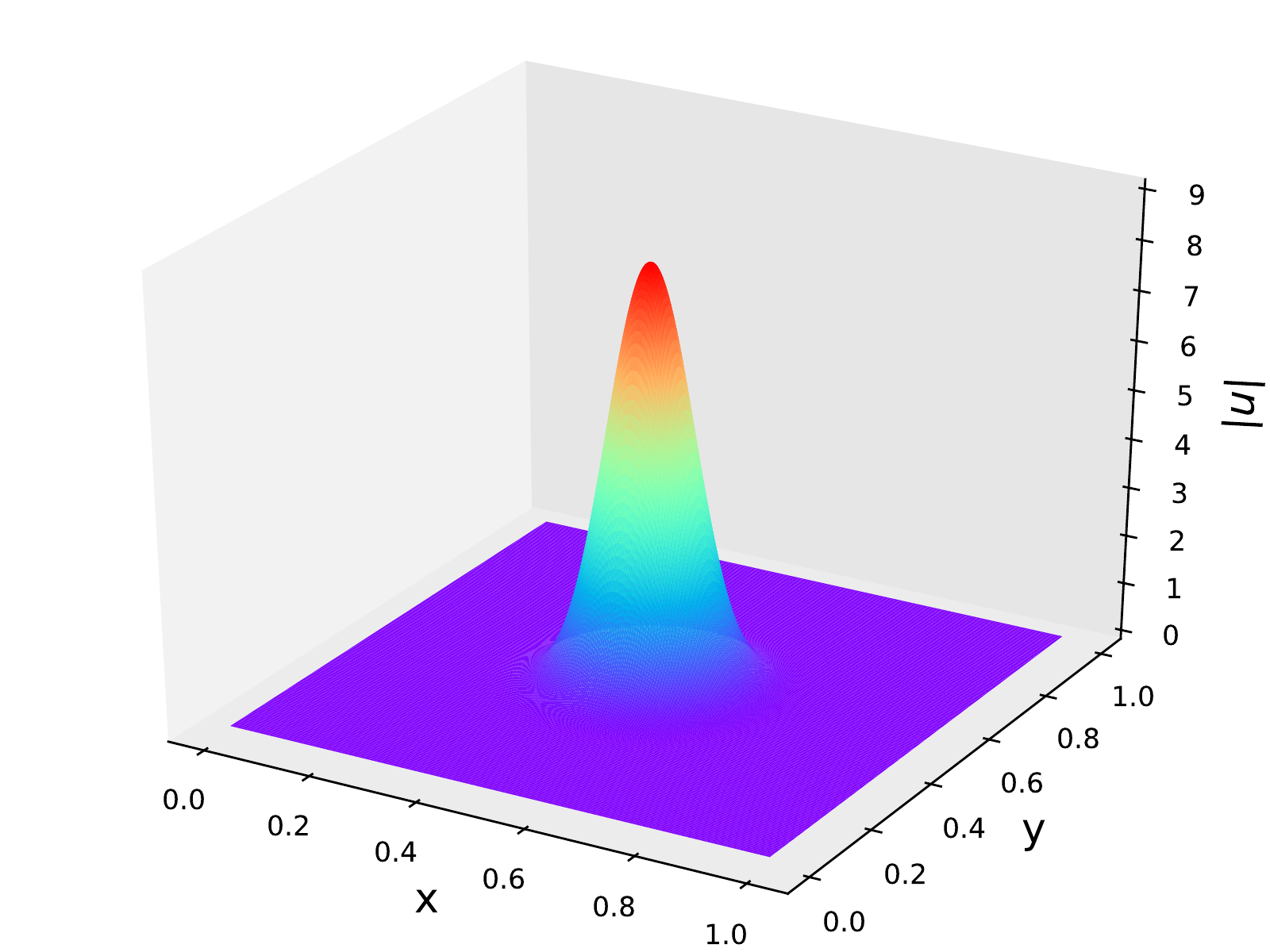}
			\end{minipage}
		}
		\subfigure[active elements at $t=0$]{
			\begin{minipage}[b]{0.46\textwidth}
				\includegraphics[width=1\textwidth]{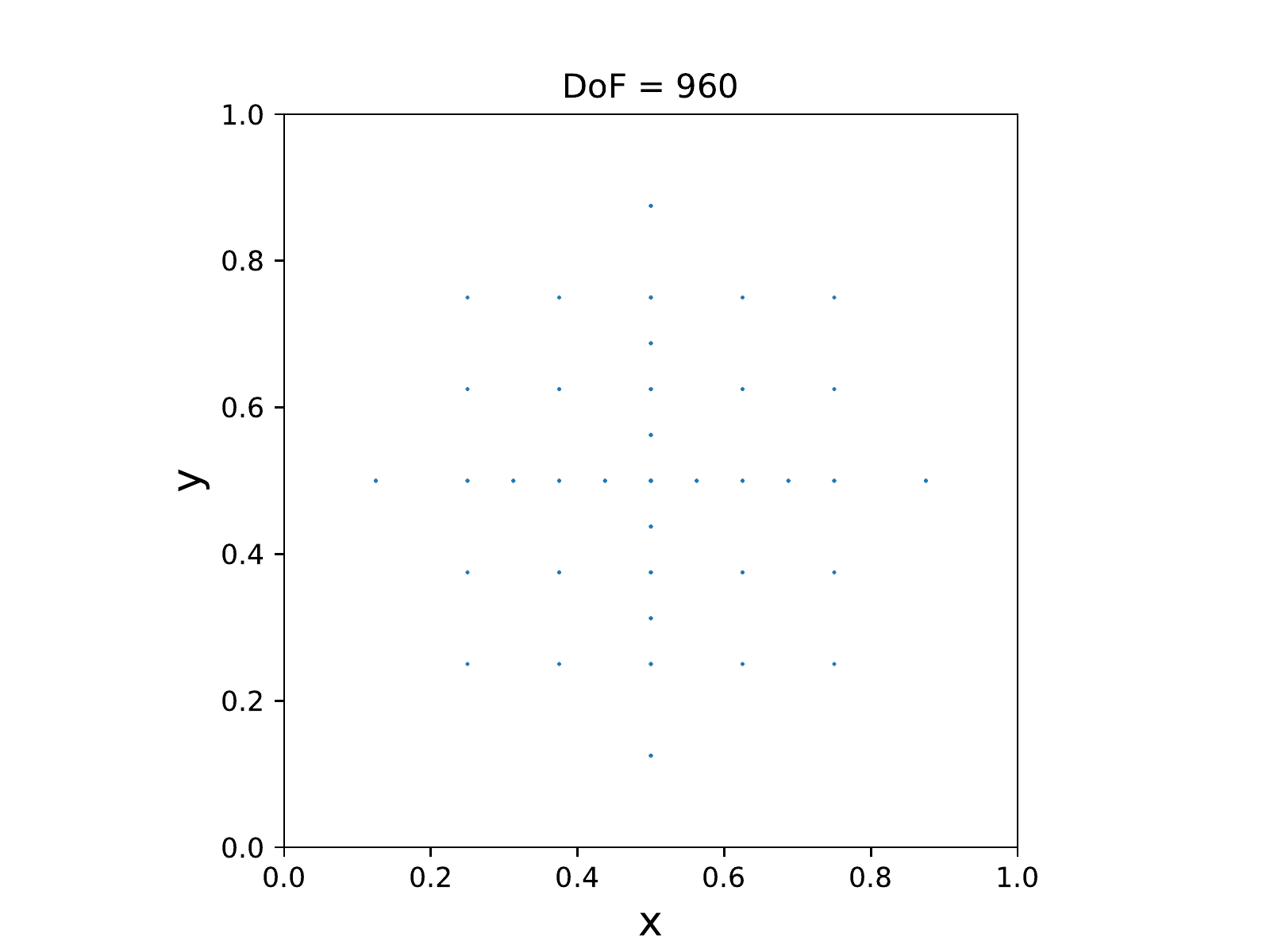}
			\end{minipage}
		}
		\bigskip
%		\subfigure[contour of $|u|$ at $t=0$]{
%			\begin{minipage}[b]{0.46\textwidth}
%				\includegraphics[width=1\textwidth]{fig/2D_blow_1/contour_init.eps}
%			\end{minipage}
%		}
%		\subfigure[contour of $|u|$ at $t=0.04$]{
%			\begin{minipage}[b]{0.46\textwidth}    
%				\includegraphics[width=1\textwidth]{fig/2D_blow_1/contour_t04.eps}
%			\end{minipage}
%		}
%		\bigskip
		\subfigure[surface of $|u|$ at $t=0.04$]{
			\begin{minipage}[b]{0.46\textwidth}    
				\includegraphics[width=1\textwidth]{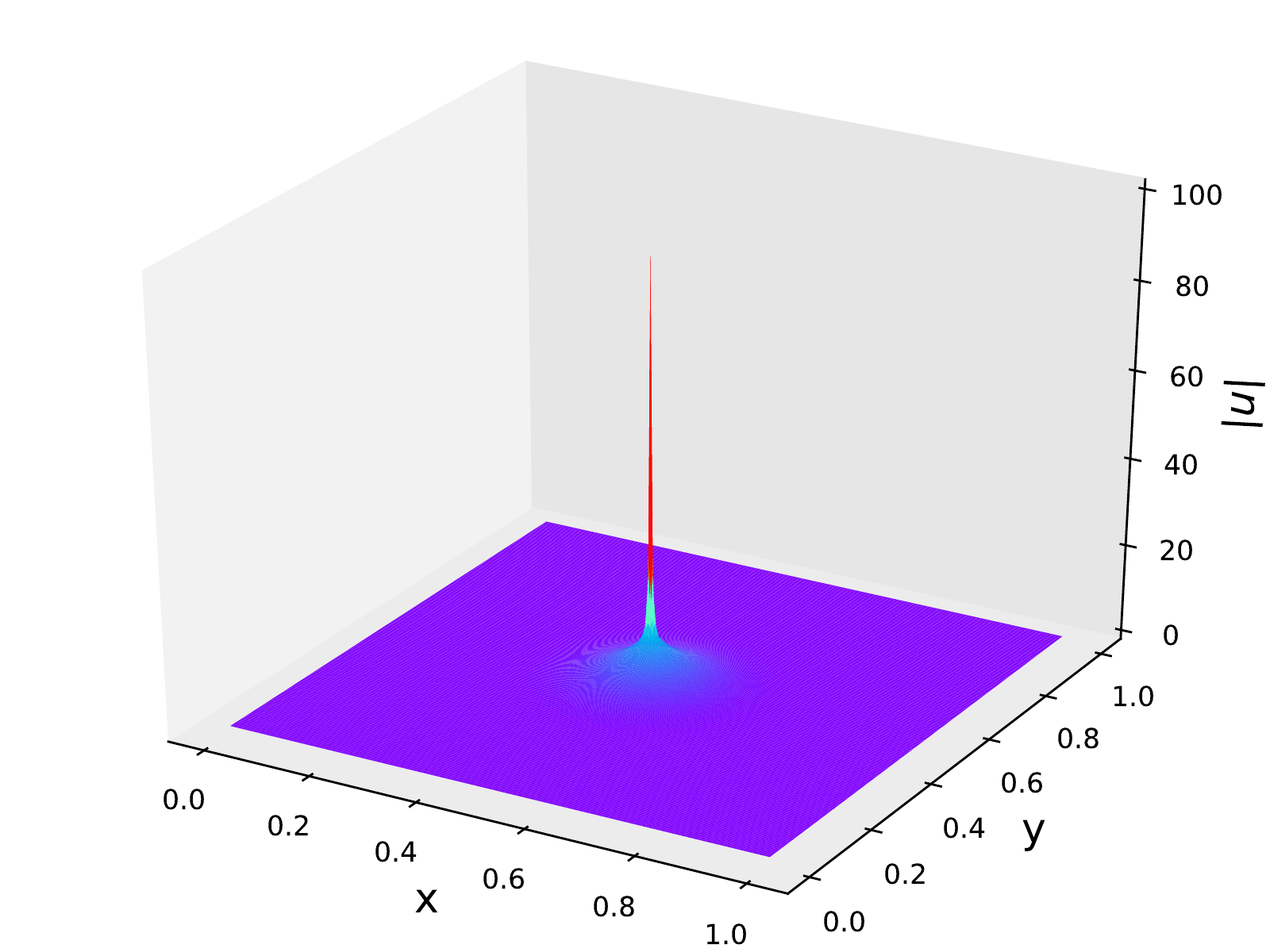}
			\end{minipage}
		}
		\subfigure[active elements at $t=0.04$]{
			\begin{minipage}[b]{0.46\textwidth}    
				\includegraphics[width=1\textwidth]{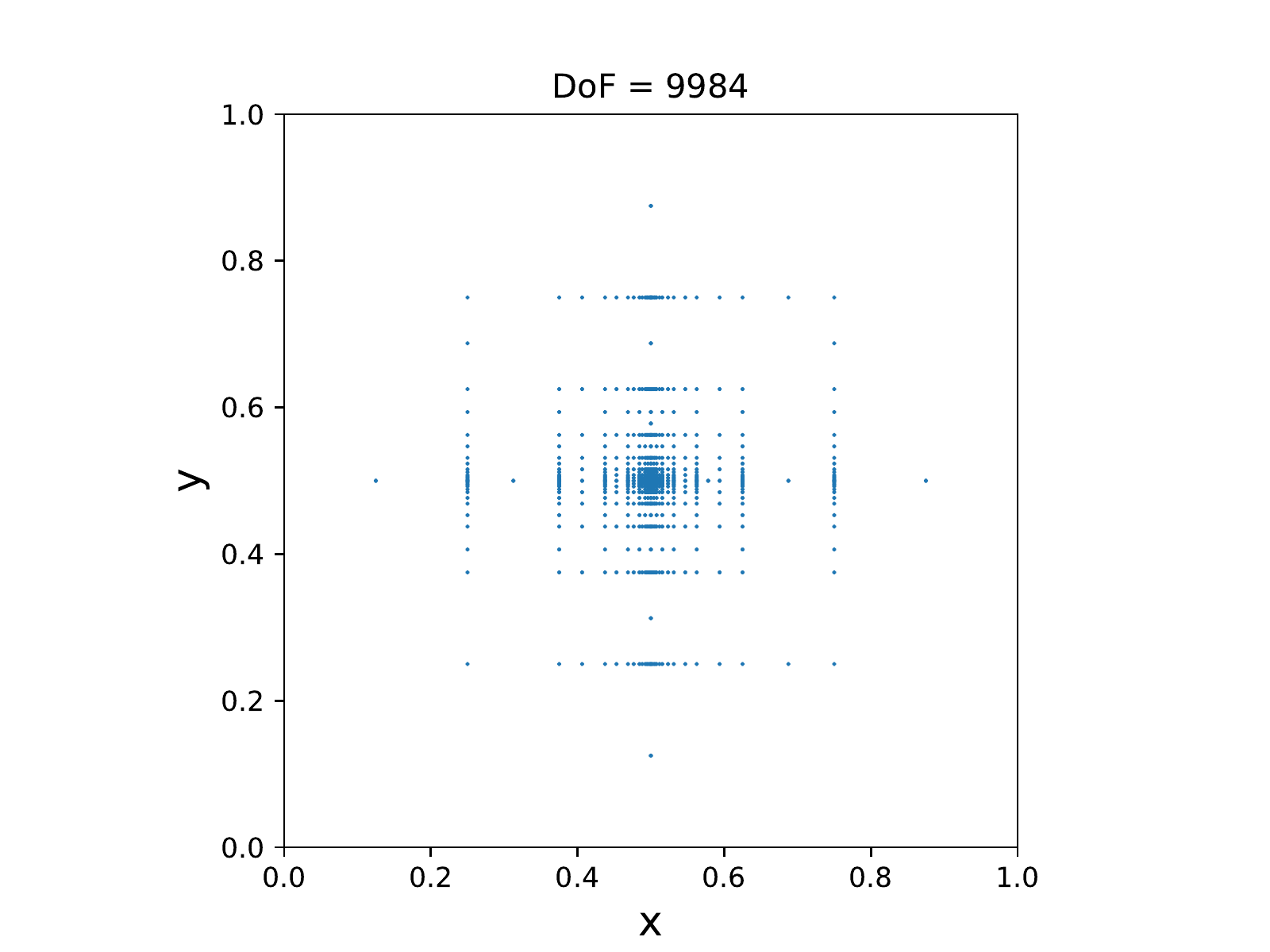}
			\end{minipage}
		} 
		\caption{Example \ref{exam:2D_blow_up_2}:  single blow up solution in 2D NLS equation. Left: numerical solutions; right: active elements. $t=0$ and 0.04.
		$N=9, k=3$, $\epsilon=10^{-3}$ and $\eta=10^{-4}$.}   
%		\caption{Example \ref{exam:2D_blow_up_2} for 2D NLS equation. Left: $t=0$; right: $t=0.04$. Top: Surface of $|u|$;  middle: contour of  $|u|$; bottom: active elements.
%			$N=9, k=3$, $\epsilon=10^{-3}$ and $\eta=10^{-4}$.}
		\label{fig:2D_blow_2}
	\end{figure}

\end{exam}

% section 6
\section{Conclusion}

In this paper, we propose an adaptive multiresolution ultra-weak DG method to solve nonlinear Schr\"{o}dinger equations. The adaptive multiwavelets are applied to achieve the multiresolution. The Alpert’s multiwavelets are used to express the DG solution and the interpolatory multiwavelets are exploited to compute the nonlinear source term.
Various numerical experiments are presented to demonstrate the excellent capability of capturing the soliton waves and the blow-up phenomenon. The code generating the results in this paper can be found at the GitHub link: \url{https://github.com/JuntaoHuang/adaptive-multiresolution-DG}.

\section*{Conflict of interest}
On behalf of all authors, the corresponding author states that there is no conflict of interest.

% section 7
\section*{Acknowledgment}

We would like to thank Qi Tang and Kai Huang for the assistance and discussion in code implementation.

\bibliographystyle{abbrv}
\bibliography{ref_schrodinger}

\end{document}